\documentclass[11pt]{amsart}

\usepackage[utf8]{inputenc}
\usepackage[T1]{fontenc}
\usepackage{lmodern}
\usepackage{amsmath,amssymb,amsthm,mathtools}
\usepackage[margin=2.8cm]{geometry}
\usepackage{microtype}
\usepackage{enumitem}
\usepackage{array}
\usepackage{hyperref}

\DeclareMathOperator{\Id}{Id}
\DeclareMathOperator{\ord}{ord}
\DeclareMathOperator{\supp}{supp}
\DeclareMathOperator{\Aut}{Aut}
\newcommand{\cE}{\mathcal{E}}
\newcommand{\EE}{\mathrm{E}}
\newcommand{\QQ}{\mathbb{Q}}
\newcommand{\FF}{\mathbb{F}}
\newcommand{\Aff}{\mathrm{Aff}}
\newcommand{\Qout}{Q_{\mathrm{out}}}

\theoremstyle{plain}
\newtheorem{theorem}{Theorem}[section]
\newtheorem{proposition}[theorem]{Proposition}
\newtheorem{lemma}[theorem]{Lemma}
\newtheorem{corollary}[theorem]{Corollary}
\theoremstyle{definition}
\newtheorem{definition}[theorem]{Definition}
\newtheorem{example}[theorem]{Example}
\newtheorem{question}[theorem]{Question}
\theoremstyle{remark}
\newtheorem{remark}[theorem]{Remark}

\begin{document}

\title[Elementary length]{Elementary length of polynomial automorphisms in dimension
at least three}

\author{El\.zbieta Adamus}
\address{Faculty of Applied Mathematics, AGH University of Krakow,\\
al. Mickiewicza 30, 30-059 Krak\'ow, Poland}
\email{esowa@agh.edu.pl}
\urladdr{https://orcid.org/0000-0003-3533-6279}

\author{Zbigniew Hajto}
\address{Faculty of Mathematics and Computer Science, Jagiellonian University,
ul. \L ojasiewicza 6, 30-348 Krak\'ow, Poland}
\email{zbigniew.hajto@uj.edu.pl}
\urladdr{https://orcid.org/0000-0002-1103-4404}


\begin{abstract}
Let $\EE(K,n)$ be the group generated by the elementary automorphisms of $K[X]$,\linebreak[1] $X=(X_1,\dots,X_n)$, the
maps which add to one coordinate a polynomial in the other coordinates. It is the group
of tame automorphisms with Jacobian determinant $1$. We study the word length $\ell$ on
$\EE(K,n)$ with respect to the elementary maps, for $n\ge3$ and in comparison with the
plane, where the Jung--van der Kulk amalgam provides a normal form. The decompositions of
$F$ of length $\#\supp(F)$, the number of coordinates moved by $F$, are determined by
$F^{-1}$, and $\ell(F\times\Id_r)=\ell(F)$ whenever $\ell(F)\le\#\supp(F)+2$.
Vasyunin's quadratic automorphism of $K[X_1,\dots,X_5]$ has exactly six minimal decompositions, which
differ in the degrees of their factors, and its dynamical degree is $\sqrt2$. In
characteristic zero a map $X+H\in K[X_1,X_2,X_3]^3$ with $H(0)=0$ and nilpotent Jacobian matrix is
Pascal finite if and only if the components of $H$ are linearly dependent. Over every field
the Anick automorphism $A$ of $K[X_1,\dots,X_4]$ is not a product of fewer than $7$ elementary maps,
while $\ell(A\times\Id_1)=6$. The lower bound rests on Cohn's theorem on
$\mathrm{GL}_2(K[X_3,X_4])$, and the decomposition of length $6$ is a polynomial form of
the commutator identity which places Cohn's matrix in $\mathrm{E}_3(K[X_3,X_4])$.
\end{abstract}

\maketitle

\section{Introduction}\label{sec:intro}

Let $K$ be a field and $n\ge2$. For $1\le j\le n$ let $U_j$ be the group of automorphisms
$F\in\Aut_KK[X_1,\dots,X_n]^n$ which change only $X_j$, by adding a polynomial in the other
coordinates. It is isomorphic to the additive group of $K[X_i:i\ne j]$. We use the
notation of van den Essen \cite[Chapter~5]{vdE}. Let $\EE(K,n)$ be the group generated by
$U_1,\dots,U_n$, and $\mathrm{TA}_n(K)$ the tame subgroup of $\Aut_KK[X_1,\dots,X_n]^n$, generated
by $\EE(K,n)$ and the affine group $\Aff(K,n)$. Then $\EE(K,n)$ is the group of tame
automorphisms with Jacobian determinant $1$ (Proposition~\ref{prop:Gamma}). It is a polynomial analogue of
the elementary subgroup $\mathrm{E}_n(R)$ of $\mathrm{SL}_n(R)$, which is generated by
its root subgroups, and the elements of $U_j\setminus\{\Id\}$, the \emph{elementary
maps}, play the role of the elementary matrices. Let $\cE(K,n)$ be the set of elementary
maps. This paper studies the word length on
$\EE(K,n)$ with respect to the elementary maps. For $F\in\EE(K,n)$ let $\ell(F)$
be the least number of elementary maps whose composition is $F$. We call it the
\emph{elementary length} of $F$, and a decomposition of $F$ into $\ell(F)$ elementary
maps a \emph{minimal decomposition}.

In the plane the word problem is governed by an amalgamated product. By the theorem of
Jung and van der Kulk, $\Aut_KK[X_1,X_2]^2=\Aff(K,2)\ast_{B(K,2)}J(K,2)$, where $J(K,2)$ is
the de~Jonqui\`eres group and $B(K,2)=\Aff(K,2)\cap J(K,2)$ \cite[Section~5.1]{vdE}, in the same
way as $\mathrm{SL}_2(K[t])$ is an amalgamated product by the theorem of Nagao
\cite[Ch.~II, \S1.6]{Ser}. The reduced word of an automorphism $F$ is determined by $F$
up to the insertion of elements of $B(K,2)$ between its letters. The canonical normal
form of \cite[Theorem~18]{AH_canon} removes this freedom by a choice of coset
representatives, and the number of de~Jonqui\`eres letters and the multidegree are
invariants of $F$ \cite[Proposition~23]{AH_canon}. The number of de~Jonqui\`eres letters
is the length of $F$ in the sense of Furter \cite{Fu}. It is half the displacement of a
vertex of the Bass--Serre tree of the amalgam, and it controls $\ell(F)$
(Theorem~\ref{thmC}).

For $n\ge3$ no such normal form is available. For $n=3$ and $\operatorname{char}K=0$ the
tame group is the generalized amalgamated product of three subgroups along their pairwise
intersections \cite{Wr}, a consequence of Umirbaev's defining relations \cite{Um}. Unlike
an amalgamated product of two groups, a generalized amalgamated product of three or more
groups need not admit a normal form for words \cite[Section~2]{Wr}, and for $n\ge4$ it is
not known whether the tame group has such a structure \cite{Wr}. For $n=3$ and
$\operatorname{char}K=0$ the tame group acts on a complete CAT(0) space \cite{LP}, and it
is acylindrically hyperbolic \cite{LP19}. For $n=3$ it is a proper subgroup of the automorphism group, since the
Nagata automorphism is wild if $\operatorname{char}K=0$ \cite{SU} or
$\operatorname{char}K\ge7$ \cite{Kur26}.

We ask how much of a decomposition of $F\in\EE(K,n)$ into elementary maps is
determined by $F$. Let $\supp(F)$ be the set of indices $i$ with $F_i\ne X_i$. For $r\ge0$ we write
$F\times\Id_r=(F_1,\dots,F_n,X_{n+1},\dots,X_{n+r})\in K[X_1,\dots,X_{n+r}]^{n+r}$ for the
automorphism obtained by the addition of $r$ variables. Every
decomposition of $F$ changes each coordinate of $\supp(F)$ at least once, so
$\ell(F)\ge\#\supp(F)$. The answer depends on the \emph{defect} $\ell(F)-\#\supp(F)$,
which counts the repetitions of directions in a minimal decomposition and the uses of
coordinates not moved by $F$ (Lemma~\ref{lem:defect}). The defect takes every value
$d\ge0$ already in the plane (Section~\ref{sec:plane}). We describe the decompositions
of defect $0$ completely, reduce those of defect $1$ to one unknown tail, and study three
maps from the literature, of defects $0$, $1$ and at least $5$, which show three
different phenomena. Table~\ref{tab:examples} in Section~\ref{sec:questions} collects
them.

\subsection*{Decompositions of defect zero}
The first result says that a shortest word in which no direction is repeated is
determined by the order in which the directions occur, and that its letters are read off
from $F^{-1}$.

\begin{theorem}[Decompositions of length $\#\supp(F)$]\label{thmA}
Let $F\in\EE(K,n)$, $F\ne\Id$, and $k=\#\supp(F)$.
\begin{enumerate}[label=\textup{(\roman*)},leftmargin=2em]
\item $\ell(F)\ge k$, and if $n\ge2$, then $F$ has decompositions of every length
$l\ge\ell(F)$.
\item The decompositions of $F$ of length $k$ are in bijection with the orderings of
$\supp(F)$ along which left elementary factors can be extracted successively
(Theorem~\ref{thm:min-decomp}). In
particular there are at most $k!$ of them, and the equality $\ell(F)=k$ reduces to at
most $k\cdot k!$ explicit polynomial tests, namely $k$ tests for each of the $k!$
orderings.
\item Let $r\ge1$. Every decomposition of $F\times\Id_r$ of length $k$ is of the
form $(Q_k\times\Id_r)\circ\dots\circ(Q_1\times\Id_r)$ for a decomposition
$Q_k\circ\dots\circ Q_1$ of $F$ of length $k$, and conversely. Moreover
$\ell(F\times\Id_r)\ge\min(\ell(F),k+2)$, so $\ell(F\times\Id_r)=\ell(F)$ if
$\ell(F)\le k+2$.
\end{enumerate}
\end{theorem}

Part (ii) is algorithmic. Given $F^{-1}$, each test asks whether a polynomial obtained
from the inverse of the remaining map involves a given variable, and the inverse of the
next remaining map is obtained by composition with one elementary map. In the language of
the groups $U_j$, the product map $U_{j_k}\times\dots\times U_{j_1}\to\EE(K,n)$ is
injective whenever $j_1,\dots,j_k$ are pairwise distinct (Corollary~\ref{cor:no-relation}).
This is the analogue of the uniqueness of the decomposition of a unipotent triangular
matrix into a product of elements of root subgroups taken in a fixed order. So every
relation between the groups $U_j$, other than the relations inside each $U_j$, repeats a
direction. Theorem~\ref{thmA} does not provide a normal form. It determines all
decompositions of defect $0$ in terms of $F^{-1}$, and the length is taken with respect
to the standard coordinates, so it is not invariant under linear conjugation
(Remark~\ref{rem:degrees}).

The example is Vasyunin's map $F=X+H$,
\[
H=\bigl(0,\;X_1X_3,\;X_1X_4+\tfrac12X_2^2,\;X_1X_5-X_2X_3,\;\tfrac12X_3^2\bigr).
\]
Its coefficients lie in $\mathbb{Z}[\tfrac12]$, so it is defined over every field $K$ of
characteristic $\ne2$, and we write $F_K$ when the field matters. It is a tame
automorphism by the decompositions of Section~\ref{ssec:three}, see also
\cite[Theorem~4.6.8]{Bt} for $\operatorname{char}K=0$. Its Jacobian matrix $J_H$ is
nilpotent and not strongly nilpotent, and in \cite{AH_sn} this map shows that strong
nilpotency cannot be replaced by nilpotency in the main theorem of that paper. If a map
$X+H$, with $H_i=0$ or $\ord H_i\ge2$ for every $i$, lies in a triangular subgroup, then
$\ell(F)=\#\supp(F)$ and a minimal decomposition is read off the dependency digraph of
$F$ (Lemma~\ref{lem:dag}). In dimension $n\le4$, $\operatorname{char}K=0$, every
quadratic homogeneous automorphism is linearly triangularizable, and $n=5$ is the least
dimension in which this can fail (Section~\ref{ssec:vasyunin-def}). Vasyunin's map is not
conjugate in $\Aut_KK[X_1,\dots,X_5]^5$ to any triangular map (Theorem~\ref{thmD}), and still
$\ell(F)=\#\supp(F)$. So its minimal decompositions are not obtained from a triangular
form, and Theorem~\ref{thmA} determines them.

\begin{theorem}[Minimal decompositions of Vasyunin's map]\label{thmB}
Let $F$ be Vasyunin's map over a field of characteristic $\ne2$.
\begin{enumerate}[label=\textup{(\roman*)},leftmargin=2em]
\item $\ell(F)=\#\supp(F)=4$, and $F$ has exactly six decompositions of length $4$. In
all of them the factor applied last is the same cubic map
$\Qout\colon X_3\mapsto X_3+X_1X_4+\tfrac12X_2^2-X_1^2X_5$, and the multiset of factor
degrees is $\{2,2,2,3\}$ for two of them and $\{2,2,3,3\}$ for the other four.
\item The least length of a decomposition of $F$ into elementary maps of degree at
most $2$ is $6$.
\item $F$ lies in no triangular subgroup attached to an order of the coordinates,
although $\ell(F)=\#\supp(F)$.
\end{enumerate}
\end{theorem}

So for $n\ge3$ the degrees of the factors of a minimal decomposition are not determined
by $F$, while in the plane the multidegree is an invariant.

\subsection*{Decompositions of defect one}
For decompositions of length $\#\supp(F)+1$ the same argument leaves one tail undetermined
(Proposition~\ref{prop:length-k1}). The maps of this kind which we study come from the
dependence problem for nilpotent Jacobian matrices in dimension three. He and Yan
\cite{HY} found a family of maps $X+H\in K[X_1,X_2,X_3]^3$ with nilpotent $J_H$, depending on a
polynomial $g$ of degree $d\ge1$, and by the classification of Wang and Yan
\cite[Theorem~4.1]{WY_nilp} every such map with $H(0)=0$ and linearly independent
components is linearly conjugate to a member of this family. Up to a triangular
conjugation of degree $2$ these maps are skew products over plane automorphisms of
H\'enon type (Lemma~\ref{lem:heyan-conj}). They have no left elementary factor, so
$\ell(F)=\#\supp(F)+1=4$, and their dynamical degree is $2d$
(Proposition~\ref{prop:heyan}). Over an infinite field they have infinitely many minimal
decompositions, as every map of positive defect, but up to two natural moves these form
exactly two classes (Proposition~\ref{prop:heyan-min}), while the six minimal
decompositions of Vasyunin's map form one class. Together with the classification of Wang
and Yan this gives the following dichotomy (Corollary~\ref{cor:dim3-dichotomy}). In
characteristic zero, a map $X+H\in K[X_1,X_2,X_3]^3$ with $H(0)=0$ and nilpotent $J_H$ is Pascal
finite if and only if the components of $H$ are linearly dependent, and otherwise its
dynamical degree is $\tfrac12\deg H$.

Pascal finiteness enters the paper in two ways. An automorphism $F$ is \emph{Pascal
finite} if $\Delta_F^m(X)=0$ for some $m$, where $\Delta_F=\sigma_F-\Id$ and
$\sigma_F(P)=P\circ F$ \cite{ABCH,ABCH2}. In characteristic zero these are exactly the
exponentials of locally nilpotent derivations \cite[Corollary~2]{AH_JvdK}. Elementary maps are the
simplest Pascal finite maps, of Pascal depth $2$ \cite{ABCH2}, and the class is not closed
under composition \cite{ABCH2}. So a decomposition into elementary maps is a factorization
into Pascal finite factors of depth $2$. For $n=2$ every automorphism with Jacobian
determinant $1$ is a composition of Pascal finite maps \cite{AH_JvdK}, and such a
factorization can be chosen canonically \cite{AH_canon}. Theorem~\ref{thmB} describes all minimal
factorizations of this kind of a map which is not Pascal finite. On the other hand,
for the maps $X+H$ with $H(0)=0$ and nilpotent $J_H(0)$, Pascal finiteness is equivalent
to the boundedness of the degrees of the iterates (Lemma~\ref{lem:pf-lf}). It separates
the maps whose iterates have bounded degree from the others, and in dimension three, by
the dichotomy above, it is decided by a linear condition on $H$.

\subsection*{The Anick automorphism}
In dimension four we treat the automorphism
$A=(X_1+\delta X_4,\,X_2+\delta X_3,\,X_3,\,X_4)$, $\delta=X_1X_3-X_2X_4$, proposed by
Anick as a candidate for a wild automorphism \cite[p.~343]{Coh} and shown to be tame by
Miasnikov and Shestakov if $\operatorname{char}K\ne2$ \cite{MS}.

\begin{theorem}[The Anick automorphism]\label{thmE}
\begin{enumerate}[label=\textup{(\roman*)},leftmargin=2em]
\item Over every field $K$, no word of length at most $6$ in $\cE(K,4)$ represents $A$.
If $\operatorname{char}K\ne2$, then $A\in\EE(K,4)$ and $7\le\ell(A)\le13$.
\item Over every field $K$ and for every $r\ge1$, $\ell(A\times\Id_r)=6$. Every
decomposition of $A\times\Id_r$ of length $6$ changes $X_1$, $X_2$ and exactly one of
the coordinates $X_5,\dots,X_{4+r}$, each of them exactly twice.
\item Every decomposition of $A$ into elementary maps changes each of $X_1$ and $X_2$ at
least twice, and $X_3$ or $X_4$ at least twice.
\end{enumerate}
\end{theorem}

So the addition of one variable shortens the minimal decompositions of $A$, while
$\ell(F\times\Id_r)=\ell(F)$ for every $F$ of defect at most $2$ (Theorem~\ref{thmA}).
The mechanism is linear. Over $R=K[X_3,X_4]$ the map $A$ acts on $(X_1,X_2)$ by the matrix
$M_A=I+vw^{\mathsf T}$, where $v=(X_4,X_3)^{\mathsf T}$, $w=(X_3,-X_4)^{\mathsf T}$ and
$w^{\mathsf T}v=0$, and $M_A$ does not lie in the subgroup $\mathrm{E}_2(R)$ generated by
elementary matrices, by a theorem of Cohn \cite{Coh66}. Part (iii) reduces to this fact,
and the lower bound in (i) combines it with comparisons of leading coefficients in
$X_1,X_2$ over $R$. In $\mathrm{E}_3(R)$ the matrix $M_A\oplus1$ is a commutator of
two block unipotent matrices,
\[
\begin{pmatrix}M_A&0\\0&1\end{pmatrix}
=\begin{pmatrix}I&v\\0&1\end{pmatrix}^{-1}
\begin{pmatrix}I&0\\w^{\mathsf T}&1\end{pmatrix}^{-1}
\begin{pmatrix}I&v\\0&1\end{pmatrix}
\begin{pmatrix}I&0\\w^{\mathsf T}&1\end{pmatrix},
\]
and the decomposition of length $6$ in (ii), due to Smith \cite{Sm}, is the polynomial
form of this identity (Section~\ref{sec:anick}). The lower block matrix is
a product of two elementary matrices, but it corresponds to a single elementary map,
which adds $\delta=w^{\mathsf T}(X_1,X_2)^{\mathsf T}$ to the new coordinate. Thus the difference between $\mathrm{E}_2(R)$ and $\mathrm{E}_3(R)$
discovered by Cohn reappears as a strict decrease of the elementary length under the
addition of one variable.

\subsection*{The plane, Steinberg relations, and iteration}
For $n\ge3$ the commutator of two elementary maps can be elementary
(Lemma~\ref{lem:commutator}). For linear maps this is the Steinberg relation
$[x_{ab}(f),x_{bk}(c)]=x_{ak}(fc)$, and $\EE(K,n)$ satisfies polynomial analogues of
these relations, which shorten words. They require a third index, and in the plane they
do not exist.

\begin{theorem}[Length in the plane]\label{thmC}
\begin{enumerate}[label=\textup{(\roman*)},leftmargin=2em]
\item Let $n=2$ and let $\ell_J(F)$ be the number of de~Jonqui\`eres letters of
$F\in\EE(K,2)$. Then $\EE(K,2)=\{F\in\Aut_KK[X_1,X_2]^2:\det J_F=1\}$, and
\[
\ell_J(F)\le\ell(F)\le4\,\ell_J(F)+9,\qquad \ell_J(F)\le\log_2\deg F .
\]
\item If $n=2$ and $Q,Q'\in\cE(K,2)$ are non-affine and change different coordinates,
then $\ell_J(Q\circ Q'\circ Q^{-1}\circ Q'^{-1})=4$, so the commutator is not elementary.
For every $n\ge3$ there are non-affine $Q,Q'\in\cE(K,n)$ changing different coordinates
whose commutator is elementary and non-affine.
\item For every $m\ge1$ there is $F_m\in\EE(K,2)$ with $\#\supp(F_m)=2$,
$\deg F_m=4^m$ and $\ell(F_m)=2m$.
\end{enumerate}
\end{theorem}

Since $\ell(F\circ G)\le\ell(F)+\ell(G)$, the limit
$d\ell(F)=\lim_m\frac1m\ell(F^{(m)})$ exists. It is the analogue of the dynamical length
in the Cremona group \cite{BF}. In the plane $d\ell(F)>0$ if and only if the degrees of
the iterates of $F$ grow exponentially, and the possible values of the pair
$\bigl(\lim_m\frac1m\ell_J(F^{(m)}),\lambda(F)\bigr)$ are described by an elementary
arithmetic condition (Remark~\ref{rem:dyn-length}). For $n\ge3$ the elementary length bounds the
displacement of a vertex of the simplicial complex of Lamy \cite{Lam19}, so for $n=3$
every element acting as a loxodromic isometry of the hyperbolic complex of
\cite{LP19} has $d\ell(F)>0$ (Remark~\ref{rem:complex}). For Vasyunin's map the degrees
of the iterates grow exponentially, and their exact values are known.

\begin{theorem}[Iterates of Vasyunin's map]\label{thmD}
Let $\operatorname{char}K\ne2$. The degrees of the iterates of $F_K$ are
$\deg F_K^{(1)}=2$, and for $j\ge1$
\[
\deg F_K^{(2j)}=2^{j+2}-4,\qquad \deg F_K^{(2j+1)}=6\cdot2^{j}-4 .
\]
For every invertible diagonal matrix $L$ the iterates of $L\circ F_K$ have the same
degrees. Consequently
$F_K$ has dynamical degree $\sqrt2$, and $L\circ F_K$ is not conjugate in $\Aut_KK[X_1,\dots,X_5]^5$
to a triangular map, nor to a linear one. For $K=\FF_p$, $p$ odd, the map
$F_K$ is a tame automorphism of infinite order in $\Aut_{\FF_p}\FF_p[X_1,\dots,X_5]^5$, and its
$J_H$ is nilpotent and not strongly nilpotent.
\end{theorem}

In \cite{AH_sn} only the unboundedness of these degrees over $\QQ$ was shown.
Theorem~\ref{thmD} gives the whole degree sequence in every characteristic $\ne2$. The
dynamical degree $\sqrt2$ also occurs for quadratic automorphisms of $K[X_1,X_2,X_3]$
\cite[Theorem~2]{BvS}. For $m\ge2$ and $\operatorname{char}K=0$ the iterates have no left
elementary factor, so $5\le\ell(F^{(m)})\le4m$ (Proposition~\ref{prop:iterate-left}).
Whether $\ell(F^{(m)})$ grows linearly in $m$ remains open
(Question~\ref{q:dyn-length}).

\subsection*{Organization}
Section~\ref{sec:prel} fixes notation and recalls the properties of Vasyunin's map.
Section~\ref{sec:left} proves Theorem~\ref{thmA} and treats decompositions of defect $1$,
direct products and moves. Theorem~\ref{thmB} is proved in Section~\ref{sec:vasyunin}
and Theorem~\ref{thmD} in Section~\ref{sec:iterates}. Section~\ref{sec:heyan} treats the
maps of He and Yan, and Section~\ref{sec:anick} the Anick automorphism, where
Theorem~\ref{thmE} is proved. Section~\ref{sec:plane} compares with the plane and proves
Theorem~\ref{thmC}. Section~\ref{sec:questions} summarizes the results, discusses the
dependence of the length on the base field and the stable length
$\lim_r\ell(F\times\Id_r)$, and states three questions.
Appendix~\ref{app:record} collects the inverse of Vasyunin's map and the polynomial
identities used in the proofs. Appendix~\ref{app:anick} shortens the word for the Anick
automorphism to $13$ factors and treats words which do not change $X_3$ and $X_4$.

\section{Preliminaries}\label{sec:prel}

Let $X=(X_1,\dots,X_n)$. As in \cite{AH_sn}, a polynomial map is the formal object, a tuple
$F=(F_1,\dots,F_n)\in K[X]^n$, and two polynomial maps are equal when their components
have the same coefficients. Polynomial maps are composed as $(F\circ G)(X)=F(G(X))$, and
$F^{(m)}$ is the $m$-th iterate. For $r\ge0$, $F\times\Id_r$ is the tuple
$(F_1,\dots,F_n,X_{n+1},\dots,X_{n+r})$, and we write $F\times\Id$ when $r$ is fixed.
A polynomial map $F$ is an \emph{automorphism} if there is $G\in K[X]^n$ with
$F\circ G=G\circ F=X$, and the automorphisms form the group $\Aut_KK[X]^n$ under
composition. The pullback $\sigma_F(P)=P\circ F$ is a $K$-algebra automorphism of $K[X]$
exactly when $F$ is an automorphism, and $\sigma_{F\circ G}=\sigma_G\circ\sigma_F$, so
$F\mapsto\sigma_F$ is an anti-isomorphism of $\Aut_KK[X]^n$ onto $\Aut_KK[X]$. Invertibility
always means the existence of a polynomial inverse, never the bijectivity of the
evaluation map $x\mapsto F(x)$ of $K^n$. The distinction matters here, since the paper
also works over finite fields. The order of an automorphism is its order in
$\Aut_KK[X]^n$. For finite $K$ it may differ from the order of the evaluation map, which
is always finite (Corollary~\ref{cor:modp}).

The degree $\deg F$ is the largest degree of a component of $F$. For $0\ne P\in K[X]$ let $\ord P$ be the least degree of a monomial
occurring in $P$. For an automorphism $F$ one has
$\deg F^{(m+k)}\le\deg F^{(m)}\deg F^{(k)}$, so the sequence $\log\deg F^{(m)}\ge0$ is
subadditive, and the limit
\[
\lambda(F)=\lim_{m\to\infty}\bigl(\deg F^{(m)}\bigr)^{1/m}
\]
exists by Fekete's lemma. It is the \emph{dynamical degree} of $F$. For $H\in K[X]^n$ let
$J_H$ be its Jacobian matrix. It is \emph{strongly nilpotent} if
$J_H(Y^{(1)})\cdots J_H(Y^{(n)})=0$, where $Y^{(1)},\dots,Y^{(n)}$ are $n$ disjoint sets
of $n$ new variables \cite{EH}. Write $J_H=\sum_\alpha X^\alpha A_\alpha$ with
$A_\alpha\in M_n(K)$. Then
$J_H(Y^{(1)})\cdots J_H(Y^{(n)})=\sum(Y^{(1)})^{\alpha_1}\cdots(Y^{(n)})^{\alpha_n}
A_{\alpha_1}\cdots A_{\alpha_n}$, and the monomials of this sum are pairwise distinct. So
$J_H$ is strongly nilpotent if and only if every product of $n$ of the matrices $A_\alpha$
vanishes, that is, if and only if the algebra $\mathcal A_H\subseteq M_n(K)$, without
unit, generated by the $A_\alpha$ is nilpotent. Indeed, if $\mathcal A_H^m=0$, the
subspaces $\mathcal A_H^iK^n$ decrease strictly until they vanish, so
$\mathcal A_H^nK^n=0$. Ordinary nilpotency of $J_H$ only requires the vanishing of the
sums $\sum A_{\alpha_1}\cdots A_{\alpha_n}$ over the tuples with a fixed value of
$\alpha_1+\dots+\alpha_n$. We write
$\Delta_F=\sigma_F-\Id$, and $\tau_K(F)$ for the Pascal depth, the least $m$ with
$\Delta_F^m(X)=0$ \cite{ABCH,ABCH2}. The map $F$ is \emph{Pascal finite} if
$\tau_K(F)<\infty$.

\begin{definition}\label{def:qt}
A map $F=X+H$, $H\ne0$, is a \emph{quasi-translation} if it is an automorphism with
inverse $X-H$. It is an \emph{elementary map} if there is an index $j$ with $H_i=0$ for
$i\ne j$ and $H_j\in K[X_i:i\ne j]$. Elementary maps are quasi-translations. We say that $F$ \emph{changes} $X_j$, and call $H_j$ its
\emph{tail}. It is \emph{affine} if $\deg H_j\le1$.
\end{definition}

The map $X+H$ is a quasi-translation if and only if $H\circ F=H$, if and only if
$\tau_K(F)=2$ \cite{ABCH2}. Every elementary map is a quasi-translation. A triangular
map with all diagonal coefficients $1$ is Pascal finite, Pascal finiteness is invariant
under conjugation, and a composition of Pascal finite maps need not be Pascal finite
\cite[Theorems~2.1, 2.2, Remark~3.1]{ABCH2}. The length defined below does not bound the
Pascal depth. For $d\ge1$ the map $F_d=(X_1+X_2^d,\,X_2+X_3,\,X_3)$ is the composition of
two elementary maps, while $\Delta_{F_d}$ acts on $K[X_2,X_3]$ as the difference operator
with step $X_3$ in $X_2$, so $\Delta_{F_d}^{d+1}(X_1)=d!\,X_3^d$ and
$\Delta_{F_d}^{d+2}(X)=0$. Hence $\tau_K(F_d)=d+2$ if $\operatorname{char}K=0$ or
$\operatorname{char}K>d$.

A polynomial map $F$ is \emph{locally finite} if the iterates $\sigma_F^m(P)$, $m\ge0$,
span a finite-dimensional subspace of $K[X]$ for every $P\in K[X]$ \cite{FM07}. In
characteristic zero, for the maps $X+H$ considered below, local finiteness and Pascal
finiteness coincide.

\begin{lemma}\label{lem:pf-lf}
Let $\operatorname{char}K=0$ and $F=X+H$ with $H(0)=0$ and $J_H(0)$ nilpotent. The
following are equivalent.
\textup{(i)} $F$ is Pascal finite.
\textup{(ii)} $F$ is locally finite.
\textup{(iii)} The degrees of the iterates $F^{(m)}$ are bounded.
\end{lemma}

\begin{proof}
A Pascal finite map is locally finite \cite{ABCH2}, and (ii) and (iii) are equivalent by
\cite[Theorem~1.1]{FM07}. Let $F$ be locally finite. Then $\det J_F$ is constant
\cite[Corollary~1.1]{FM07}, equal to $\det(I+J_H(0))=1$, so $F$ is an automorphism
\cite[Corollary~1.2]{FM07}. Since $F(0)=0$ and its linear part $I+J_H(0)$ is
unipotent, $F$ is unipotent \cite[Proposition~2.1]{FM07}, hence the exponential of a
locally nilpotent derivation \cite[Theorem~2.3]{FM07}, and therefore Pascal finite
\cite[Corollary~2]{AH_JvdK}.
\end{proof}

If $\ord H_i\ge2$ for every $i$, the equivalence of (i) and (ii) is
\cite[Proposition~2.1]{ABCH}.

The lemma applies whenever $J_H$ is nilpotent and $H(0)=0$, in particular to Vasyunin's
map and to the maps of Section~\ref{sec:heyan}. For them the statements ``not Pascal
finite'', ``not locally finite'' and ``iterates of unbounded degree'' are equivalent, and
each implies infinite order, since a map of finite order is locally finite.

Let $\cE=\cE(K,n)$ be the set of elementary maps and $\EE(K,n)$ the group they generate.
For $F\in\EE(K,n)$ put
\[
\ell(F)=\min\{l:F=Q_l\circ\dots\circ Q_1,\ Q_i\in\cE\},\qquad
\supp(F)=\{i:F_i\ne X_i\}.
\]
A word $Q_l\circ\dots\circ Q_1$ with $Q_i\in\cE$ is a \emph{decomposition} of $F$ of
length $l$, and $Q_1$ is applied first.

The group $\EE(K,n)$ has the following description, the case $R=K$ of
\cite[Proposition~3.20]{BEW}. It rests on the facts that $\mathrm{SL}_n(K)$ is generated
by elementary matrices and that conjugation by a diagonal map preserves elementary maps.

\begin{proposition}\label{prop:Gamma}
Let $n\ge2$. Then $\EE(K,n)=\{F\in\mathrm{TA}_n(K):\det J_F=1\}$.
\end{proposition}

If $\operatorname{char}K=0$, then $\EE(K,3)$ is the special tame group
$\mathrm{STame}(K^3)$ of Lamy and Przytycki \cite{LP19}, the kernel of the Jacobian
determinant on $\mathrm{TA}_3(K)$.

By Proposition~\ref{prop:Gamma} and Corollary~\ref{cor:modp}, the reductions modulo odd
primes of Vasyunin's map lie in $\EE(\FF_p,5)$.

For $\sigma\in S_n$ let $J_\sigma(K,n)$ be the group of automorphisms $G$ such that
\[
G_{\sigma(i)}=a_iX_{\sigma(i)}+f_i(X_{\sigma(1)},\dots,X_{\sigma(i-1)}),\qquad a_i\in K^*,
\]
and $J^1_\sigma(K,n)$ its subgroup defined by $a_1=\dots=a_n=1$. An elementary map
changing $X_j$ lies in $J^1_\sigma$ if and only if every variable of its tail precedes
$X_j$ in the order $X_{\sigma(1)},\dots,X_{\sigma(n)}$. For $F=X+H$ the
\emph{dependency digraph} $D_F$ has the vertices $1,\dots,n$ and an edge $i\to j$
whenever $X_i$ occurs in $H_j$, loops included. If $\sigma(i)=n+1-i$, then
$J_\sigma(K,n)$ is the de~Jonqui\`eres subgroup $J(K,n)$ of \cite[Chapter~5]{vdE}, and
every $J_\sigma(K,n)$ is conjugate to $J(K,n)$ by a permutation of the coordinates.

\begin{lemma}\label{lem:dag}
Let $F=X+H$ be an automorphism such that, for every $i$, $H_i=0$ or $\ord H_i\ge2$.
Then $F\in J_\sigma$
for some $\sigma$ if and only if $D_F$ is acyclic, and then $F\in J^1_\sigma$ and
$\ell(F)=\#\supp(F)$.
\end{lemma}

The inequality $\ell(F)\ge\#\supp(F)$ used here is Lemma~\ref{lem:support-bound} below.

\begin{proof}
If $D_F$ is acyclic, let $\sigma$ be a topological order of its vertices. Since there
are no loops, every variable of $H_{\sigma(i)}$ precedes $X_{\sigma(i)}$, so
$F\in J^1_\sigma$. Moreover $F=E_{\sigma(1)}\circ\dots\circ E_{\sigma(n)}$, where
$E_{\sigma(i)}$ adds $H_{\sigma(i)}$ to $X_{\sigma(i)}$ and factors with
$H_{\sigma(i)}=0$ are omitted. Here $E_{\sigma(n)}$ is applied first and
$E_{\sigma(1)}$ last, so the variables of $H_{\sigma(i)}$, which precede $X_{\sigma(i)}$,
are not yet changed when $E_{\sigma(i)}$ is applied. Conversely, if $F\in J_\sigma$, then
$H_{\sigma(i)}=(a_i-1)X_{\sigma(i)}+f_i$, and either $H_{\sigma(i)}=0$ or
$\ord H_{\sigma(i)}\ge2$. In both cases $a_i=1$ and $H_{\sigma(i)}=f_i$. So every edge of $D_F$ points forward in the order $\sigma$.
\end{proof}

The assumption on $\ord H_i$ cannot be omitted. If $\operatorname{char}K\ne2$, the map
$(2X_1,X_2)=X+(X_1,0)$ lies in $J_{\mathrm{id}}(K,2)$, while $D_F$ has a loop at $1$.

\subsection{Vasyunin's map}\label{ssec:vasyunin-def}

Throughout, $F=X+H$ denotes the map
\begin{equation}\label{ve}
H_1=0,\quad H_2=X_1X_3,\quad H_3=X_1X_4+\tfrac12X_2^2,\quad H_4=X_1X_5-X_2X_3,\quad
H_5=\tfrac12X_3^2 .
\end{equation}
By \cite[Section~5.3]{AH_sn} it is an automorphism of $\QQ^5$ with an explicit inverse
$G=F^{-1}$, whose coordinates have the degrees $(1,4,3,7,6)$. The matrix $J_H$ is
nilpotent of index $5$. It is not strongly nilpotent, since the algebra
$\mathcal A_H$, generated by the matrices $J_H(e_i)$, contains a non-nilpotent element.
Indeed all powers of
$J_H(e_1)J_H(e_2)=\mathrm{diag}(0,1,-1,0,0)$ are nonzero. A second witness is
$J_H(e_3)J_H(e_2)J_H(e_1)^3=E_{55}$, which has integer entries. Moreover $F$ is not
Pascal finite \cite[Section~5.3]{AH_sn}. A \emph{quadratic homogeneous} map is a map $X+H$ with
all $H_i$ homogeneous of degree $2$. De~Bondt proved that every quadratic homogeneous
automorphism of $K[X_1,\dots,X_5]$, $\operatorname{char}K=0$, is tame \cite[Theorem~4.6.8]{Bt}, see
also \cite{Sun}. In dimension $n\le4$, $\operatorname{char}K=0$, every quadratic homogeneous automorphism is
linearly triangularizable \cite{MO}, see \cite[Proposition~1.1]{Sun} and
\cite[Example~7.4.10]{vdE}. So $n=5$ is the
least dimension in which a quadratic homogeneous automorphism can fail to be linearly
triangularizable, and $F$ is such a map by Corollary~\ref{cor:shifted}. Sun classified
the quadratic homogeneous automorphisms of $\mathbb{C}[X_1,\dots,X_5]$ up to linear conjugation
\cite{Sun14}. Among those which are not linearly triangularizable are the maps
\[
\begin{aligned}
S_2=\bigl(&X_1,\;X_2+\alpha X_4^2+o_1X_1^2,\;X_3-2\alpha X_4X_5+X_1X_2+o_2X_1^2,\\
&X_4+\alpha X_5^2+X_1X_3+o_3X_1^2,\;X_5+X_1X_4+o_4X_1^2\bigr),
\end{aligned}
\]
$\alpha\ne0$ \cite[Theorem~3.2]{Sun14}, see \cite[Lemma~4.3]{SS}. Exchanging $X_2$ with
$X_5$ and $X_3$ with $X_4$ transforms $F$ into $S_2$ with $\alpha=\tfrac12$ and
$o_1=\dots=o_4=0$, as a direct substitution shows.

The coefficients $\tfrac12$ are determined by the other coefficients of $H$, in the
following sense.

\begin{proposition}\label{prop:coefficients}
Let $\operatorname{char}K\ne2$, $c,d\in K$ and
$H^{c,d}=(0,\,X_1X_3,\,X_1X_4+cX_2^2,\,X_1X_5-X_2X_3,\,dX_3^2)$. Then
\[
\det(tI-J_{H^{c,d}})=t^5+(1-2c)X_1X_2\,t^3+(1-2d)X_1^2X_3\,t^2 ,
\]
and the following conditions are equivalent. \textup{(i)} $J_{H^{c,d}}$ is nilpotent.
\textup{(ii)} $\det J_{X+H^{c,d}}=1$. \textup{(iii)} $X+H^{c,d}$ is an automorphism.
\textup{(iv)} $c=d=\tfrac12$.
\end{proposition}

\begin{proof}
The first row of $J=J_{H^{c,d}}$ is zero, and the other rows are
\[
(X_3,0,X_1,0,0),\quad(X_4,2cX_2,0,X_1,0),\quad(X_5,-X_3,-X_2,0,X_1),\quad(0,0,2dX_3,0,0).
\] The nonzero principal
minors of order $2$ are those on $\{2,3\}$ and $\{3,4\}$, equal to $-2cX_1X_2$ and
$X_1X_2$, and those of order $3$ are on $\{2,3,4\}$ and $\{3,4,5\}$, equal to $-X_1^2X_3$
and $2dX_1^2X_3$. Principal minors containing the index $1$ vanish, since the first row
of $J$ is zero, and the minor on
$\{2,3,4,5\}$ vanishes by expansion along its last column. This gives the
characteristic polynomial. Since $J$ is linear in $X$, $\det(I+tJ(X))=\det(I+J(tX))$,
so $\det(I+J)=1$ holds if and only if all non-leading coefficients of the characteristic
polynomial vanish. By the formula above these coefficients are $(1-2c)X_1X_2$ and
$(1-2d)X_1^2X_3$, and they vanish exactly when $c=d=\tfrac12$. Hence (i), (ii) and (iv)
are equivalent. An automorphism has constant Jacobian determinant, which equals $1$
here since $J$ vanishes at $0$, so (iii) implies (ii). For $c=d=\tfrac12$ the map is
$F$, an automorphism by Section~\ref{ssec:three}, so (iv) implies (iii).
\end{proof}

So within the family $H^{c,d}$, $cd\ne0$, the support of $H$ and the digraph $D_F$ are
constant, while only one member is an automorphism.

\section{Left elementary factors and minimal length}\label{sec:left}

\begin{lemma}\label{lem:support-bound}
If $F=Q_l\circ\dots\circ Q_1$ with $Q_i\in\cE$ changing $X_{m_i}$, then
$\supp(F)\subseteq\{m_1,\dots,m_l\}$. In particular $\ell(F)\ge\#\supp(F)$, and if
$l=\#\supp(F)$ then the $m_i$ are pairwise distinct and $\{m_1,\dots,m_l\}=\supp(F)$.
\end{lemma}

\begin{proof}
A coordinate changed by no $Q_i$ is fixed by $F$.
\end{proof}

\begin{lemma}\label{lem:word}
Let $F=Q_l\circ\dots\circ Q_1$ with $Q_i\in\cE(K,n)$ changing $X_{m_i}$.
\begin{enumerate}[label=\textup{(\roman*)},leftmargin=2em]
\item If $X_j$ is changed by exactly one factor, then $F_j\ne X_j$. Hence every index
outside $\supp(F)$ occurs in $(m_1,\dots,m_l)$ either never or at least twice.
\item If $m_i=m_{i+1}$ for some $i$, then $Q_{i+1}\circ Q_i$ is elementary or equal to
$\Id$, and $F$ has a decomposition of length at most $l-1$. In particular, adjacent
factors of a decomposition of length $\ell(F)$ change different coordinates.
\end{enumerate}
\end{lemma}

\begin{proof}
(i) Let $Q_i$ be the only factor changing $X_j$, by a tail $t$. Then $t\ne0$, since
$Q_i\ne\Id$. The factors $Q_{i+1},\dots,Q_l$ fix $X_j$, and $G=Q_{i-1}\circ\dots\circ Q_1$
satisfies $G_j=X_j$. Hence $F_j=(Q_i)_j\circ G=X_j+t\circ G$, and $t\circ G\ne0$, since
composition with the automorphism $G$ is injective on $K[X]$.
(ii) If $Q_i$ and $Q_{i+1}$ change $X_j$ by tails $t$ and $t'$, which do not involve
$X_j$, then $Q_{i+1}\circ Q_i$ changes $X_j$ by $t+t'$ and fixes the other coordinates.
\end{proof}

\begin{proposition}\label{prop:length-spectrum}
Let $n\ge2$ and $\Id\ne F\in\EE(K,n)$. Then $F$ has decompositions of every length
$l\ge\ell(F)$.
\end{proposition}

\begin{proof}
It suffices to pass from a decomposition of length $l$ to one of length $l+1$. Let $Q$
be a factor of a decomposition of length $l$, changing $X_j$ by $h\ne0$. Choose $k\ne j$
and $h'\in\{X_k,X_k+1\}$ with $h'\ne h$. Then $h'\ne0$ and $h-h'\ne0$, and both lie in
$K[X_i:i\ne j]$. So the maps $Q'$ and $Q''$ changing $X_j$ by $h'$ and by $h-h'$ belong
to $\cE$. Since neither of them changes the variables of the tail of the other,
$Q'\circ Q''$ changes $X_j$ by $h'+(h-h')=h$, that is $Q'\circ Q''=Q$. Replacing $Q$ by
$Q'\circ Q''$ gives a decomposition of length $l+1$.
\end{proof}

So the set of lengths of decompositions carries no information beyond $\ell(F)$. If
$\deg h\ge2$, the factors $Q'$ and $Q''$ have degrees at most $\deg h$, so the same
holds for decompositions into factors of degree at most $d$, for every $d\ge2$.

If $F=Q_l\circ\dots\circ Q_1$ and $Q_l$ is the only factor changing $X_j$, then
$Q_l$ can be removed from the left without changing the other coordinates of $F$. The
following definition isolates this situation.

\begin{definition}\label{def:left-factor}
Let $F\in\Aut_KK[X]^n$ and $1\le j\le n$. An elementary map $Q$ changing $X_j$ is a
\emph{left elementary factor} of $F$ in the direction $X_j$ if $F=Q\circ G$ for an
automorphism $G$ with $G_j=X_j$.
\end{definition}

\begin{lemma}\label{lem:left-factor}
Let $F\in\Aut_KK[X]^n$ and $1\le j\le n$. Then $F$ has a left elementary factor in the
direction $X_j$ if and only if
\[
h:=X_j-(F^{-1})_j\in K[X_i:i\ne j],\qquad h\ne0 .
\]
In that case the left elementary factor is unique, it is the map $Q$ changing $X_j$ by
$h$, and $G=Q^{-1}\circ F$ satisfies $\supp(G)=\supp(F)\setminus\{j\}$.
\end{lemma}

\begin{proof}
Let $F=Q\circ G$ with $Q$ changing $X_j$ by $h$ and $G_j=X_j$. Then $h\ne0$, since
$Q\ne\Id$, and $h\in K[X_i:i\ne j]$ by definition of an elementary map. Moreover
$(G^{-1})_j=X_j$, since $X_j=G_j\circ G^{-1}=(G^{-1})_j$, and
$(F^{-1})_j=(G^{-1}\circ Q^{-1})_j=X_j-h$. Conversely, let $h$ be as stated, $Q$ the map
changing $X_j$ by $h$ and $G=Q^{-1}\circ F$. Then $G_i=F_i$ for $i\ne j$, and
$G_j=F_j-h(F)=X_j$, by composing $(F^{-1})_j=X_j-h$ with $F$.
\end{proof}

In particular $h=0$ for $j\notin\supp(F)$, so left elementary factors exist only in the
directions of $\supp(F)$.

\begin{theorem}\label{thm:min-decomp}
Let $F\in\Aut_KK[X]^n$ and $k=\#\supp(F)\ge1$. Call an ordering $(j_k,\dots,j_1)$ of
$\supp(F)$ \emph{admissible} if left elementary factors can be extracted successively
in the directions $X_{j_k},\dots,X_{j_1}$. That is, $G^{(k)}=F$, and for
$i=k,\dots,1$ the map $G^{(i)}$ has a left elementary factor $Q_i$ in the direction
$X_{j_i}$, and $G^{(i-1)}=Q_i^{-1}\circ G^{(i)}$. The decompositions
$F=Q_k\circ\dots\circ Q_1$ of length $k$ are in bijection with the admissible orderings.
Given $(G^{(i)})^{-1}$, the test in the direction $X_j$ is whether
$h=X_j-((G^{(i)})^{-1})_j$ is nonzero and does not involve $X_j$, and then
$(G^{(i-1)})^{-1}=(G^{(i)})^{-1}\circ Q_i$, so each step uses only the inverse of the
previous one.
\end{theorem}

The extraction starts at the left end of the word. The factor $Q_k$, which is applied
last, is removed first, and $Q_1$, which is applied first, is removed last.

\begin{proof}
By Lemma~\ref{lem:support-bound} the factors of a decomposition of length $k$ change the
coordinates of $\supp(F)$, each exactly once. So $Q_{k-1}\circ\dots\circ Q_1$ fixes
$X_{j_k}$, and $Q_k$ is the left elementary factor of $F$ in the direction $X_{j_k}$,
which is unique by Lemma~\ref{lem:left-factor}. The map $Q_k^{-1}\circ F$ is again an
automorphism, as a composition of automorphisms. By induction every decomposition of length $k$ arises from an admissible
ordering, which determines it. Conversely, if an ordering is admissible, then
$\supp(G^{(0)})=\emptyset$ by Lemma~\ref{lem:left-factor}, so $G^{(0)}=\Id$.
\end{proof}

For $1\le j\le n$ let $U_j\subseteq\Aut_KK[X]^n$ be the group of the maps changing only $X_j$
by an element of $K[X_i:i\ne j]$. It is abelian and isomorphic to the additive group of
$K[X_i:i\ne j]$, and $\EE(K,n)$ is generated by $U_1,\dots,U_n$. In this language
Theorem~\ref{thm:min-decomp} says the following. If $(j_k,\dots,j_1)$ is an ordering of
$\supp(F)$, $k=\#\supp(F)$, then $F$ has at most one representation $F=Q_k\circ\dots\circ Q_1$
with $Q_i\in U_{j_i}$, and $Q_k,\dots,Q_1$ are obtained from $F^{-1}$ by the tests of the
theorem.

\begin{corollary}\label{cor:stab}
Let $F\in\EE(K,n)$ and $r\ge1$. The decompositions of $F\times\Id_r$ of length
$\#\supp(F)$ are the words $(Q_k\times\Id)\circ\dots\circ(Q_1\times\Id)$ for the
decompositions $Q_k\circ\dots\circ Q_1$ of $F$ of length $k=\#\supp(F)$. In particular
$\ell(F\times\Id_r)=k$ if and only if $\ell(F)=k$, and $\ell(F)=k+1$ implies
$\ell(F\times\Id_r)=k+1$.
\end{corollary}

\begin{proof}
Clearly $\supp(F\times\Id)=\supp(F)$ and $(F\times\Id)^{-1}=F^{-1}\times\Id$. For $j\le n$
the polynomial $X_j-((F\times\Id)^{-1})_j=X_j-(F^{-1})_j$ involves only $X_1,\dots,X_n$.
So each test of Theorem~\ref{thm:min-decomp} has the same outcome for $F\times\Id$ as for
$F$, the factor obtained is $Q\times\Id$, and the remaining map is again of the form
$G\times\Id$. For the last statement, $\ell(F\times\Id)\le\ell(F)$, and
$\ell(F\times\Id)>k$ by the first statement.
\end{proof}

\begin{proposition}\label{prop:stab-bound}
Let $F\in\EE(K,n)$, $F\ne\Id$, and $r\ge1$. Then
\[
\ell(F\times\Id_r)\ge\min\bigl(\ell(F),\,\#\supp(F)+2\bigr).
\]
In particular $\ell(F\times\Id_r)=\ell(F)$ if $\ell(F)\le\#\supp(F)+2$.
\end{proposition}

\begin{proof}
Let $F\times\Id_r=Q_l\circ\dots\circ Q_1$ with $Q_i\in\cE(K,n+r)$. Suppose first that
no $Q_i$ changes a coordinate $X_k$ with $k>n$. On the group of automorphisms of
$K[X_1,\dots,X_{n+r}]$ fixing $X_{n+1},\dots,X_{n+r}$ the map
$\varphi\mapsto\bar\varphi=(\varphi_1,\dots,\varphi_n)|_{X_{n+1}=\dots=X_{n+r}=0}$ is a
homomorphism to $\Aut_KK[X]^n$. Each $\bar Q_i$ is elementary or equal to $\Id$, and
$\overline{F\times\Id}=F$, so $l\ge\ell(F)$. Suppose now that some factor changes $X_k$ with
$k>n$. Since $(F\times\Id)_k=X_k$, at least two factors change $X_k$ by
Lemma~\ref{lem:word}(i).
By Lemma~\ref{lem:support-bound} at least $\#\supp(F)$ further factors change the
coordinates of $\supp(F)$, and $l\ge\#\supp(F)+2$.
\end{proof}

In the first case of the proof the additional coordinates enter only as parameters in
the tails. A word which changes them is not compatible with the specialization
$X_{n+1}=\dots=X_{n+r}=0$. For instance, let $n=2$ and $r=1$, so that $X_3$ is the
additional coordinate. The maps $P=(X_1,X_2,X_3+X_1)$ and $Q=(X_1,X_2+X_3,X_3)$ have
first two coordinates equal to $(X_1,X_2)$ after the substitution $X_3=0$, while
$Q\circ P=(X_1,X_2+X_1+X_3,X_3+X_1)$ is not of the form $G\times\Id_1$.

The arguments above count the factors in each direction. Let
$F=Q_l\circ\dots\circ Q_1$ with $Q_i\in\cE(K,n)$, and let $n_j$ be the number of factors
changing $X_j$.

\begin{lemma}\label{lem:defect}
One has $n_j\ge1$ for $j\in\supp(F)$ and $n_j\ne1$ for $j\notin\supp(F)$, and
\[
l-\#\supp(F)=\sum_{j\in\supp(F)}(n_j-1)+\sum_{j\notin\supp(F)}n_j .
\]
In particular, if $q$ indices outside $\supp(F)$ occur, then $l\ge\#\supp(F)+2q$.
\end{lemma}

\begin{proof}
The first statement is Lemma~\ref{lem:support-bound}, and the second is
Lemma~\ref{lem:word}(i). The formula follows from
$l=\sum_jn_j$, and each index outside $\supp(F)$ which occurs contributes at least $2$.
\end{proof}

We call $l-\#\supp(F)$ the \emph{defect} of the decomposition. Decompositions of defect
$0$ are described by Theorem~\ref{thm:min-decomp}. In a decomposition of defect $1$ no
index outside $\supp(F)$ occurs and exactly one index of $\supp(F)$ occurs twice, which is
the starting point of the next proposition.

\begin{corollary}\label{cor:no-relation}
Let $j_1,\dots,j_k$ be pairwise distinct. Then the map
$U_{j_k}\times\dots\times U_{j_1}\to\EE(K,n)$, $(Q_k,\dots,Q_1)\mapsto Q_k\circ\dots\circ Q_1$,
is injective. In particular, if $Q_k\circ\dots\circ Q_1=\Id$ with $Q_i\in U_{j_i}$, then
$Q_i=\Id$ for all $i$.
\end{corollary}

\begin{proof}
Let $F=Q_k\circ\dots\circ Q_1$. Omitting the factors equal to $\Id$ gives a decomposition
of $F$ in which every index occurs at most once. By Lemma~\ref{lem:defect} its length is
$\#\supp(F)$, and $\supp(F)$ is the set of the $j_i$ with $Q_i\ne\Id$. So $F$ determines
which factors are different from $\Id$, and Theorem~\ref{thm:min-decomp}, applied to the
induced ordering of $\supp(F)$, determines these factors.
\end{proof}

So every relation between the groups $U_1,\dots,U_n$ which is not a consequence of the
group laws in each $U_j$ involves some index at least twice. The relations of
Lemma~\ref{lem:commutator} below are of this kind.

For decompositions of length $\#\supp(F)+1$ the argument leaves one additional unknown
tail, and the following proposition gives the two resulting cases.

\begin{proposition}\label{prop:length-k1}
Let $F\in\Aut_KK[X]^n$, $k=\#\supp(F)\ge1$, $\ell(F)=k+1$, and let
$F=Q_{k+1}\circ\dots\circ Q_1$ with $Q_i\in\cE$ changing $X_{m_i}$. Then
$\{m_1,\dots,m_{k+1}\}=\supp(F)$, exactly one index $c$ occurs twice, and
$G=Q_{k+1}^{-1}\circ F$ satisfies one of the following.
\begin{enumerate}[label=\textup{(\alph*)},leftmargin=2em]
\item $m_{k+1}\ne c$, $Q_{k+1}$ is the left elementary factor of $F$ in the direction
$X_{m_{k+1}}$, $\supp(G)=\supp(F)\setminus\{m_{k+1}\}$, and
$Q_k\circ\dots\circ Q_1$ is a decomposition of $G$ of length $\#\supp(G)+1$.
\item $m_{k+1}=c\ne m_k$, $\supp(G)=\supp(F)$, and $Q_k\circ\dots\circ Q_1$ is a
decomposition of $G$ of length $\#\supp(G)$, given by Theorem~\ref{thm:min-decomp}.
\end{enumerate}
If $F$ has no left elementary factor, only \textup{(b)} occurs. In \textup{(b)}, if
$Q_{k+1}$ changes $X_c$ by $t$, then $G^{-1}=F^{-1}\circ Q_{k+1}$, and the tests of
Theorem~\ref{thm:min-decomp} for $G$ are conditions on $t$.
\end{proposition}

\begin{proof}
By Lemma~\ref{lem:support-bound}, $\supp(F)\subseteq\{m_1,\dots,m_{k+1}\}$. By
Lemma~\ref{lem:word}(i) an index outside $\supp(F)$ occurs at least twice. Since there are
$k+1$ factors, all $m_i$ lie in $\supp(F)$ and exactly one index $c$ occurs twice. If
$m_{k+1}$ occurs once, then $Q_k\circ\dots\circ Q_1$ fixes $X_{m_{k+1}}$, which gives
(a). Let $m_{k+1}=c$. Then $m_k\ne c$ by Lemma~\ref{lem:word}(ii). Moreover $G_j=F_j$ for $j\ne c$, and $G_c\ne X_c$ by Lemma~\ref{lem:word}(i), since
$X_c$ is changed by exactly one factor of $Q_k\circ\dots\circ Q_1$. Hence
$\supp(G)=\supp(F)$, which gives (b).
\end{proof}

In case (b) the tail $t$ is not determined by $F$. The proof of
Proposition~\ref{prop:anick}(ii) below applies it for $k=2$.

For $F\in\Aut_KK[X]^n$ and $G\in\Aut_KK[Y]^m$ let $F\times G\in\Aut_KK[X,Y]^{n+m}$ act by $F$ on
$X=(X_1,\dots,X_n)$ and by $G$ on $Y=(Y_1,\dots,Y_m)$. If $\ell(F)=\#\supp(F)$, let $N(F)$
be the number of decompositions of $F$ of length $\#\supp(F)$.

\begin{proposition}\label{prop:products}
Let $F\in\EE(K,n)$ and $G\in\EE(K,m)$ with $k_F=\#\supp(F)\ge1$ and
$k_G=\#\supp(G)\ge1$. Then $k_F+k_G\le\ell(F\times G)\le\ell(F)+\ell(G)$, and the
following hold.
\begin{enumerate}[label=\textup{(\roman*)},leftmargin=2em]
\item $\ell(F\times G)=k_F+k_G$ if and only if $\ell(F)=k_F$ and $\ell(G)=k_G$. In this
case the decompositions of $F\times G$ of length $k_F+k_G$ are obtained from pairs of
decompositions of $F$ and of $G$ of minimal length by interleaving their factors, and
\[
N(F\times G)=\binom{k_F+k_G}{k_F}N(F)\,N(G).
\]
\item If neither $F$ nor $G$ has a left elementary factor, then
$\ell(F\times G)\ge k_F+k_G+2$.
\end{enumerate}
\end{proposition}

\begin{proof}
The bounds follow from Lemma~\ref{lem:support-bound}, since
$\supp(F\times G)$ is the disjoint union of $\supp(F)$ and $\supp(G)$, and from the
concatenation of decompositions. Since $(F\times G)^{-1}=F^{-1}\times G^{-1}$, the test
of Lemma~\ref{lem:left-factor} for $F\times G$ in a direction $X_j$ is the test for $F$,
and the resulting factor changes $X_j$ by a polynomial in $X$.

(i) If $\ell(F)=k_F$ and $\ell(G)=k_G$, concatenation gives $\ell(F\times G)=k_F+k_G$.
Conversely, by Theorem~\ref{thm:min-decomp} a decomposition of $F\times G$ of length
$k_F+k_G$ is obtained by successive extraction of left elementary factors. By the
remark above and induction, every extracted factor is of the form $Q\times\Id$ or
$\Id\times P$, and the remaining map is again a product. Factors of different forms
commute, so the factors of the form $Q\times\Id$, in their order, give a decomposition of
$F$ of length $k_F$, and similarly for $G$. Conversely every interleaving of two such
decompositions is a decomposition of $F\times G$. A decomposition of $F\times G$ determines
the positions of its factors of each form, since their directions lie in different
blocks, and the two decompositions. This gives the formula for $N(F\times G)$.

(ii) The map $F\times G$ has no left elementary factor, so $\ell(F\times G)>k_F+k_G$ by
Theorem~\ref{thm:min-decomp}. Suppose that $\ell(F\times G)=k_F+k_G+1$. By
Proposition~\ref{prop:length-k1}, case (b) occurs. Let the last factor $Q$ change the
repeated coordinate, say in the block $X$, the other case being symmetric. Then
$G'=Q^{-1}\circ(F\times G)$ has a decomposition of length $\#\supp(G')=k_F+k_G$, obtained by
successive extraction of left elementary factors from $G'$. The last $m$ components of
$G'^{-1}=(F^{-1}\times G^{-1})\circ Q$ are $G^{-1}(Y)$, since $Q$ changes a coordinate of $X$.
As long as only factors changing coordinates of $X$ are extracted, the last $m$
components of the inverse of the remaining map stay equal to $G^{-1}(Y)$, and no
direction $Y_i$ passes the test, since $G$ has no left elementary factor. After at most
$k_F$ extractions no direction is left, while $k_G\ge1$ factors remain to be extracted.
This is a contradiction.
\end{proof}

So the defect $\ell(F)-\#\supp(F)$ is subadditive under direct products and vanishes on a
product if and only if it vanishes on both factors. For Vasyunin's map $F$
(Theorem~\ref{thm:V-minimal}) the $r$-fold product satisfies $\ell(F^{\times r})=4r$ and
$N(F^{\times r})=6^r(4r)!/(4!)^r$, for instance $N(F\times F)=2520$. For two maps $F,G$ of
the family of He and Yan in Section~\ref{sec:heyan}, which have no left elementary factor
and length $4$, Proposition~\ref{prop:products}(ii) gives $\ell(F\times G)=8$, so the
defects add. For such a map $F$ and an elementary map $E$ in new variables,
$\ell(F\times E)=5$ by (i), and the decomposition of length $5$ with $E$ applied last is of
type (a) in Proposition~\ref{prop:length-k1}. For Vasyunin's map $V$ and such a map $F$,
(i) gives $\ell(V\times F)\ge8$, hence $\ell(V\times F)=8$. In the situation of (i) the
minimal decompositions of $F\times G$ are the interleavings, so the largest degree of a
factor, minimized over them, is $\max(\omega_0(F),\omega_0(G))$, with $\omega_0$ as in
Section~\ref{sec:plane}. Finally $(F\times G)^{(m)}=F^{(m)}\times G^{(m)}$, so
$\lambda(F\times G)=\max(\lambda(F),\lambda(G))$, while the length is only subadditive.

Let $F=Q_l\circ\dots\circ Q_1$ with $Q_i$ changing $X_{m_i}$. We say that a second
decomposition of $F$ of length $l$ is obtained by a \emph{move} if it arises in one of
the following ways.
\begin{enumerate}[label=\textup{(M\arabic*)},leftmargin=2.6em]
\item Two adjacent factors $Q_{i+1},Q_i$ are replaced by elementary maps
$Q'_{i+1},Q'_i$ with $Q'_{i+1}\circ Q'_i=Q_{i+1}\circ Q_i$.
\item For $i<j$ with $m_i=m_j=c$ and an elementary map $P$ changing $X_c$ by a
polynomial $p$ which involves none of the $X_{m_k}$, $i<k<j$, the factors $Q_j$, $Q_k$
for $i<k<j$, and $Q_i$ are replaced by $Q_j\circ P$, $P^{-1}\circ Q_k\circ P$ and
$P^{-1}\circ Q_i$.
\end{enumerate}
In (M2), if $Q_k$ changes $X_{m_k}$ by $h$, then $P^{-1}\circ Q_k\circ P$ changes
$X_{m_k}$ by $h\circ P$, since $p$ does not involve $X_{m_k}$. If $l=\ell(F)$, no new
factor equals $\Id$, so moves map decompositions of length $\ell(F)$ to decompositions
of length $\ell(F)$. The classes of the equivalence relation generated by moves are
called \emph{move classes}.

If $\ell(F)>\#\supp(F)$, then in every decomposition of length $\ell(F)$ some index $c$
occurs twice, by Lemmas~\ref{lem:support-bound} and~\ref{lem:word}(i), and at
non-adjacent positions by Lemma~\ref{lem:word}(ii). The move (M2) with the translations
$P\colon X_c\mapsto X_c+\varepsilon$, $\varepsilon\in K$, then gives pairwise distinct
decompositions of length $\ell(F)$. So for infinite $K$ such an $F$ has infinitely many
minimal decompositions in the move class of each of them, and the number of minimal decompositions
is meaningful only up to moves. If $\ell(F)=\#\supp(F)$, (M2) does not apply, and there
are at most $k!$ minimal decompositions by Theorem~\ref{thm:min-decomp}.

Theorem~\ref{thmA} follows from Lemma~\ref{lem:support-bound},
Proposition~\ref{prop:length-spectrum}, Theorem~\ref{thm:min-decomp},
Corollary~\ref{cor:stab} and Proposition~\ref{prop:stab-bound}.

\begin{example}\label{ex:strict}
Let $\operatorname{char}K\ne2$ and $F=\bigl(X_1+X_2+(X_1+2X_2)^2,\;X_1+2X_2\bigr)$. Then
$F=Q_3\circ Q_2\circ Q_1$ with $Q_1=(X_1+X_2,X_2)$, $Q_2=(X_1,X_2+X_1)$,
$Q_3=(X_1+X_2^2,X_2)$, and
$F^{-1}=(2X_1-X_2-2X_2^2,\;X_2-X_1+X_2^2)$. The polynomial $X_1-(F^{-1})_1$ involves
$X_1$ and $X_2-(F^{-1})_2$ involves $X_2$. So no ordering is admissible, and
$\ell(F\times\Id_r)=\ell(F)=3$ for every $r\ge0$ by Corollary~\ref{cor:stab}. The
same argument applies to the linear map $M=(5X_1+3X_2,\,3X_1+2X_2)$, with
$M^{-1}=(2X_1-3X_2,\,-3X_1+5X_2)$ and
$M=(X_1+\tfrac43X_2,X_2)\circ(X_1,X_2+3X_1)\circ(X_1+\tfrac13X_2,X_2)$, so
$\ell(M\times\Id_r)=3$ when $\operatorname{char}K\notin\{2,3\}$. Here
$X_2-(M^{-1})_2=3X_1-4X_2$ involves $X_2$ only because $4\ne0$. If
$\operatorname{char}K=2$, then $M=(X_1+X_2,X_1)=(X_1,X_2+X_1)\circ(X_1+X_2,X_2)$, and
$\ell(M)=2$.
\end{example}

\section{Vasyunin's map}\label{sec:vasyunin}

\subsection{Three decompositions}\label{ssec:three}

Let
\begin{align*}
Q_1&\colon X_5\mapsto X_5+\tfrac12X_3^2, &
Q_2&\colon X_4\mapsto X_4+X_1X_5-X_2X_3-\tfrac12X_1X_3^2,\\
Q_3&\colon X_2\mapsto X_2+X_1X_3, &
Q_4&\colon X_3\mapsto X_3+X_1X_4+\tfrac12X_2^2-X_1^2X_5,
\end{align*}
each map fixing the other coordinates. Then $F=Q_4\circ Q_3\circ Q_2\circ Q_1$. Indeed
$Q_3\circ Q_2\circ Q_1=(X_1,\,X_2+X_1X_3,\,X_3,\,X_4+X_1X_5-X_2X_3,\,X_5+\tfrac12X_3^2)$,
and $Q_4$ then adds
$\tfrac12(X_2+X_1X_3)^2+X_1(X_4+X_1X_5-X_2X_3)-X_1^2(X_5+\tfrac12X_3^2)
=X_1X_4+\tfrac12X_2^2$ to $X_3$.

Let $R_1,\dots,R_7$ add, in this order, $\tfrac12X_3^2$ to $X_5$, $X_1X_3$ to $X_2$,
$\tfrac12X_1X_3^2$ to $X_4$, $-X_2X_3$ to $X_4$, $\tfrac12X_2^2$ to $X_3$, $X_1X_4$ to
$X_3$ and $X_1X_5$ to $X_4$. Then $F=R_7\circ\dots\circ R_1$. Let moreover
$C'\colon X_2\mapsto X_2+\tfrac12X_3^2$ and $C\colon X_4\mapsto X_4+X_1X_2$. Then
$C\circ C'\circ C^{-1}\circ C'^{-1}=R_3$ by Lemma~\ref{lem:commutator} below, so
\begin{equation}\label{eq:quad-decomposition}
F=R_7\circ R_6\circ R_5\circ R_4\circ C\circ C'\circ C^{-1}\circ C'^{-1}\circ R_2\circ R_1
\end{equation}
is a decomposition of length $10$ into factors of degree $2$. The factors of the
three decompositions have the degrees $(2,3,2,3)$, $(2,2,3,2,2,2,2)$ and $(2,\dots,2)$.
For $n=2$ the multidegree of the normal form is an invariant
\cite[Proposition~23]{AH_canon}. For Vasyunin's map not even the multiset of the degrees of
the factors of a minimal decomposition is determined by $F$ (Theorem~\ref{thm:V-minimal}).

\subsection{Decompositions of minimal length}\label{ssec:minimal}

Let $K$ be a field with $\operatorname{char}K\ne2$, and let $F$ be Vasyunin's map over
$K$. Let $\Qout$, $P_2$, $P_5$ and $P_4^{(c)}$, $c\in K$, be the elementary maps
\[
\begin{aligned}
\Qout&\colon X_3\mapsto X_3+X_1X_4+\tfrac12X_2^2-X_1^2X_5, &
P_2&\colon X_2\mapsto X_2+X_1X_3,\\
P_5&\colon X_5\mapsto X_5+\tfrac12X_3^2, &
P_4^{(c)}&\colon X_4\mapsto X_4+X_1X_5-X_2X_3+c\,X_1X_3^2 .
\end{aligned}
\]

\begin{theorem}\label{thm:V-minimal}
One has $\ell(F)=4$, and the decompositions of $F$ of length $4$ are exactly the six
words $F=\Qout\circ W$, where $W$ is a product of $P_2$, $P_5$ and $P_4^{(c)}$ in any
order, and $c=\varepsilon_2-\tfrac12\varepsilon_5$. Here $\varepsilon_2=1$ if $P_2$ is
applied before $P_4^{(c)}$ and $\varepsilon_2=0$ otherwise, and $\varepsilon_5$ is
defined in the same way for $P_5$.
\end{theorem}

\begin{proof}
The factor applied last in a decomposition of length $4$ is a left elementary factor of
$F$, and the inverse of $F$ shows that $X_3$ is the only direction in which $F$ has one.
Here $\supp(F)=\{2,3,4,5\}$, so $\ell(F)=4$ by Lemma~\ref{lem:support-bound} and
Section~\ref{ssec:three}, whose factors have coefficients in $\mathbb{Z}[\tfrac12]$. We
apply Theorem~\ref{thm:min-decomp} with the inverse $G=F^{-1}$. By
Section~\ref{ssec:three}, $G=Q_1^{-1}\circ Q_2^{-1}\circ Q_3^{-1}\circ Q_4^{-1}$, and
applying $Q_4^{-1},Q_3^{-1},Q_2^{-1},Q_1^{-1}$ in turn gives $G_1=X_1$ and
\begin{equation}\label{eq:inverse-compact}
\begin{gathered}
G_3=X_3-X_1X_4-\tfrac12X_2^2+X_1^2X_5,\qquad G_2=X_2-X_1G_3,\\
G_4=X_4-X_1X_5+G_2G_3+\tfrac12X_1G_3^2,\qquad G_5=X_5-\tfrac12G_3^2
\end{gathered}
\end{equation}
over every field of characteristic $\ne2$. This is the inverse computed in
\cite[Section~5.3]{AH_sn}, whose expanded form is recalled in
Appendix~\ref{ssec:inverse-data}. By the proof of Theorem~\ref{thm:min-decomp}, the factor
applied last in a decomposition of length $4$ is a left elementary factor of $F$, in a direction
$X_j$ with $j\in\supp(F)$ for which $X_j-G_j$ does not involve $X_j$. The polynomial
$X_2-G_2=X_1G_3$ contains the monomial
$X_1X_2^2$ with coefficient $-\tfrac12$. The monomial $X_1^4X_4X_5$ occurs in $G_2G_3$
with coefficient $2$ and in $\tfrac12X_1G_3^2$ with coefficient $-1$, so $X_4-G_4$
contains it with coefficient $-1$. Finally $X_5-G_5=\tfrac12G_3^2$ contains $X_1^4X_5^2$
with coefficient $\tfrac12$. These coefficients are nonzero in $K$,
whereas $X_3-G_3=X_1X_4-X_1^2X_5+\tfrac12X_2^2$. So the factor applied last changes $X_3$, and it is
$\Qout$ by Lemma~\ref{lem:left-factor}. Then
\[
\Qout^{-1}\circ F=\bigl(X_1,\;X_2+X_1X_3,\;X_3,\;X_4+X_1X_5-X_2X_3,\;
X_5+\tfrac12X_3^2\bigr).
\]
This map is of the form $M=(X_1,X_2+u_2,X_3,X_4+u_4,X_5+u_5)$ with
$u_2,u_5\in K[X_1,X_3]$ and $u_4\in K[X_1,X_2,X_3,X_5]$. Every such map is an
automorphism, with inverse
\[
M^{-1}=\bigl(X_1,\;X_2-u_2,\;X_3,\;X_4-u_4(X_1,X_2-u_2,X_3,X_5-u_5),\;X_5-u_5\bigr).
\]
So the polynomials
\[
\begin{gathered}
h_2=X_2-(M^{-1})_2=u_2,\qquad h_5=X_5-(M^{-1})_5=u_5,\\
h_4=X_4-(M^{-1})_4=u_4(X_1,X_2-u_2,X_3,X_5-u_5)
\end{gathered}
\]
do not involve $X_2$, $X_4$, $X_5$ respectively, and $h_j\ne0$ exactly when $u_j\ne0$.
By Lemma~\ref{lem:left-factor}, for each $j$ with $u_j\ne0$ the map $P$ changing $X_j$
by $h_j$ is a left elementary factor of $M$, and $P^{-1}\circ M$ agrees with $M$ except
in the $j$-th coordinate, which is $X_j+u_j-h_j(M)$. Since $M$ fixes $X_1$ and $X_3$,
$h_2(M)=u_2$ and $h_5(M)=u_5$. Since $M$ maps $X_2-u_2$ to $X_2$ and $X_5-u_5$ to
$X_5$, $h_4(M)=u_4$. Hence $P^{-1}\circ M$ is the map of the same form with $u_j$
replaced by $0$. So all six orderings of $\{2,4,5\}$
are admissible. When $P_4^{(c)}$ is applied, $X_2$ has been replaced by $X_2+X_1X_3$ if
$\varepsilon_2=1$ and $X_5$ by $X_5+\tfrac12X_3^2$ if $\varepsilon_5=1$, so the tail
of $P_4^{(c)}$ becomes
$X_1X_5-X_2X_3+(c-\varepsilon_2+\tfrac12\varepsilon_5)X_1X_3^2$. This matches the
fourth coordinate of $\Qout^{-1}\circ F$ exactly when $c$ is as stated.
\end{proof}

Explicitly, listing the factors of $W$ in the order of application, $c=0$ for
$(P_4^{(c)},P_2,P_5)$ and $(P_4^{(c)},P_5,P_2)$, $c=1$ for $(P_2,P_4^{(c)},P_5)$,
$c=-\tfrac12$ for $(P_5,P_4^{(c)},P_2)$, and $c=\tfrac12$ for $(P_2,P_5,P_4^{(c)})$ and
$(P_5,P_2,P_4^{(c)})$. For the two orders with $c=0$ the factor
degrees are $\{2,2,2,3\}$. For the other four orders they are $\{2,2,3,3\}$. The word
$Q_4\circ Q_3\circ Q_2\circ Q_1$ of Section~\ref{ssec:three} corresponds to the order
$P_5,P_4^{(-1/2)},P_2$.

\subsection{Elementary factors of degree at most two}\label{ssec:ell2}

In this subsection, as in Section~\ref{ssec:minimal}, $K$ is a field with
$\operatorname{char}K\ne2$. Let $\ell_2(F)$ be the least length of a decomposition of
$F$ into elementary maps of degree at most $2$. By \eqref{eq:quad-decomposition}, $\ell_2(F)\le10$.

\begin{proposition}\label{prop:ell2}
$\ell_2(F)=6$. With $S\colon X_4\mapsto X_4-X_1X_5$ and
$P_3'\colon X_3\mapsto X_3+X_1X_4+\tfrac12X_2^2$,
\begin{equation}\label{eq:six-quadratic}
F=S^{-1}\circ P_3'\circ S\circ P_2\circ P_5\circ P_4^{(0)} .
\end{equation}
\end{proposition}

\begin{proof}
\emph{Upper bound.} The map $S$ replaces $X_4$ by $X_4-X_1X_5$, then $P_3'$ adds
$X_1(X_4-X_1X_5)+\tfrac12X_2^2$ to $X_3$, and $S^{-1}$ restores $X_4$. So
$\Qout=S^{-1}\circ P_3'\circ S$, and \eqref{eq:six-quadratic} follows from
Theorem~\ref{thm:V-minimal} with $P_4^{(0)}$ applied first.

\emph{Lower bound.} Every decomposition of length $4$ contains $\Qout$, of degree $3$. It
remains to exclude length $5$. In a word of length $5$ one direction is repeated, and
for every possible position of the repetition a monomial of degree at least $3$ appears
in a tail which must have degree at most $2$. Let $F=Q_5\circ\dots\circ Q_1$ with factors of degree at
most $2$, $Q_i$ changing $X_{m_i}$. By Lemma~\ref{lem:support-bound} the indices
$m_1,\dots,m_5$ contain $2,3,4,5$. If they are pairwise distinct, then one of them is
$1$, and exactly one factor changes $X_1$, which contradicts $F_1=X_1$ by
Lemma~\ref{lem:word}(i). Hence exactly one index $c\in\{2,3,4,5\}$
occurs twice among the $m_i$, and $1$ does not occur.

\emph{Step 1.} $m_5=c$. Otherwise $Q_4\circ\dots\circ Q_1$ fixes $X_{m_5}$, so $Q_5$ is
a left elementary factor of $F$. By the proof of Theorem~\ref{thm:V-minimal}, $X_3$ is the
only direction in which $F$ has a left elementary factor, and this factor is $\Qout$, of
degree $3$.

\emph{Step 2.} Let $Q_5$ change $X_c$ by $h$, with $\deg h\le2$, and let $j=m_4$. The
indices $m_1,\dots,m_4$ are $2,3,4,5$ in some order, so $Q_3\circ Q_2\circ Q_1$ fixes
$X_j$, and $Q_4$ is a left elementary factor of $Q_5^{-1}\circ F=Q_4\circ\dots\circ Q_1$
in the direction $X_j$. By Lemma~\ref{lem:left-factor}, applied with
$(Q_5^{-1}\circ F)^{-1}=F^{-1}\circ Q_5$, the tail of $Q_4$ is
$X_j-(F^{-1}\circ Q_5)_j$. It must not involve $X_j$ and must have degree at most $2$.
For the $14$ pairs $(c,j)$ other than
$(2,3)$ and $(4,3)$, Table~\ref{tab:monomials} in Appendix~\ref{app:record} lists a monomial
of degree at least $3$ whose coefficient in $X_j-(F^{-1}\circ Q_5)_j$ is a nonzero
constant, independent of $h$. These constants are $\pm1$, $\pm\tfrac12$ and $\tfrac18$,
nonzero in $K$. Each entry of the table is obtained from Lemma~\ref{lem:monomial-coeff},
which reduces it to an inspection of the terms of $(F^{-1})_j$ involving $X_c$. For $(c,j)=(2,3)$, $h\in K[X_1,X_3,X_4,X_5]$ and
\[
X_3-(F^{-1}\circ Q_5)_3=X_1X_4-X_1^2X_5+\tfrac12(X_2+h)^2 .
\]
Its homogeneous component of degree $4$ is $\frac12h_2^2$, where $h_2$ is the quadratic part of $h$.
So $h_2^2=0$, hence $h_2=0$, and then the coefficient of $X_1^2X_5$ is $-1$. For
$(c,j)=(4,3)$, $h\in K[X_1,X_2,X_3,X_5]$ and
\[
X_3-(F^{-1}\circ Q_5)_3=X_1X_4+X_1h-X_1^2X_5+\tfrac12X_2^2 .
\]
Its cubic component $X_1(h_2-X_1X_5)$ vanishes, so $h_2=X_1X_5$, and $X_3$ does not occur,
so the linear part of $h$ does not involve $X_3$. Hence
$h=X_1X_5+a_0+a_1X_1+a_2X_2+a_5X_5$, $a_i\in K$, and after extraction of $Q_4$ the
remaining map is
\[
G'=\bigl(X_1,\;X_2+X_1X_3,\;X_3,\;X_4+\tau,\;X_5+\tfrac12X_3^2\bigr),
\]
where
\[
\tau=-\tfrac12X_1X_3^2-X_2X_3-a_2X_1X_3-\tfrac12a_5X_3^2-a_1X_1-a_2X_2-a_5X_5-a_0 .
\]

\emph{Step 3.} The word $Q_3\circ Q_2\circ Q_1$ for $G'$ changes each of $X_2,X_4,X_5$ exactly once.
The map $G'$ is of the form $M$ in the proof of Theorem~\ref{thm:V-minimal}, with
$u_2=X_1X_3$, $u_5=\tfrac12X_3^2$ and $u_4=\tau$. By that proof, in every order of
extraction the factors changing $X_2$ and $X_5$ add $X_1X_3$ and $\tfrac12X_3^2$, and the factor
changing $X_4$ has a tail $t\in K[X_1,X_2,X_3,X_5]$ of degree at most $2$ with
$t(X+q)=\tau$, where $q=(0,\,\varepsilon_2X_1X_3,\,0,\,0,\,\tfrac{\varepsilon_5}2X_3^2)$ and
$\varepsilon_2,\varepsilon_5\in\{0,1\}$ record whether these factors are applied earlier. Write
$t=t_0+t_1+t_2$ with $t_i$ homogeneous of degree $i$. Then
\[
t(X+q)=t(X)+\Bigl(t_1(q)+\sum_iq_i\,\partial_{X_i}t_2\Bigr)+t_2(q),
\]
where $t_1(q)$ is a quadratic form in $X_1,X_3$, the sum is homogeneous of degree $3$
and $t_2(q)$ is homogeneous of degree $4$. Hence every monomial of degree $2$ which
involves $X_2$ or $X_5$ has the same coefficient in $t(X+q)$ as in $t_2$, and the cubic
part of $t(X+q)$ is
$\varepsilon_2X_1X_3\,\partial_{X_2}t_2+\tfrac{\varepsilon_5}2X_3^2\,\partial_{X_5}t_2$.
Write $t_2=\sum_{i\le j}a_{ij}X_iX_j$ with $i,j\in\{1,2,3,5\}$. In $\tau$ the monomials
$X_2^2$, $X_2X_5$, $X_5^2$ and $X_1X_5$ do not occur, and $X_2X_3$ has the coefficient
$-1$. Hence $a_{22}=a_{25}=a_{55}=a_{15}=0$ and $a_{23}=-1$, so that
\[
\partial_{X_2}t_2=a_{12}X_1-X_3,\qquad \partial_{X_5}t_2=a_{35}X_3 .
\]
Put $\alpha=a_{12}$ and $\kappa=a_{35}$. The cubic part of
$t(X+q)$ is then $\varepsilon_2X_1X_3(\alpha X_1-X_3)+\tfrac{\varepsilon_5}2\kappa X_3^3$, whose
coefficient of $X_1X_3^2$ is $-\varepsilon_2\in\{0,-1\}$. In $\tau$ it is $-\tfrac12$,
which differs from $0$ and from $-1$ since $\operatorname{char}K\ne2$. This is a
contradiction.
\end{proof}

A cubic elementary map is never a product of two elementary maps of degree at most $2$.
Two factors changing the same coordinate give a tail of degree at most $2$. If the two
factors change different coordinates, their composition changes two coordinates and
hence is not elementary. The
decomposition \eqref{eq:six-quadratic} writes $\Qout$ as a conjugate of a quadratic
elementary map, and \eqref{eq:quad-decomposition} writes $R_3$ as a commutator.

\subsection{Triangular subgroups}\label{ssec:triangular}

As before, $K$ is a field with $\operatorname{char}K\ne2$.

\begin{corollary}\label{cor:no-borel}
$F$ lies in no $J_\sigma(K,5)$, and no decomposition of $F$ into elementary maps has
all its factors in one $J^1_\sigma(K,5)$.
\end{corollary}

\begin{proof}
The hypothesis of Lemma~\ref{lem:dag} holds, since $H_1=0$ and $H_2,\dots,H_5$ are
quadratic forms. $D_F$ contains the cycle $2\to3\to2$, from $X_3$ in $H_2$ and $X_2$ in
$H_3$, so $F$ lies in no $J_\sigma(K,5)$ by Lemma~\ref{lem:dag}. The second statement follows, since $J^1_\sigma$ is a
group.
\end{proof}

Thus the equality $\ell(F)=\#\supp(F)$ does not imply that $F$ lies in some
$J_\sigma(K,5)$, and
Theorem~\ref{thmB} follows from Theorem~\ref{thm:V-minimal},
Proposition~\ref{prop:ell2} and Corollary~\ref{cor:no-borel}.

\subsection{Degrees of the inverse}\label{ssec:inverse}

A decomposition $F=Q_l\circ\dots\circ Q_1$ gives $F^{-1}=Q_1^{-1}\circ\dots\circ
Q_l^{-1}$, hence $\deg F^{-1}\le\prod_i\deg Q_i$. For the three decompositions of
Section~\ref{ssec:three} this gives $36$, $192$ and $1024$, while $\deg F^{-1}=7$, within
the bound $(\deg F)^{n-1}=16$ of \cite[Theorem~1.5]{BCW}. The following recursion follows
the degrees of the individual coordinates.

\begin{proposition}\label{prop:degree-recursion}
Let $F=Q_l\circ\dots\circ Q_1$, where $Q_j$ changes $X_{m_j}$ by a tail $t_j$. For
$D\in\mathbb{Z}_{\ge0}^n$ let $\Lambda_{t}(D)=\max\sum_ka_kD_k$ over the monomials $X^a$
of $t$. Let $D^{(l)}_{m_l}=\max(1,\deg t_l)$ and $D^{(l)}_k=1$ for $k\ne m_l$, and for
$j=l-1,\dots,1$ let $D^{(j)}$ agree with $D^{(j+1)}$ except
$D^{(j)}_{m_j}=\max\bigl(D^{(j+1)}_{m_j},\Lambda_{t_j}(D^{(j+1)})\bigr)$. Then
$\deg(F^{-1})_k\le D^{(1)}_k$ for all $k$.
\end{proposition}

\begin{proof}
Put $G_j=Q_j^{-1}\circ\dots\circ Q_l^{-1}$. Then $(G_j)_k=(G_{j+1})_k$ for $k\ne m_j$ and
$(G_j)_{m_j}=(G_{j+1})_{m_j}-t_j(G_{j+1})$, and reverse induction on $j$ gives
$\deg(G_j)_k\le D^{(j)}_k$.
\end{proof}

For the decomposition $Q_4\circ Q_3\circ Q_2\circ Q_1$, for the seven-factor
decomposition, for \eqref{eq:six-quadratic} and for the six decompositions of
Theorem~\ref{thm:V-minimal} the recursion gives $D^{(1)}=(1,4,3,7,6)$, the degrees of
the coordinates of $F^{-1}$. For \eqref{eq:quad-decomposition} it gives $(1,6,3,7,6)$.
The recursion depends only on the monomials of the tails and does not detect
cancellation.

\section{Iterates of Vasyunin's map}\label{sec:iterates}

Let $\operatorname{char}K\ne2$ and $F_K^{(m)}=(S_1^{(m)},\dots,S_5^{(m)})$. Then
$S_1^{(m)}=X_1$ and
\begin{equation}\label{eq:outer-rec}
\begin{aligned}
S_2^{(m+1)}&=S_2^{(m)}+X_1S_3^{(m)}, &
S_3^{(m+1)}&=S_3^{(m)}+X_1S_4^{(m)}+\tfrac12\bigl(S_2^{(m)}\bigr)^2,\\
S_4^{(m+1)}&=S_4^{(m)}+X_1S_5^{(m)}-S_2^{(m)}S_3^{(m)}, &
S_5^{(m+1)}&=S_5^{(m)}+\tfrac12\bigl(S_3^{(m)}\bigr)^2 .
\end{aligned}
\end{equation}
For $P\ne0$ let $T(P)$ be its homogeneous component of highest degree. Then
$T(PQ)=T(P)T(Q)$, and $T(P+Q)=T(P)$ if $\deg P>\deg Q$. Let $s_0=1$,
\[
s_{2j+1}=2^{j+2}-2\ (j\ge0),\qquad s_{2j}=3\cdot2^{j}-2\ (j\ge1),
\]
so that $(s_m)=(1,2,4,6,10,14,22,\dots)$ is strictly increasing and $s_{m+1}=2+2s_{m-1}$
for $m\ge1$. Let $t_1=2$, $t_{2j}=7\cdot2^{j-1}-3$ and $t_{2j+1}=5\cdot2^{j}-3$ for
$j\ge1$.

\begin{theorem}\label{thm:exact-degrees}
For every $m\ge1$,
\begin{equation}\label{eq:top-forms}
T\bigl(S_2^{(m)}\bigr)=X_1\,T\bigl(S_3^{(m-1)}\bigr),\quad
T\bigl(S_3^{(m+1)}\bigr)=\tfrac12X_1^2\,T\bigl(S_3^{(m-1)}\bigr)^2,\quad
T\bigl(S_5^{(m)}\bigr)=\tfrac12\,T\bigl(S_3^{(m-1)}\bigr)^2,
\end{equation}
with $T(S_3^{(0)})=X_3$ and $T(S_3^{(1)})=X_1X_4+\tfrac12X_2^2$. Consequently
$\deg S_3^{(m)}=s_m$, $\deg S_2^{(m)}=1+s_{m-1}$, $\deg S_5^{(m)}=2s_{m-1}$,
$\deg S_4^{(m)}\le t_m$, and $\deg F_K^{(m)}=2s_{m-1}$.
\end{theorem}

The proof follows the top forms, that is, it computes in the associated graded ring of
$K[X]$ with respect to the total degree, in which the top forms of the iterates satisfy
the recursion \eqref{eq:top-forms}. In each sum of \eqref{eq:outer-rec} defining
$S_2^{(m+1)}$, $S_3^{(m+1)}$ and $S_5^{(m+1)}$ one summand has strictly larger degree than
the others, namely $X_1S_3^{(m)}$, $\tfrac12(S_2^{(m)})^2$ and $\tfrac12(S_3^{(m)})^2$,
so its top form is the top form of the sum. For $S_4^{(m+1)}$ only an upper bound is
needed. The comparison of degrees reduces to inequalities between the numbers $s_m$ and
$t_m$.

\begin{proof}
We prove by induction on $m$ the statement $(\mathrm{A}_m)$ that \eqref{eq:top-forms} holds
for $S_2^{(m)}$ and $S_5^{(m)}$, $\deg S_3^{(m)}=s_m$, and $\deg S_4^{(m)}\le t_m$. For
$m=1$ this is read off $F_K$. Assume $(\mathrm{A}_k)$ for $k\le m$.

In $S_2^{(m+1)}$, $\deg S_2^{(m)}=1+s_{m-1}<1+s_m=\deg X_1S_3^{(m)}$, and in
$S_5^{(m+1)}$, $\deg S_5^{(m)}=2s_{m-1}<2s_m$. This gives the top forms of
$S_2^{(m+1)}$ and $S_5^{(m+1)}$. In $S_3^{(m+1)}$ the summand $\tfrac12(S_2^{(m)})^2$ has
degree $2+2s_{m-1}=s_{m+1}$ and top form $\tfrac12X_1^2T(S_3^{(m-1)})^2\ne0$, since
$2\in K^*$. The other summands have degrees $s_m<s_{m+1}$ and at most $1+t_m$, and
$1+t_m<s_{m+1}$. Indeed $1+t_1=3<4$, $1+t_{2j}=7\cdot2^{j-1}-2<8\cdot2^{j-1}-2$ and
$1+t_{2j+1}=5\cdot2^{j}-2<6\cdot2^{j}-2$. In $S_4^{(m+1)}$ the summands have degrees at
most $t_m$, $1+2s_{m-1}$ and $1+s_{m-1}+s_m$, the last being the largest of the second
and third. For $m=1$ this gives $4=t_2$, for $m=2j$ one has
$1+s_{2j-1}+s_{2j}=t_{2j+1}\ge t_{2j}$, and for $m=2j+1$ one has
$1+s_{2j}+s_{2j+1}=t_{2j+2}\ge t_{2j+1}$. This proves $(\mathrm{A}_{m+1})$.

Finally $2s_{m-1}$ is the largest of the five degrees. For $m\ge2$,
$s_m=2+2s_{m-2}\le2s_{m-1}$ and $1+s_{m-1}<2s_{m-1}$, and
$t_{2j}\le8\cdot2^{j-1}-4=2s_{2j-1}$, with equality only for $j=1$, and
$t_{2j+1}<6\cdot2^{j}-4=2s_{2j}$.
\end{proof}

In terms of $a_m=\deg S_2^{(m)}$ and $b_m=\deg S_3^{(m)}$ the theorem gives
$a_{m+1}=b_m+1$ and $b_{m+1}=2a_m$ for $m\ge1$. This is an affine recursion whose linear
part has the matrix $\left(\begin{smallmatrix}0&1\\2&0\end{smallmatrix}\right)$, with the
eigenvalues $\pm\sqrt2$, which explains the value of the dynamical degree.

Unwinding \eqref{eq:top-forms} gives, for $j\ge1$,
$T(S_3^{(2j)})=c_jX_1^{2^{j+1}-2}X_3^{2^j}$ and
$T(S_3^{(2j+1)})=c_j'X_1^{2^{j+1}-2}(X_1X_4+\tfrac12X_2^2)^{2^j}$ with $c_j,c_j'$ powers
of $\tfrac12$.

\begin{lemma}\label{lem:bounded}
\textup{(i)} If $G\in J_\sigma(K,n)$, then $\deg G^{(m)}$ is bounded in $m$.
\textup{(ii)} If $G=\varphi^{-1}\circ G'\circ\varphi$ with $\varphi\in\Aut_KK[X]^n$, then
$\deg G^{(m)}\le\deg\varphi^{-1}\cdot\deg G'^{(m)}\cdot\deg\varphi$.
\end{lemma}

\begin{proof}
(i) We may take $\sigma=\mathrm{id}$. Then
$G^{(m)}_i=a_iG^{(m-1)}_i+f_i(G^{(m-1)}_1,\dots,G^{(m-1)}_{i-1})$. If
$\deg G^{(m)}_k\le d_k$ for $k<i$ and all $m$, induction on $m$ gives
$\deg G^{(m)}_i\le\max(1,\deg f_i\cdot\max_{k<i}d_k)$. (ii) is clear.
\end{proof}

\begin{corollary}\label{cor:shifted}
Let $L=\mathrm{diag}(\lambda_1,\dots,\lambda_5)$, $\lambda_i\in K^*$. Then
$\deg(L\circ F_K)^{(m)}=\deg F_K^{(m)}$ for all $m$, and
$\lim_m(\deg F_K^{(m)})^{1/m}=\sqrt2$. The map $L\circ F_K$ has infinite order, is not
locally finite, and is not conjugate in $\Aut_KK[X_1,\dots,X_5]^5$ to any automorphism whose iterates
have bounded degree, such as an element of some $J_\sigma(K,5)$ or a linear map. In
particular $F_K$ is not Pascal finite.
\end{corollary}

\begin{proof}
Write $(L\circ F_K)^{(m)}=(S_1^{(m)},\dots,S_5^{(m)})$. Since
$(L\circ F_K)^{(m+1)}=L\circ F_K\circ(L\circ F_K)^{(m)}$, one has
$S_i^{(m+1)}=\lambda_i\bigl(S_i^{(m)}+H_i(S^{(m)})\bigr)$, so $S_1^{(m)}=\lambda_1^mX_1$ and
the iterates satisfy \eqref{eq:outer-rec} with every summand multiplied by a nonzero
constant. This changes the top forms in
\eqref{eq:top-forms} by nonzero constant factors and does not change the degrees of the
summands. The proof of Theorem~\ref{thm:exact-degrees} compares only these degrees and
uses only that the top forms are nonzero, so it applies. The limit follows from the closed formulas.
A map of finite order, or a locally finite map, has iterates of bounded degree, so
$L\circ F_K$ has neither property. If $\Delta_{F_K}^m(X)=0$, then
$\sum_{i=0}^m(-1)^{m-i}\binom mi F_K^{(i+k)}=0$ for every $k\ge0$, by composition with
$F_K^{(k)}$, so every iterate has degree at most $\max_{i<m}\deg F_K^{(i)}$. Hence $F_K$
is not Pascal finite. The rest follows from Lemma~\ref{lem:bounded}, since a linear map
has iterates of degree $1$.
\end{proof}

\begin{corollary}\label{cor:modp}
Let $p$ be an odd prime and $F_p$ the reduction of $F$ modulo $p$. Then $F_p$ is a tame
element of infinite order of $\Aut_{\FF_p}\FF_p[X_1,\dots,X_5]^5$, while its evaluation map
$x\mapsto F_p(x)$ is a permutation of the finite set $\FF_p^{\,5}$ and has finite order.
Its $J_H$ is nilpotent and not strongly nilpotent.
\end{corollary}

\begin{proof}
The decompositions of Section~\ref{ssec:three} have coefficients in $\mathbb{Z}[\tfrac12]$
and reduce modulo $p$. Infinite order is Corollary~\ref{cor:shifted}. The identity
$J_H^5=0$ holds over $\mathbb{Z}[\tfrac12]$, and $J_H(e_3)J_H(e_2)J_H(e_1)^3=E_{55}$ is
nonzero modulo $p$.
\end{proof}

This proves Theorem~\ref{thmD}.

\begin{remark}\label{rem:modp-contrast}
In characteristic $p$ an automorphism $G$ is Pascal finite if and only if its order is a
power of $p$, since $\Delta_G^{p^e}=\sigma_G^{p^e}-\Id$ \cite[Lemma~5]{AH_JvdK}. The Nagata
automorphism satisfies $\Delta^3(X)=0$ over $\mathbb{Z}$ \cite[Example~2.1]{ABCH2}, so its
reductions modulo odd primes have order $p$. The
reductions of Vasyunin's map have infinite order.
\end{remark}

\begin{remark}\label{rem:context}
The dynamical degree of a plane automorphism is $1$ or an integer $\ge2$ \cite{FM,F}.
Blanc and van Santen determine the dynamical degrees of automorphisms of $K[X_1,X_2,X_3]$
of degree $2$ over every field as $\{1,\sqrt2,\tfrac{1+\sqrt5}2,2\}$
\cite[Theorem~2]{BvS}, extending the case $K=\mathbb{C}$ due to Maegawa
\cite[Theorem~3.1]{Mae}. In dimension five, Shao and Sun show that the dynamical degree
of the maps $S_2$ of Section~\ref{ssec:vasyunin-def} is an algebraic integer of degree at
most $4$ \cite[Proposition~4.7]{SS}, in accordance with the conjecture of Favre that dynamical
degrees of endomorphisms of $K[X_1,\dots,X_n]$ are algebraic integers of degree at most $n$,
recalled in \cite[Conjecture~1.1]{SS}. For the member of this family which is linearly conjugate to
Vasyunin's map, Corollary~\ref{cor:shifted} gives the value $\sqrt2$. Maubach
and Poloni show that for the Nagata automorphism $N$ and $s\in K^*\setminus\{\pm1\}$ the
map $sN$ is conjugate to $sI$ \cite[Theorem~3.3]{MP}. By Corollary~\ref{cor:shifted} no
diagonal shift of Vasyunin's map is linearizable.
\end{remark}

Let $\operatorname{char}K=0$. For Vasyunin's map $\ell(F)=4$, while $\supp(F^{(m)})=\{2,3,4,5\}$ and
$\deg F^{(m)}=2s_{m-1}$ for all $m\ge1$. Indeed, $S_2^{(m)}$, $S_3^{(m)}$ and $S_5^{(m)}$
have degree at least $2$ by Theorem~\ref{thm:exact-degrees}, and putting $X_1=0$ in
\eqref{eq:outer-rec} gives $S_2^{(m)}=X_2$, $S_3^{(m)}=X_3+\tfrac m2X_2^2$ and
$S_4^{(m)}=X_4-mX_2X_3-\tfrac{m(m-1)}4X_2^3$, so $S_4^{(m)}\ne X_4$.

\begin{proposition}\label{prop:iterate-left}
Let $\operatorname{char}K=0$ and $m\ge2$. Then $F^{(m)}$ has no left elementary factor,
and $5\le\ell(F^{(m)})\le4m$.
\end{proposition}

\begin{proof}
The upper bound follows from Theorem~\ref{thmB}. Let $G^{(m)}=(F^{(m)})^{-1}$. Since
$F^{(k+1)}=F^{(k)}+H\circ F^{(k)}$, one has $F^{(m)}=X+\sum_{k=0}^{m-1}H\circ F^{(k)}$,
and composition with $G^{(m)}$ gives
\begin{equation}\label{eq:inverse-iterates}
G^{(m)}=X-\sum_{i=1}^mH\circ G^{(i)} .
\end{equation}
All these maps fix $X_1$. For $P\in K[X]$ write $P=\sum_aX_1^aP_{[a]}$ with
$P_{[a]}\in K[X_2,\dots,X_5]$, and let $G^{(i)}=(X_1,G_2^{(i)},\dots,G_5^{(i)})$. By
\eqref{eq:inverse-iterates} the polynomials $h_j=X_j-G_j^{(m)}$ are
\begin{align*}
h_2&=X_1\sum_{i=1}^mG_3^{(i)},&
h_3&=\sum_{i=1}^m\bigl(X_1G_4^{(i)}+\tfrac12(G_2^{(i)})^2\bigr),\\
h_4&=\sum_{i=1}^m\bigl(X_1G_5^{(i)}-G_2^{(i)}G_3^{(i)}\bigr),&
h_5&=\tfrac12\sum_{i=1}^m(G_3^{(i)})^2 .
\end{align*}
The same formulas with $i$ in place of $m$ determine the parts of $G^{(i)}$ of each
degree in $X_1$ recursively. In degree $0$,
\begin{align*}
G_{2[0]}^{(i)}&=X_2,\qquad G_{3[0]}^{(i)}=X_3-\tfrac i2X_2^2,\qquad
G_{4[0]}^{(i)}=X_4+iX_2X_3-\tfrac{i(i+1)}4X_2^3,\\
G_{5[0]}^{(i)}&=X_5-\tfrac12\sum_{l=1}^i\bigl(X_3-\tfrac l2X_2^2\bigr)^2 .
\end{align*}
So $G_{4[0]}^{(i)}$ does not involve $X_5$, and $G_{5[0]}^{(i)}$ does not involve $X_4$
and has the coefficients $1$ at $X_5$ and $-\tfrac i2$ at $X_3^2$. In degree $1$,
\begin{align*}
G_{2[1]}^{(i)}&=-\sum_{l=1}^iG_{3[0]}^{(l)}=-iX_3+\tfrac{i(i+1)}4X_2^2,\qquad
G_{3[1]}^{(i)}=-\sum_{l=1}^i\bigl(G_{4[0]}^{(l)}+X_2G_{2[1]}^{(l)}\bigr),\\
G_{4[1]}^{(i)}&=-\sum_{l=1}^i\bigl(G_{5[0]}^{(l)}-X_2G_{3[1]}^{(l)}-G_{2[1]}^{(l)}G_{3[0]}^{(l)}\bigr).
\end{align*}
Hence $G_{3[1]}^{(i)}$ does not involve $X_5$ and has the coefficient $-i$ at $X_4$, and
$G_{4[1]}^{(i)}$ has the coefficient $-\sum_{l\le i}1=-i$ at $X_5$ and
$-\sum_{l\le i}(-\tfrac l2+l)=-\tfrac{i(i+1)}4$ at $X_3^2$. In degree $2$, the
coefficient of $X_5$ in
$G_{3[2]}^{(i)}=-\sum_{l\le i}\bigl(G_{4[1]}^{(l)}+\tfrac12(G_{2[1]}^{(l)})^2+X_2G_{2[2]}^{(l)}\bigr)$
is $\sum_{l\le i}l=\tfrac{i(i+1)}2$.

These formulas give the following coefficients of monomials involving $X_j$ in $h_j$.
In $h_2$ the coefficient of $X_1X_2^2$ is $-\sum_i\tfrac i2=-\tfrac{m(m+1)}4$. In $h_4$
the coefficient of $X_1X_2X_4$ comes only from $-X_2G_{3[1]}^{(i)}$ and equals
$\sum_ii=\tfrac{m(m+1)}2$. In $h_3$ the coefficient of $X_1^2X_3^2$ comes from
$G_{4[1]}^{(i)}$ and from $\tfrac12(G_{2[1]}^{(i)})^2$ and equals
$\sum_i\bigl(-\tfrac{i(i+1)}4+\tfrac{i^2}2\bigr)=\tfrac{(m-1)m(m+1)}{12}$. In $h_5$ the
coefficient of $X_1^2X_3X_5$ comes only from $G_{3[0]}^{(i)}G_{3[2]}^{(i)}$, since
$G_{3[1]}^{(i)}$ does not involve $X_5$, and equals
$\sum_i\tfrac{i(i+1)}2=\tfrac{m(m+1)(m+2)}6$. For $m\ge2$ these four numbers are nonzero,
so no $h_j$ lies in $K[X_i:i\ne j]$, and $F^{(m)}$ has no left elementary factor by
Lemma~\ref{lem:left-factor}. Since $\#\supp(F^{(m)})=4$, Theorem~\ref{thm:min-decomp}
gives $\ell(F^{(m)})\ge5$.
\end{proof}

For $m=1$ the coefficient $\tfrac{(m-1)m(m+1)}{12}$ vanishes, in accordance with the
left elementary factor $\Qout$ of $F$ in the direction $X_3$.

\section{Maps with nilpotent Jacobian in dimension three}\label{sec:heyan}

Let $\operatorname{char}K=0$. He and Yan classify the maps
$H=(u(X_1,X_2),\,v(X_1,X_2,X_3),\,h(X_1,X_2,X_3))$ with $H(0)=0$, linearly
independent components and nilpotent $J_H$ for
$\deg_{X_3}v\ge2\deg_{X_3}h$
\cite[Theorems~2.3, 3.2]{HY}, and, up to conjugation by an element of
$\mathrm{GL}_3(K)$, for $\deg_{X_3}h\le3$ \cite[Theorem~3.3]{HY}. In the normal forms
obtained in these two cases, after the linear change of coordinates in the second one,
$H$ is of the form
\begin{equation}\label{eq:heyan}
H_1=g(W),\qquad H_2=v_1X_3-a^{-1}b'(X_1)\,g(W)-v_1l_2X_1,\qquad
H_3=c_0\,g(W)^2+l_2\,g(W),
\end{equation}
where $W=aX_2+b(X_1)$, $b(X_1)=v_1c_0aX_1^2+l_1X_1+\tilde l_2$, $a,v_1,c_0\in K^*$,
$l_1,l_2,\tilde l_2\in K$, and $g\in K[t]$ satisfies $g(0)=0$ and $d:=\deg g\ge1$
\cite[Theorem~2.3]{HY}. Then $J_H$ is nilpotent, and $F=X+H$ satisfies $\det J_F=1$.
Since $W(0)=\tilde l_2$, one has $H_1(0)=g(\tilde l_2)$, and $H(0)=0$ holds if and only if
$g(\tilde l_2)=0$. The results of this section do not use this condition, and we work with
the family \eqref{eq:heyan} without it.
The case $a=1$, $l_2=\tilde l_2=0$ is the normal form obtained by Chamberland and van den
Essen for the maps $(u(X_1,X_2),\,v(X_1,X_2,X_3),\,h(u,v))$ with $H(0)=0$, nilpotent
$J_H$ and linearly independent components, up to a linear change of coordinates
\cite[Theorem~2.1, Proposition~3.3]{CE}. It contains the counterexample of van den Essen
to the dependence problem in dimension three,
$H=\bigl(X_2-X_1^2,\,X_3+2X_1(X_2-X_1^2),\,-(X_2-X_1^2)^2\bigr)$ \cite[Section~1]{CE},
which is the case $v_1=1$, $c_0=-1$, $l_1=0$ and $g(t)=t$.

Wang and Yan remove the restrictions on the degrees in $X_3$. Let
\[
N_P=\bigl(P(X_2+X_1^2),\;X_3-2X_1P(X_2+X_1^2),\;P(X_2+X_1^2)^2\bigr),\qquad
P\in K[t]\setminus K,\ P(0)=0 .
\]
A map $H\in K[X_1,X_2,X_3]^3$ with $H(0)=0$ and linearly independent components has
nilpotent $J_H$ if and only if $T\circ H\circ T^{-1}=N_P$ for some $T\in\mathrm{GL}_3(K)$ and
some $P$ \cite[Theorem~4.1]{WY_nilp}. This answers \cite[Problem~4.1]{HY}, see
\cite[Remark~4.3]{WY_nilp}. The map $N_P$ is the case $a=v_1=c_0=1$, $b=X_1^2$,
$l_1=l_2=\tilde l_2=0$, $g=P$ of \eqref{eq:heyan}, and $\deg N_P=4\deg P$
\cite[Section~6]{WY_nilp}. Moreover $X+H$ is tame for every $H$ with nilpotent $J_H$ in
dimension three, and $X+sN_P$ has an explicit factorization into four elementary maps over
$K[s]$ \cite[Theorem~4.6, (4.20)]{WY_nilp}.

Since $\ell$ is not invariant under linear
conjugation (Remark~\ref{rem:degrees}), we work with the explicit form \eqref{eq:heyan}. Throughout this section
$F=X+H$ with $H$ as in \eqref{eq:heyan}, and $\theta=av_1$. The proofs below use only the
formulas \eqref{eq:heyan}. They rest on Lemma~\ref{lem:heyan-conj}, by which $F$ is
conjugate to a skew product over a plane automorphism of H\'enon type.

\begin{lemma}\label{lem:heyan-conj}
Let $\varphi=(X_1,\,W,\,X_3-l_2X_1)$,
\[
\Psi_1=\bigl(X_1+g(X_2),\;X_2,\;X_3+c_0\,g(X_2)^2\bigr),\qquad
\Psi_2=\bigl(X_1,\;X_2+\theta X_3,\;X_3\bigr).
\]
Then $\varphi$ is a triangular automorphism of degree $2$, with inverse
$\varphi^{-1}=(X_1,\,a^{-1}(X_2-b(X_1)),\,X_3+l_2X_1)$ of degree $2$, and
$\varphi\circ F=\Psi_2\circ\Psi_1\circ\varphi$.
\end{lemma}

\begin{proof}
Put $S=X_3-l_2X_1$, so that $\varphi=(X_1,W,S)$. The coordinates of $\varphi\circ F$
are $X_1\circ F$, $W\circ F$ and $S\circ F$. First $X_1\circ F=X_1+g(W)$ and
$S\circ F=F_3-l_2F_1=S+c_0\,g(W)^2$. Since $b$ is quadratic with leading coefficient
$v_1c_0a$,
\[
b\bigl(X_1+g(W)\bigr)=b(X_1)+b'(X_1)\,g(W)+v_1c_0a\,g(W)^2,
\]
and therefore
\[
W\circ F=aF_2+b(F_1)=W+\theta\bigl(S+c_0\,g(W)^2\bigr).
\]
On the other hand $\Psi_1\circ\varphi=(X_1+g(W),\,W,\,S+c_0g(W)^2)$, and $\Psi_2$ adds $\theta$
times the third coordinate to the second. This gives the same three coordinates.
\end{proof}

So $F$ is conjugate to $\Psi=\Psi_2\circ\Psi_1$, which acts on the coordinates $(X_1,W,S)$. The
coordinates $(W,S)$ are transformed by the plane automorphism
\[
\Psi'=(W+\theta(S+c_0g(W)^2),\;S+c_0g(W)^2),
\]
the composition of the de~Jonqui\`eres map $(W,S+c_0g(W)^2)$ of degree $2d$ with the
affine map $(W+\theta S,S)$, and $X_1$ is translated by $g(W)$. Thus $F$ is conjugate to
a skew product over a plane automorphism of H\'enon type.

\begin{proposition}\label{prop:heyan}
Let $F=X+H$ with $H$ as in \eqref{eq:heyan}.
\begin{enumerate}[label=\textup{(\roman*)},leftmargin=2em]
\item $F=E_2'\circ E_1\circ E_2\circ E_3$, where $E_3$ changes $X_3$ by $H_3$, $E_2$ changes $X_2$ by
$a^{-1}b(X_1)$, $E_1$ changes $X_1$ by $g(aX_2)$, and $E_2'$ changes $X_2$ by
$v_1(X_3-l_2X_1)-a^{-1}b(X_1)$. In particular $F\in\EE(K,3)$.
\item $\supp(F)=\{1,2,3\}$, $F$ has no left elementary factor, and
$\ell(F\times\Id_r)=\ell(F)=4$ for every $r\ge0$.
\item $\deg F=4d$, and $(2d)^m\le4\deg F^{(m)}$ and $\deg F^{(m)}\le4(2d)^m$ for every
$m\ge1$. Hence $\lambda(F)=2d$.
\item $F$ has infinite order, is not locally finite and not Pascal finite, and is not
conjugate in $\Aut_KK[X_1,X_2,X_3]^3$ to any automorphism whose iterates have bounded degree, such as
a triangular or a linear map. The matrix $J_H$ is nilpotent and not strongly nilpotent.
\end{enumerate}
\end{proposition}

\begin{proof}
(i) We apply $E_3$, $E_2$, $E_1$, $E_2'$ in this order. After $E_3$ and $E_2$ the coordinates are
$(X_1,\,X_2+a^{-1}b(X_1),\,X_3+H_3)$, and $a$ times the second one is $W$. So $E_1$
replaces $X_1$ by $X_1+g(W)=F_1$. Then $E_2'$ adds to the second coordinate
\[
v_1\bigl(X_3+H_3-l_2X_1-l_2g(W)\bigr)-a^{-1}b\bigl(X_1+g(W)\bigr).
\]
By the expansion of $b(X_1+g(W))$ in the proof of Lemma~\ref{lem:heyan-conj} and
$v_1H_3=v_1c_0g(W)^2+v_1l_2g(W)$, the second coordinate becomes
$X_2+v_1X_3-v_1l_2X_1-a^{-1}b'(X_1)g(W)=F_2$. The first and third coordinates are
$F_1$ and $F_3$. So $F$ is a product of four elementary maps.

(ii) $H_1=g(W)\ne0$, $H_3=g(W)(c_0g(W)+l_2)\ne0$, and $H_2$ contains $v_1X_3$, so
$\supp(F)=\{1,2,3\}$. By Lemma~\ref{lem:heyan-conj},
$F^{-1}=\varphi^{-1}\circ\Psi_1^{-1}\circ\Psi_2^{-1}\circ\varphi$. Put
$R=aX_2+b(X_1)-\theta(X_3-l_2X_1)$. Evaluating, one finds
\[
\begin{gathered}
(F^{-1})_1=X_1-g(R),\qquad (F^{-1})_3=X_3-c_0\,g(R)^2-l_2\,g(R),\\
(F^{-1})_2=X_2-v_1(X_3-l_2X_1)+a^{-1}b'(X_1)\,g(R)-v_1c_0\,g(R)^2 .
\end{gathered}
\]
Let $h_j=X_j-(F^{-1})_j$. Since $g'\ne0$,
\[
\begin{gathered}
\partial_{X_1}h_1=g'(R)\bigl(b'(X_1)+\theta l_2\bigr),\qquad
\partial_{X_3}h_3=-\theta\,g'(R)\bigl(2c_0g(R)+l_2\bigr),\\
\partial_{X_2}h_2=g'(R)\bigl(2\theta c_0\,g(R)-b'(X_1)\bigr)
\end{gathered}
\]
are nonzero. Indeed $b'(X_1)=2v_1c_0aX_1+l_1$ is not constant, and $g(R)$ involves
$X_2$ and $X_3$, since $R$ does and $g$ is not constant. By Lemma~\ref{lem:left-factor}, $F$ has no left
elementary factor. By Theorem~\ref{thm:min-decomp} it has no decomposition of length
$3=\#\supp(F)$, so $\ell(F)\ge4$, and $\ell(F)=4$ by (i). Corollary~\ref{cor:stab} gives
$\ell(F\times\Id_r)=4$.

(iii) The coordinate $F_3=X_3+c_0g(W)^2+l_2g(W)$ has degree $4d$, while $F_1$ and $F_2$
have degree at most $2d+1\le4d$. Let $\Psi=\Psi_2\circ\Psi_1$ and $\Psi^{(m)}=(x_m,w_m,s_m)$, so that
$(x_0,w_0,s_0)=(X_1,X_2,X_3)$ and
\[
s_{m+1}=s_m+c_0\,g(w_m)^2,\qquad w_{m+1}=w_m+\theta s_{m+1},\qquad
x_{m+1}=x_m+g(w_m).
\]
We show $\deg w_m=\deg s_m=(2d)^m$ for $m\ge1$. This holds for $m=1$. If it holds for
$m$, then $c_0g(w_m)^2$ has degree $2d(2d)^m$, which exceeds $\deg s_m$, so
$\deg s_{m+1}=(2d)^{m+1}$. Since $\theta\ne0$ and $\deg w_m<\deg s_{m+1}$, also
$\deg w_{m+1}=(2d)^{m+1}$. Moreover $\deg x_m\le d(2d)^{m-1}<(2d)^m$, so
$\deg \Psi^{(m)}=(2d)^m$. Finally $F^{(m)}=\varphi^{-1}\circ \Psi^{(m)}\circ\varphi$ and
$\Psi^{(m)}=\varphi\circ F^{(m)}\circ\varphi^{-1}$, and Lemma~\ref{lem:bounded}(ii) with
$\deg\varphi=\deg\varphi^{-1}=2$ gives both inequalities.

(iv) By (iii) the degrees of the iterates are unbounded. Since $H(0)=0$ is not assumed
in \eqref{eq:heyan}, we argue directly. A map of finite order, or a locally finite map,
has iterates of bounded degree, and a Pascal finite map is locally finite \cite{ABCH2}.
Lemma~\ref{lem:bounded} excludes conjugation to a map with iterates of bounded degree.
Finally $J_H$ is nilpotent by \cite[Theorem~2.3]{HY}. If $J_H$ were strongly
nilpotent, then $F$ would be linearly conjugate to a triangular map with all diagonal
coefficients $1$ \cite[Theorem~1.6]{EH}, see also \cite[Theorem~7.4.4]{vdE}, hence Pascal finite, which was just excluded.
\end{proof}

The degree of $F$ is $4d$, and the factors in (i) have degrees $4d$, at most $2$, $d$
and at most $2$. By the proof of (iii), $\deg \Psi^{(m)}=(2d)^m$ exactly. For $F=X+N_P$
the word in (i) is a rearrangement of the factorization \cite[(4.20)]{WY_nilp} at $s=1$,
which also has four elementary factors, and the inverse in the proof of (ii) is the
case $s=1$ of \cite[(4.23)]{WY_nilp}. The new information in
Proposition~\ref{prop:heyan} is the minimality in (ii) and the growth of the iterates in
(iii).

\begin{example}\label{ex:heyan}
Let $a=v_1=c_0=1$, $l_1=l_2=\tilde l_2=0$, so that $b=X_1^2$, $W=X_2+X_1^2$ and
$H=N_g=\bigl(g(W),\,X_3-2X_1g(W),\,g(W)^2\bigr)$. For $g(t)=t$ one obtains
$H=(W,\,X_3-2X_1W,\,W^2)$, with $\deg F=4$ and $\lambda(F)=2$. For $g(t)=t^2$ one
obtains $H=(W^2,\,X_3-2X_1W^2,\,W^4)$, with $\deg F=8$ and $\lambda(F)=4$. In both cases
$H(0)=0$, the components of $H$ are linearly independent, $\ell(F)=4$, and $F$ has no left
elementary factor. By contrast, for $H=(0,\,X_3,\,0)$, whose components are linearly
dependent, $F$ is elementary and Pascal finite, in accordance with
Corollary~\ref{cor:dim3-dichotomy} below.
\end{example}

The minimal decompositions of $F$ can be determined with
Proposition~\ref{prop:length-k1}. Put $W=aX_2+b(X_1)$ and $\phi=c_0g^2+l_2g$, so that
$H_3=\phi(W)$. We list a decomposition $Q_4\circ Q_3\circ Q_2\circ Q_1$ by its direction
sequence $(m_1,\dots,m_4)$, where $Q_i$ changes $X_{m_i}$, and by the tails of
$Q_1,\dots,Q_4$.

\begin{proposition}\label{prop:heyan-min}
Let $F=X+H$ with $H$ as in \eqref{eq:heyan}. The decompositions
$F=Q_4\circ Q_3\circ Q_2\circ Q_1$ of length $4$ are exactly those of
Table~\ref{tab:heyan-min}, where $Q_i$ changes $X_{m_i}$ by the tail $t_i$, with
$s\in K[X_3]$ in \textup{(A)}, $\sigma\in K$ in \textup{(B)} and $\tau\in K$ in \textup{(C)}.
\begin{table}[ht]
\centering
\small
\renewcommand{\arraystretch}{1.15}
\begin{tabular}{llll}
& $(m_1,m_2,m_3,m_4)$ & $i$ & $t_i$\\ \hline
(A) & $(3,2,1,2)$ & $1$ & $\phi(W)$\\
& & $2$ & $v_1X_3-s+a^{-1}b(X_1)$\\
& & $3$ & $g(aX_2+as-\theta X_3)$\\
& & $4$ & $s-v_1l_2X_1-a^{-1}b(X_1)$\\ \hline
(B) & $(2,3,1,2)$ & $1$ & $a^{-1}b(X_1)-\sigma$\\
& & $2$ & $\phi(aX_2+a\sigma)$\\
& & $3$ & $g(aX_2+a\sigma)$\\
& & $4$ & $v_1X_3+\sigma-v_1l_2X_1-a^{-1}b(X_1)$\\ \hline
(B) & $(2,1,3,2)$ & $1$ & $a^{-1}b(X_1)-\sigma$\\
& & $2$ & $g(aX_2+a\sigma)$\\
& & $3$ & $\phi(aX_2+a\sigma)$\\
& & $4$ & $v_1X_3+\sigma-v_1l_2X_1-a^{-1}b(X_1)$\\ \hline
(C) & $(3,2,1,3)$ & $1$ & $-\theta^{-1}\bigl(b(X_1)+b'(X_1)g(W)\bigr)-l_2X_1-\tau$\\
& & $2$ & $v_1X_3+v_1\tau+a^{-1}b(X_1)$\\
& & $3$ & $g(aX_2-\theta X_3-\theta\tau)$\\
& & $4$ & $\theta^{-1}b(X_1)+l_2X_1+\tau$
\end{tabular}
\caption{The decompositions of length $4$ of the maps \eqref{eq:heyan}. The factor
$Q_1$ is applied first.}\label{tab:heyan-min}
\end{table}

These decompositions form exactly two move classes, \textup{(A)}$\,\cup\,$\textup{(B)}
and \textup{(C)}. The decomposition of Proposition~\ref{prop:heyan}\textup{(i)} is
\textup{(A)} with $s=v_1X_3$. Moreover $\omega_0(F)=2d$, and this value is attained
exactly by the decompositions \textup{(B)}.
\end{proposition}

The two move classes are not two decompositions. Each of them contains, over an infinite
field, infinitely many decompositions, parametrized by $s$, $\sigma$ or $\tau$, and the
class (A)$\,\cup\,$(B) contains three direction sequences.

\begin{proof}
By Proposition~\ref{prop:heyan}(ii) and Proposition~\ref{prop:length-k1}, every
decomposition is of type (b). So $Q_4$ changes some $X_c$ by $t\in K[X_i:i\ne c]$,
$G=Q_4^{-1}\circ F$ has the inverse $G^{-1}=F^{-1}\circ Q_4$, and $Q_3,Q_2,Q_1$ are the
left elementary factors of $G$ extracted along an ordering of $\{1,2,3\}$ which does not
begin with $c$. We use the inverse $F^{-1}$ from the proof of
Proposition~\ref{prop:heyan}(ii), the identity
\[
b(X_1+y)=b(X_1)+b'(X_1)\,y+v_1c_0a\,y^2 ,\tag{$\ast$}
\]
and the fact that, for $\operatorname{char}K=0$ and a non-constant polynomial $\psi$,
$\psi(u)$ does not involve $X_i$ if and only if $u$ does not, since
$\partial_{X_i}\psi(u)=\psi'(u)\,\partial_{X_i}u$ and $\psi'(u)\ne0$ for non-constant
$u$. All polynomials $u$ below involve $X_2$ or $X_3$, so they are non-constant.

\emph{Case $c=2$.} Here $t\in K[X_1,X_3]$, and in $G^{-1}=F^{-1}\circ Q_4$ the argument $R$
of $g$ becomes $R+at$. If $Q_3$ changes $X_1$, the test requires that $R+at$ does
not involve $X_1$, that is, $t=s(X_3)-v_1l_2X_1-a^{-1}b(X_1)$ with $s\in K[X_3]$. Then
$R+at=u=aX_2+as-\theta X_3$, $Q_3$ changes $X_1$ by $g(u)$, and by $(\ast)$
\[
G^{-1}\circ Q_3=\bigl(X_1,\;X_2+s-v_1X_3-a^{-1}b(X_1),\;X_3-\phi(u)\bigr).
\]
Extracting $X_2$ next gives $Q_2$ with tail $v_1X_3-s+a^{-1}b(X_1)$ for every $s$, and
the remaining map changes $X_3$ by $\phi(W)$. This is (A). Extracting $X_3$ next
requires that $u$ does not involve $X_3$, that is, $s=v_1X_3+\sigma$, and gives the
sequence $(2,3,1,2)$ of (B). If $Q_3$ changes $X_3$, the test requires that $R+at$ does
not involve $X_3$, that is, $t=v_1X_3+\mu(X_1)$. Then
$R+at=u=aX_2+b(X_1)+\theta l_2X_1+a\mu(X_1)$, and $G^{-1}\circ Q_3$ has the second
coordinate $X_2+\mu+v_1l_2X_1+a^{-1}b'(X_1)g(u)-v_1c_0g(u)^2$, of degree $2d$ in $X_2$.
So $Q_2$ changes $X_1$, the test requires that $u$ does not involve $X_1$, that is,
$\mu=\sigma-v_1l_2X_1-a^{-1}b(X_1)$, and this gives the sequence $(2,1,3,2)$ of (B).

\emph{Case $c=3$.} Here $t\in K[X_1,X_2]$, and the argument $R$ of $g$ becomes $R-\theta t$. If
$Q_3$ changes $X_1$, then $t=\theta^{-1}b(X_1)+l_2X_1+\tau(X_2)$, $R-\theta t=u=
aX_2-\theta X_3-\theta\tau(X_2)$, and by $(\ast)$
\[
G^{-1}\circ Q_3=\bigl(X_1,\;X_2-v_1X_3-v_1\tau-a^{-1}b(X_1),\;
X_3+\theta^{-1}b(X_1)+l_2X_1+\tau+\theta^{-1}b'(X_1)g(u)\bigr).
\]
The third coordinate minus $X_3$ involves $X_3$ through $g(u)$, so $Q_2$ changes $X_2$,
which requires $\tau\in K$ and gives (C). If $Q_3$ changes $X_2$, put
$\rho=R-\theta t$. The tail $X_2-(G^{-1})_2$ has the derivative
\[
v_1t_{X_2}+(1-v_1t_{X_2})\,g'(\rho)\bigl(2\theta c_0\,g(\rho)-b'(X_1)\bigr)
\]
with respect to $X_2$. Since $\rho$ has degree $1$ in $X_3$, the second term has degree
$2d-1\ge1$ in $X_3$ unless $v_1t_{X_2}=1$, and then the derivative is $1$. So this case
does not occur.

\emph{Case $c=1$.} Here $t\in K[X_2,X_3]$. Put $w=X_1+t$ and
$r=aX_2+b(w)-\theta X_3+\theta l_2w$, so that $G^{-1}$ is obtained from $F^{-1}$ by
replacing $X_1$ with $w$. Since $(X_1,X_2,X_3)\mapsto(w,X_2,X_3)$ is an automorphism,
polynomial identities may be read in the variables $w,X_2,X_3$. If $Q_3$ changes $X_3$,
the test requires $\partial_{X_3}r=r_wt_{X_3}-\theta=0$ with $r_w=b'(w)+\theta l_2$ of
degree $1$ in $w$, which is impossible. If $Q_3$ changes $X_2$, the tail is
$h=\mathcal H(w,X_2,X_3)$ with
$\mathcal H=v_1(X_3-l_2w)-a^{-1}b'(w)g(r)+v_1c_0g(r)^2$, and $\partial_{X_2}h=0$ reads
\[
\mathcal H_{X_2}\bigl(1+a^{-1}t_{X_2}r_w\bigr)=t_{X_2}\bigl(v_1l_2+2v_1c_0\,g(r)\bigr),
\qquad \mathcal H_{X_2}=g'(r)\bigl(2\theta c_0\,g(r)-b'(w)\bigr).
\]
In $w$, $\mathcal H_{X_2}$ has degree $4d-2$ and a constant leading coefficient, and
$g(r)$ has degree $2d$. If $t_{X_2}\ne0$, the left side has degree $4d-1>2d$, and if
$t_{X_2}=0$, the equation reads $\mathcal H_{X_2}=0$. So this case does not occur.

Conversely, each list is a decomposition of $F$, as the extraction above shows. The
factors of (A) include $\phi(W)$, of degree $4d$, those of (C) include a factor of degree
$2d+1$, and those of (B) have degrees $2$, $2d$, $d$ and $2$. Hence $\omega_0(F)=2d$.

It remains to determine the move classes. In (A), the move (M2) with $i=2$, $j=4$ and
$P\colon X_2\mapsto X_2+s'-s$ replaces the parameter $s$ by $s'$, since $s'-s\in
K[X_3]$ does not involve $X_1$. For $s=v_1X_3+\sigma$, the first two factors of (A)
apply $X_3\mapsto X_3+\phi(W)$ and then $X_2\mapsto X_2+a^{-1}b(X_1)-\sigma$, and the
first two factors of the sequence $(2,3,1,2)$ of (B) apply
$X_2\mapsto X_2+a^{-1}b(X_1)-\sigma$ and then $X_3\mapsto X_3+\phi(aX_2+a\sigma)$. Both
products equal $(X_1,\,X_2+a^{-1}b(X_1)-\sigma,\,X_3+\phi(W))$, and the last two factors
coincide, so (M1) relates the two decompositions.
The two sequences of (B) differ by the exchange of the commuting maps
$X_3\mapsto X_3+\phi(aX_2+a\sigma)$ and $X_1\mapsto X_1+g(aX_2+a\sigma)$. In (C), the
move (M2) with $i=1$, $j=4$ and $P\colon X_3\mapsto X_3+\varepsilon$ replaces $\tau$ by
$\tau+\varepsilon$. Thus (A) and (B) form one class, and (C) is contained in one class.
It remains to show that no move leads from (C) to (A) or (B). Every decomposition of
length $4$ occurs in Table~\ref{tab:heyan-min}, and a decomposition belongs to (C) if and
only if its direction sequence is $(3,2,1,3)$. A move (M2) preserves the direction
sequence. A move (M1) replaces two adjacent factors and leaves the other two unchanged,
so it changes at most two adjacent entries of the direction sequence. Starting from
$(3,2,1,3)$, a change of the first two entries gives a sequence ending with $3$, and a
change of the middle two entries gives a sequence beginning and ending with $3$. No
sequence of (A) or (B) has either property. So a move (M1) from (C) to (A) or (B) changes
the last two entries, keeps the first two factors, and leads to the sequence $(3,2,1,2)$
of (A). This is
impossible, since in (C) the product of the last two factors changes $X_3$, by
$\theta^{-1}b(X_1)+l_2X_1+\tau$ composed with an automorphism, while in (A) it fixes $X_3$.
\end{proof}

\begin{corollary}\label{cor:dim3-dichotomy}
Let $\operatorname{char}K=0$, let $H\in K[X_1,X_2,X_3]^3$ satisfy $H(0)=0$ with $J_H$
nilpotent, and let $F=X+H$.
\begin{enumerate}[label=\textup{(\alph*)},leftmargin=2em]
\item If the components of $H$ are linearly dependent over $K$, then $F$ is Pascal
finite.
\item If they are linearly independent, then $\deg H=4d$ for some $d\ge1$,
$\lambda(F)=2d=\tfrac12\deg H$, $F$ is not Pascal finite and has infinite order, and
$J_H$ is not strongly nilpotent.
\end{enumerate}
Hence the conditions of Lemma~\ref{lem:pf-lf} are also equivalent to the linear
dependence of the components of $H$.
\end{corollary}

\begin{proof}
(b) By \cite[Theorem~4.1]{WY_nilp} there are $T\in\mathrm{GL}_3(K)$ and $P$ with
$T\circ H\circ T^{-1}=N_P$. Since $T$ is linear, $T\circ F\circ T^{-1}=X+N_P$, and
$\deg H=\deg N_P=4d$ with $d=\deg P$. Proposition~\ref{prop:heyan}(iii) and (iv)
apply to $X+N_P$. Pascal finiteness is invariant under conjugation
\cite[Theorem~2.2]{ABCH2}, and $\lambda(F)=\lambda(X+N_P)$ by
Lemma~\ref{lem:bounded}(ii), since $\deg T=\deg T^{-1}=1$. The matrix $J_H$ is not
strongly nilpotent by Proposition~\ref{prop:heyan}(iv) applied to $X+N_P$, since strong
nilpotency is preserved by linear conjugation.

(a) By the proof of \cite[Theorem~4.6]{WY_nilp}, after a linear conjugation
$H=(h_1,h_2,0)$, and there is a tame automorphism $S$ such that $S\circ F\circ S^{-1}$ is
of the form $(X_1+u(X_3),\,X_2+w(X_1,X_3),\,X_3)$ \cite[(4.21)]{WY_nilp} at $s=1$.
This is a triangular map with all diagonal coefficients $1$, hence Pascal finite, and so
is $F$ \cite[Theorems~2.1, 2.2]{ABCH2}.

The last statement follows from (a) and (b).
\end{proof}

By Lemma~\ref{lem:pf-lf}, the equivalence of Pascal finiteness, local finiteness and
boundedness of the degrees of the iterates holds in every dimension. The content of
Corollary~\ref{cor:dim3-dichotomy} is that in dimension three these properties are
decided by the linear dependence of the components of $H$, and that otherwise the degrees
grow exactly like $(\tfrac12\deg H)^m$ up to a factor $4$. In characteristic zero
strong nilpotency of $J_H$ implies Pascal finiteness in every dimension, by
\cite[Theorem~1.6]{EH}, see also \cite[Theorem~7.4.4]{vdE}, and since triangular maps with all diagonal coefficients $1$ are
Pascal finite, see also \cite[Theorem~5.1]{AH_sn}. In dimension three, for
$H(0)=0$, nilpotency of $J_H$ together with linear independence of the components
excludes it.

\begin{remark}\label{rem:heyan-compare}
Vasyunin's map and the maps of Proposition~\ref{prop:heyan} both have nilpotent and not
strongly nilpotent $J_H$, and neither is Pascal finite. They differ in three respects.
\begin{enumerate}[label=(\arabic*),leftmargin=2em]
\item Vasyunin's map satisfies $\ell(F)=\#\supp(F)$, and Theorem~\ref{thm:min-decomp}
lists all its minimal decompositions. They are six, and they form one move class,
since $P_2$ and $P_5$ commute and any two orders of $P_2,P_5,P_4^{(c)}$ which differ by an
adjacent transposition give decompositions related by (M1). The maps of He and Yan satisfy
$\ell(F)=\#\supp(F)+1$, so they have infinitely many minimal decompositions, and these
form two move classes (Proposition~\ref{prop:heyan-min}). By
Proposition~\ref{prop:stab-bound} the addition of variables does not shorten them.
\item The dynamical degree of Vasyunin's map is $\sqrt2$, while for the maps of He and
Yan it is the even integer $2d$. The reason is Lemma~\ref{lem:heyan-conj}. Up to
conjugation, the dynamics of $F$ is that of the plane automorphism $\Psi'$, whose reduced
word consists of one affine and one de~Jonqui\`eres letter, so that
$\lambda(\Psi')=2d$ as in Remark~\ref{rem:dyn-length}.
\item For $n\le4$ every quadratic homogeneous automorphism is linearly triangularizable
(Section~\ref{ssec:vasyunin-def}), so $n=5$ is the least dimension of a quadratic
homogeneous map $X+H$ with nilpotent and not strongly nilpotent $J_H$, and Vasyunin's map
is such a map. Without homogeneity the maps of He and Yan give such maps for $n=3$.
\end{enumerate}
Since $F^{(m)}=(E_2'\circ E_1\circ E_2\circ E_3)^{(m)}$, one has $\ell(F^{(m)})\le4m$. By
Remark~\ref{rem:dyn-length} the plane automorphism $\Psi'$ satisfies $d\ell(\Psi')>0$, so the
maps of He and Yan are test cases for Question~\ref{q:dyn-length} in dimension three in which
the corresponding plane question is answered.
\end{remark}

\section{The Anick automorphism}\label{sec:anick}

Let $K$ be a field, $\delta=X_1X_3-X_2X_4$ and
\[
A=\bigl(X_1+\delta X_4,\;X_2+\delta X_3,\;X_3,\;X_4\bigr).
\]
This automorphism was proposed by Anick as a candidate for a wild automorphism,
see \cite[p.~343]{Coh}. Miasnikov and Shestakov \cite{MS} show that it is tame if
$\operatorname{char}K\ne2$. They write it as a product of sixteen elementary maps
fixing $X_3$, after the linear change of coordinates
$V=\tfrac12(X_2+X_4)$, $U=\tfrac12(X_2-X_4)$, and their identities have coefficients in
$\mathbb{Z}[\tfrac12]$. Since the change of coordinates also fixes $X_3$, $A$ is tame as an
automorphism of $K[X_3][X_1,X_2,X_4]$ over $K[X_3]$ \cite[Section~5]{MS}. In
characteristic zero $A=\exp(\delta D_0)$ with
$D_0=X_4\partial_{X_1}+X_3\partial_{X_2}$, and $D_0(\delta)=0$.

\begin{proposition}\label{prop:anick}
\begin{enumerate}[label=\textup{(\roman*)},leftmargin=2em]
\item $A=X+H$ with $H=(\delta X_4,\delta X_3,0,0)$ is a quasi-translation, and
$A^{(m)}=X+mH$ for $m\ge1$. Hence $\tau_K(A)=2$, $A$ has order $p$ if
$\operatorname{char}K=p>0$, and infinite order if $\operatorname{char}K=0$.
\item $\supp(A)=\{1,2\}$, $A$ has no left elementary factor, and no word of length at
most $3$ in $\cE(K,4)$ represents $A$.
\item For every field $K$, $A\times\Id_1\in\EE(K,5)$ and $\ell(A\times\Id_1)\le6$,
with equality by Theorem~\ref{thm:anick-seven}.
More precisely, $A\times\Id_1=T^{-1}\circ E^{-1}\circ T\circ E$, where $E$ changes $X_5$
by $\delta$ and $T=(X_1+X_5X_4,\,X_2+X_5X_3,\,X_3,\,X_4,\,X_5)$ is the product of the
elementary maps changing $X_1$ by $X_5X_4$ and $X_2$ by $X_5X_3$.
\item If $\operatorname{char}K\ne2$, then $A\in\EE(K,4)$ and $\ell(A)\le13$.
\end{enumerate}
\end{proposition}

\begin{proof}
(i) One has $\delta\circ A=(X_1+\delta X_4)X_3-(X_2+\delta X_3)X_4=\delta$, and $A$ fixes
$X_3$ and $X_4$. Hence $H\circ A=H$, so $A$ is a quasi-translation, and
$A^{(m+1)}=A^{(m)}+H\circ A^{(m)}=X+(m+1)H$ by induction. The equality $\tau_K(A)=2$
holds for every quasi-translation \cite{ABCH2}. Since $H\ne0$ has coefficients in
$\mathbb{Z}$ and content $1$, $mH=0$ in $K[X]^4$ if and only if
$\operatorname{char}K$ divides $m$.

(ii) Clearly $\supp(A)=\{1,2\}$. By (i), $A^{-1}=X-H$, so
$X_1-(A^{-1})_1=X_1X_3X_4-X_2X_4^2$ involves $X_1$ and
$X_2-(A^{-1})_2=X_1X_3^2-X_2X_3X_4$ involves $X_2$. By Lemma~\ref{lem:left-factor}, $A$ has
no left elementary factor. A word of length at most $2$ for $A$ would have length
$\#\supp(A)$ by Lemma~\ref{lem:support-bound}, which is excluded by
Theorem~\ref{thm:min-decomp}. Suppose that $A=Q_3\circ Q_2\circ Q_1$, with $Q_i$ changing
$X_{m_i}$. Since words of length at most $2$ are excluded, $\ell(A)=3$, and since $A$ has no left elementary factor,
Proposition~\ref{prop:length-k1}(b) gives $m_1=m_3=c$ and $m_2=j$ with
$\{c,j\}=\{1,2\}$.

Let $c=1$, so that $Q_1$ changes $X_1$ by $a\in K[X_2,X_3,X_4]$ and $Q_2$ changes $X_2$
by $b\in K[X_1,X_3,X_4]$. Since $Q_3$ fixes $X_2$, comparison of the second coordinates
gives $b(X_1+a,X_3,X_4)=X_1X_3^2-X_2X_3X_4$. Substituting $X_1-a$ for $X_1$ gives
$b(X_1,X_3,X_4)=X_1X_3^2-aX_3^2-X_2X_3X_4$. The left-hand side does not involve $X_2$,
so $aX_3^2+X_2X_3X_4=\beta\in K[X_1,X_3,X_4]$. Comparing the coefficients of $X_2$,
as polynomials in $X_2$ over $K[X_1,X_3,X_4]$, gives $X_3^2\mid X_3X_4$, a contradiction.
For $c=2$ the same argument, with the first coordinate
$b(X_2+a,X_3,X_4)=X_1X_3X_4-X_2X_4^2$, gives $X_4^2\mid X_3X_4$, a contradiction.

(iii) Apply $E$, $T$, $E^{-1}$ and $T^{-1}$ in this order. After $E$ and $T$ the
coordinates are $(X_1+(X_5+\delta)X_4,\,X_2+(X_5+\delta)X_3,\,X_3,\,X_4,\,X_5+\delta)$, and $\delta$
evaluated at them is $\delta$. So $E^{-1}$ restores $X_5$, and $T^{-1}$ subtracts
$X_5X_4$ and $X_5X_3$, which gives $A\times\Id_1$. The coefficients lie in $\mathbb{Z}$.
This is the construction of Smith \cite{Sm}. For the Nagata automorphism an explicit
decomposition of this kind into six elementary maps, over every commutative ring, is given
in \cite{Sp}.

(iv) Let $\Lambda=(X_1,\,\tfrac12(X_2+X_4),\,X_3,\,X_4-X_2)$ and
$D=\mathrm{diag}(1,1,1,-\tfrac12)$. Then $D\circ\Lambda$ is the change of coordinates of
\cite{MS}, so $A=\Lambda^{-1}\circ D^{-1}\circ\tilde A\circ D\circ\Lambda$, where
$\tilde A$ is the product of the sixteen elementary maps of \cite[Section~5]{MS}.
Conjugation by $D$ preserves $\cE(K,4)$. The map $\Lambda$ is the product of two
elementary maps. One first changes $X_4$ by $-X_2$ and then $X_2$ by $\tfrac12X_4$, in
accordance with the identity
$\left(\begin{smallmatrix}1/2&1/2\\-1&1\end{smallmatrix}\right)
=\left(\begin{smallmatrix}1&1/2\\0&1\end{smallmatrix}\right)
\left(\begin{smallmatrix}1&0\\-1&1\end{smallmatrix}\right)$.

The word of \cite{MS} can be shortened. In the coordinates $(x,v,z,u)$ of \cite{MS}, for
$h\in K[u,z]$ let $T_h$ change $v$ by $hz$ and then $x$ by $2hv-h^2z$, and for
$k\in K[v,z]$ let $R_k$ change $u$ by $kz$ and then $x$ by $-2ku+k^2z$. A direct
computation gives
\[
T_h=\bigl(x+2hv+h^2z,\;v+hz,\;z,\;u\bigr),\qquad R_k=\bigl(x-2ku-k^2z,\;v,\;z,\;u+kz\bigr).
\]
Since $T_{h'}$ does not change $h$ and $R_{k'}$ does not change $k$, it follows that
$T_h\circ T_{h'}=T_{h+h'}$ and $R_k\circ R_{k'}=R_{k+k'}$. The factors $5$ to $8$ of
\cite[Section~5]{MS} are $T_{-u}$ followed by $T_{u^2/2}$, and the factors $13$ to $16$ are
$R_{-v}$ followed by $R_{-v^2/2}$. Replacing them by $T_{u^2/2-u}$ and $R_{-v-v^2/2}$ gives
a word of twelve elementary maps for $\tilde A$. In the notation of formula (1) of
\cite[Section~3]{MS} this merges $T_u^{-1}$ with $T_{u^2/2}$ and $R_v^{-1}$ with
$R_{-v^2/2}$. Hence $\ell(A)\le2+12+2=16$. Two further mergings, which use the
commutation of the maps $T_f$ and a different factorization of the change of coordinates,
give a word of $13$ elementary maps (Proposition~\ref{prop:app-13} in
Appendix~\ref{app:anick}). So $\ell(A)\le13$.
\end{proof}

The degree of $A$ is $3$, and the elementary maps of the shortened word have degrees up
to $5$. No two adjacent factors of the shortened word change the same coordinate. The
proof of the upper bound in (iv) uses the coefficient $\tfrac12$. In characteristic
$2$ we do not know whether $A$ is tame, while by (iii) it is stably tame over every
field.

\begin{remark}\label{rem:anick-iterates}
Let $m\in\mathbb{Z}$ with $m\ne0$ in $K$, and let $D=\mathrm{diag}(1,m,m,1)$. Since
$\delta\circ D=m\delta$, one has $D^{-1}\circ A\circ D=X+mH=A^{(m)}$. Conjugation by a
diagonal map preserves $\cE(K,4)$, so $\ell(A^{(m)})=\ell(A)$. Hence $d\ell(A)=0$ if
$\operatorname{char}K=0$. Moreover $\deg A^{(m)}=3$, so $\lambda(A)=1$. By
Theorem~\ref{thm:anick-seven} below the iterates of $A$ have constant length at least
$7$ and bounded degree, while Vasyunin's map has $\ell(F)=4$ and $\lambda(F)=\sqrt2$. So
the length of an automorphism does not bound the growth of its iterates, nor conversely.
\end{remark}

\begin{remark}\label{rem:anick-AS}
Let $L=K(X_1,X_2,X_3,X_4)$, $B=K(X_3,X_4,\delta)$ and $\alpha=X_1/(\delta X_4)\in L$. Since
$\delta\circ A=\delta$, one has $\alpha\circ A=\alpha+1$, and $A$ acts on $L$ over $B$.
Moreover $X_1=\alpha\delta X_4$ and $X_2=\alpha\delta X_3-\delta/X_4$, so $L=B(\alpha)$, and
$\alpha$ is transcendental over $B$, since $X_3,X_4,\delta$ are algebraically independent.
Thus $A$ acts on the rational function field $B(\alpha)$ by $\alpha\mapsto\alpha+1$. If
$\operatorname{char}K=0$, this automorphism has infinite order and its fixed field is $B$.
If $\operatorname{char}K=p>0$, it has order $p$, and by the theorem of Artin its fixed
field is $B(\alpha^p-\alpha)$, since $\alpha^p-\alpha$ is fixed and
$[B(\alpha):B(\alpha^p-\alpha)]=p$. So $L$ is the Artin--Schreier extension of the fixed
field defined by $T^p-T=\alpha^p-\alpha$. If $\operatorname{char}K=0$, then
$A=\exp(\delta D_0)$, and $\delta D_0(\alpha)=1$.
\end{remark}

The computation in (iii) is an instance of the following identity. Let
$V\in K[X]^n$ and $\delta\in K[X]$ satisfy $V(X+cV(X))=V(X)$ and $\delta(X+cV(X))=\delta(X)$
for a new variable $c$, and in $K[X,u]^{n+1}$ let
$T(X,u)=(X+uV(X),u)$ and $E(X,u)=(X,u+\delta(X))$. Then $T^{-1}(X,u)=(X-uV(X),u)$, and
\[
T^{-1}\circ E^{-1}\circ T\circ E=\bigl(X+\delta(X)V(X),\,u\bigr),
\]
since $E^{-1}$ subtracts $\delta(X+(u+\delta)V)=\delta$ from $u+\delta$, and $T^{-1}$ then
subtracts $uV$. For the Anick automorphism $V=(X_4,X_3,0,0)$, and this is the
construction of Smith \cite{Sm}.

The lower bounds below rest on the theorem of Cohn. Relative to $K[X_3,X_4]$ the map $A$
is wild, so every decomposition of $A$ changes $X_3$ or $X_4$
(Proposition~\ref{prop:app-base}). The remaining configurations of length at most $6$ are
excluded by comparisons of coefficients in $X_1,X_2$ over $K[X_3,X_4]$.

\begin{lemma}\label{lem:anick-sequences}
Let $A=Q_4\circ Q_3\circ Q_2\circ Q_1$ with $Q_i\in\cE(K,4)$ changing $X_{m_i}$. Then
$(m_1,m_2,m_3,m_4)$ is one of
\[
(1,2,1,2),\quad(2,1,2,1),\quad(c,1,2,c),\quad(c,2,1,c),\qquad c\in\{3,4\}.
\]
\end{lemma}

\begin{proof}
Words of length at most $3$ are excluded by Proposition~\ref{prop:anick}(ii), so adjacent
directions differ by Lemma~\ref{lem:word}(ii). By Lemma~\ref{lem:support-bound}, $X_1$ and
$X_2$ are each changed at least once, and by Lemma~\ref{lem:word}(i) each of $X_3$, $X_4$
is changed either never or at least twice, so at most one of them occurs. If $m_4$ occurred
only once, $Q_4$ would be a left elementary factor of $A$. If $m_1$ occurred only once,
$Q_1^{-1}$ would be a left elementary factor of $A^{-1}=Q_1^{-1}\circ\dots\circ Q_4^{-1}$,
which is excluded by Lemma~\ref{lem:left-factor}, since $X_1-A_1=-\delta X_4$ involves
$X_1$ and $X_2-A_2=-\delta X_3$ involves $X_2$. So the first and the last direction occur
at least twice. If neither $3$ nor $4$ occurs, the multiset of directions is
$\{1,1,2,2\}$, since three equal directions cannot be pairwise non-adjacent in a
sequence of length $4$, and the sequence alternates. Otherwise some $c\in\{3,4\}$ occurs
exactly twice, it is the only repeated direction, and $m_1=m_4=c$.
\end{proof}

\begin{lemma}\label{lem:anick-tools}
Let $\rho=(X_1,-X_2,-X_3,X_4)$, $\pi=(X_2,X_1,X_4,X_3)$, $A^{-1}=(P,Q,X_3,X_4)$ and
$\{c,c'\}=\{3,4\}$.
\begin{enumerate}[label=\textup{(\roman*)},leftmargin=2em]
\item $\rho\circ A\circ\rho=\pi\circ A\circ\pi=A^{-1}$. If $A$ has a decomposition of
length $l$ with direction sequence $(m_1,\dots,m_l)$, then it has decompositions of length
$l$ with the direction sequences $(m_l,\dots,m_1)$ and
$(\varsigma(m_1),\dots,\varsigma(m_l))$, where $\varsigma$ exchanges $1$ with $2$ and
$3$ with $4$.
\item If $h\in K[X_1,X_2,X_{c'}]$ and $h(P,Q,X_{c'})$ does not involve $X_c$, then
$h\in K[X_{c'}]$.
\item Let $A=Q_l\circ\dots\circ Q_1$ with $Q_i\in\cE(K,4)$ changing $X_{m_i}$. If
$j\in\{1,2\}$ occurs exactly once in $(m_1,\dots,m_l)$, at the position $i$, then
$i\ge3$.
\end{enumerate}
\end{lemma}

\begin{proof}
(i) One has $\delta\circ\rho=\delta\circ\pi=-\delta$, which gives the first statement by
direct substitution. The maps $\rho$ and $\pi$ are linear with monomial matrices, so
conjugation by $\rho$ maps an elementary map changing $X_m$ to an elementary map changing
$X_m$, and conjugation by $\rho\circ\pi$ maps it to one changing $X_{\varsigma(m)}$. If
$A=Q_l\circ\dots\circ Q_1$, then
$A=\rho\circ A^{-1}\circ\rho=(\rho\circ Q_1^{-1}\circ\rho)\circ\dots\circ(\rho\circ
Q_l^{-1}\circ\rho)$ and $A=(\rho\circ\pi)\circ A\circ(\rho\circ\pi)^{-1}$.

(ii) Write $h=\sum_{i,j}\gamma_{ij}X_1^iX_2^j$ with
$\gamma_{ij}\in K[X_{c'}]$.

Let $c=3$. As polynomials in $X_3$, $P$ has degree $1$ with leading coefficient $-X_1X_4$,
and $Q$ has degree $2$ with leading coefficient $-X_1$. Let $D$ be the largest value of
$i+2j$ with $\gamma_{ij}\ne0$. The coefficient of $X_3^D$ in $h(P,Q,X_4)$ is
$\sum_{i+2j=D}(-1)^{i+j}\gamma_{ij}X_4^iX_1^{i+j}$. Distinct pairs with $i+2j=D$ have
distinct values of $i+j$, so these terms have distinct degrees in $X_1$ and do not
cancel. Hence $D=0$. Let $c=4$. As polynomials in $X_4$, $P$ has degree $2$ with leading
coefficient $X_2$, and $Q$ has degree $1$ with leading coefficient $X_2X_3$. With $D$ the
largest value of $2i+j$, the coefficient of $X_4^D$ is
$\sum_{2i+j=D}\gamma_{ij}X_3^jX_2^{i+j}$, and distinct pairs give distinct degrees in
$X_2$. Hence again $D=0$, that is, $h\in K[X_{c'}]$.

(iii) The factors after $Q_i$ fix $X_j$, so $A_j-X_j=h_i\circ Q_{i-1}\circ\dots\circ
Q_1$, where $h_i$ is the tail of $Q_i$ and does not involve $X_j$. If $i=1$, this
contradicts the fact that $A_1-X_1=\delta X_4$ involves $X_1$ and $A_2-X_2=\delta X_3$
involves $X_2$. Let $i=2$, and let $Q_1$ change $X_m$, $m\ne j$, by $s$. Then
$h_2=(A_j-X_j)\circ Q_1^{-1}$ is obtained from $A_j-X_j$ by replacing $X_m$ with $X_m-s$,
and we show that it involves $X_j$. Let $s_1$ denote the coefficient of $X_j$ in $s$,
regarded as a polynomial in $X_j$. For $j=1$ and $m=2,3,4$ the polynomial $h_2$ is
\[
X_1X_3X_4-(X_2-s)X_4^2,\qquad X_1(X_3-s)X_4-X_2X_4^2,\qquad
X_1X_3(X_4-s)-X_2(X_4-s)^2 .
\]
In the first, the coefficient of $X_1$ is $X_4(X_3+s_1X_4)\ne0$. In the second, the
coefficient of $X_3$ is $X_1X_4$. In the third, the coefficient of $X_4$ is
$X_1X_3+2X_2s$, whose coefficient of $X_1$ is $X_3+2X_2s_1\ne0$. For $j=2$ and $m=1,3,4$
the polynomial $h_2$ is
\[
(X_1-s)X_3^2-X_2X_3X_4,\qquad X_1(X_3-s)^2-X_2(X_3-s)X_4,\qquad X_1X_3^2-X_2X_3(X_4-s).
\]
In the first, the coefficient of $X_2$ is $-X_3(s_1X_3+X_4)\ne0$. In the second, the
coefficient of $X_3$ is $-2X_1s-X_2X_4$, whose coefficient of $X_2$ is
$-2X_1s_1-X_4\ne0$. In the third, the coefficient of $X_4$ is $-X_2X_3$. Each of the stated
coefficients contains one of the monomials $X_3$, $X_4$, $X_3X_4$, $X_1X_4$, $X_2X_3$ with
coefficient $\pm1$, while its other term is a multiple of $s_1X_4^2$, $X_2s_1$,
$s_1X_3^2$ or $X_1s_1$ and does not contain this monomial. So it is nonzero over every
field.
\end{proof}

\begin{proposition}\label{prop:anick-five}
No word of length $4$ in $\cE(K,4)$ represents $A$.
\end{proposition}

\begin{proof}
Let $A=Q_4\circ Q_3\circ Q_2\circ Q_1$, and let $h_i$ be the tail of $Q_i$. By
Lemma~\ref{lem:anick-sequences} there are two cases.

\emph{The sequences $(1,2,1,2)$ and $(2,1,2,1)$.} Let $R=K[X_3,X_4]$. All factors fix
$X_3$ and $X_4$, so we regard them as automorphisms of $R[X_1,X_2]$, and $A$ is the
linear map with matrix
\[
M=\begin{pmatrix}1+X_3X_4&-X_4^2\\X_3^2&1-X_3X_4\end{pmatrix}.
\]
Consider $(1,2,1,2)$ and write $h_1=p_1(X_2)$, $h_2=p_2(X_1)$, $h_3=p_3(X_2)$, $h_4=p_4(X_1)$
with $p_1,p_2,p_3,p_4\in R[t]$. Since $Q_4$ fixes $X_1$,
\[
A_1=X_1+p_1(X_2)+p_3\bigl(X_2+p_2(X_1+p_1(X_2))\bigr)=\lambda X_1-X_4^2X_2,\qquad
\lambda=1+X_3X_4 .
\]
Substituting $X_1-p_1(X_2)$ for $X_1$ gives
$p_3(X_2+p_2(X_1))=(\lambda-1)X_1-\lambda p_1(X_2)-X_4^2X_2$, and $X_2=0$ gives
$p_3(p_2(X_1))=X_3X_4X_1-\lambda p_1(0)$. Since $R$ is a domain, $\deg p_3\cdot\deg p_2=1$, so $p_2$
and $p_3$ have degree $1$. Then the previous identity shows that $\lambda p_1(X_2)$, and hence
$p_1(X_2)$, has degree at most $1$. Finally $A_2=X_2+p_2(X_1+p_1(X_2))+p_4(A_1)$ is linear, so
$p_4(A_1)$ has degree at most $1$, and since $A_1$ has degree $1$ in $X_1$, so has $p_4$. Thus
all four factors are affine maps of $R^2$ whose linear parts are elementary matrices, and
$M$, the linear part of $A$, lies in the subgroup $\mathrm{E}_2(R)$ of $\mathrm{SL}_2(R)$
generated by elementary matrices. This contradicts the theorem of Cohn that the transpose
of $M$ does not lie in $\mathrm{E}_2(R)$ \cite{Coh66}, since $\mathrm{E}_2(R)$ is closed
under transposition. The sequence $(2,1,2,1)$ is treated in the same way, with the roles
of $X_1$ and $X_2$ exchanged. Substituting $X_2-p_1(X_1)$ for $X_2$ in
$A_2=X_3^2X_1+(1-X_3X_4)X_2$ gives $p_3(p_2(X_2))=-X_3X_4X_2+\mathrm{const}$, and the rest is
the same.

\emph{The sequences $(c,1,2,c)$ and $(c,2,1,c)$.} The index at the position $2$ occurs
exactly once, which is excluded by Lemma~\ref{lem:anick-tools}(iii).
\end{proof}

The proof uses no assumption on $\operatorname{char}K$. Neither the tails of the factors
nor their degrees are bounded a priori, so the argument cannot be replaced by a finite
search.

\begin{proposition}\label{prop:anick-six}
No word of length at most $5$ in $\cE(K,4)$ represents $A$.
\end{proposition}

\begin{proof}
Let $A=Q_5\circ\dots\circ Q_1$, with $Q_i$ changing $X_{m_i}$ by $h_i$. By
Propositions~\ref{prop:anick}(ii) and~\ref{prop:anick-five}, $A$ has no word of length at
most $4$. As in the proof of Lemma~\ref{lem:anick-sequences}, by Lemma~\ref{lem:word} adjacent
directions differ,
$1$ and $2$ occur, each of $3,4$ occurs never or at least twice, and $m_1$ and $m_5$
occur at least twice. So at most one index $c\in\{3,4\}$ occurs. By
Lemma~\ref{lem:anick-tools}(iii), applied to the given word and to the reversed word of
Lemma~\ref{lem:anick-tools}(i), an index $j\in\{1,2\}$ which occurs exactly once does so
at the position $3$. If $c$ occurs three times, then at the positions $1,3,5$, and $1$ and
$2$ occur once each, at the positions $2$ and $4$. Hence either no $c$ occurs and the
sequence alternates, or $c$ occurs twice, one index $j\in\{1,2\}$ twice and the other
index $k$ once, at the position $3$. In the second case the sequence is one of
\[
(c,j,k,c,j),\quad(j,c,k,j,c),\quad(c,j,k,j,c),\quad(j,c,k,c,j).
\]
The first two are exchanged by reversal, and by Lemma~\ref{lem:anick-tools}(i) we may
assume $c=3$ in the first and the third. We also use the following observation. If
$m_{i-1}=m_{i+1}$ and $h_i$ does not involve $X_{m_{i-1}}$, then
$Q_i\circ Q_{i-1}\circ Q_i^{-1}$ is elementary in the direction $m_{i-1}$, and
$Q_{i+1}\circ Q_i\circ Q_{i-1}=\bigl(Q_{i+1}\circ Q_i\circ Q_{i-1}\circ Q_i^{-1}\bigr)\circ
Q_i$ gives a word of length at most $4$. So in this situation $h_i$ involves
$X_{m_{i-1}}$.

\emph{The sequences $(3,j,k,3,j)$.} Only $Q_3$ changes $X_k$, so
$h_3=(A_k-X_k)\circ(Q_2\circ Q_1)^{-1}$, and $h_3$ does not involve $X_k$. Let $(j,k)=(1,2)$.
The map $(Q_2\circ Q_1)^{-1}$ replaces $X_1$ by $P_1=X_1-h_2(X_2,X_3,X_4)$ and $X_3$ by
$P_3=X_3-h_1(P_1,X_2,X_4)$, so $h_3=P_3\,(P_1P_3-X_2X_4)$. Both factors are nonzero,
since $h_3\ne0$, and degrees in $X_2$ add, so neither involves $X_2$. Then $P_1$ has
degree $1$ in $X_2$, with a coefficient $-\eta\in K[X_3,X_4]$, and comparing the
coefficients of $X_2$ in $P_1P_3-X_2X_4$ gives $\eta P_3=-X_4$. So $\eta$ divides $X_4$ in
$K[X_3,X_4]$, hence $P_3\in K[X_4]$, although $P_3$ and $X_4$ are components of the
automorphism $(Q_2\circ Q_1)^{-1}$ and therefore algebraically independent. Let $(j,k)=(2,1)$. Now $(Q_2\circ Q_1)^{-1}$ replaces $X_2$ by
$P_2=X_2-h_2(X_1,X_3,X_4)$ and $X_3$ by $P_3=X_3-h_1(X_1,P_2,X_4)$, and
$h_3=X_4(X_1P_3-X_4P_2)$. So $X_1P_3-X_4P_2=w$ does not involve $X_1$. Putting $X_1=0$
and subtracting gives $X_1P_3=-X_4\bigl(h_2-h_2|_{X_1=0}\bigr)$, which does not involve
$X_2$. Since $X_2\mapsto P_2$ is an automorphism over $K[X_1,X_3,X_4]$, $P_3$ does not
involve $X_2$ only if $h_1$ does not involve its second argument, and then
$X_1\bigl(X_3-h_1(X_1,X_4)\bigr)$ is divisible by $X_4$, which is false, since its
coefficient of $X_3$ is $X_1$.

\emph{The sequences $(3,j,k,j,3)$.} Only $Q_1$ and $Q_5$ change $X_3$. With
$G=Q_4\circ Q_3\circ Q_2\circ Q_1$ one has $G_i=A_i$ for $i\ne3$, and comparing the third
coordinates gives $h_1+h_5(A_1,A_2,X_4)=0$ with $h_1,h_5\in K[X_1,X_2,X_4]$. Composing
with $A^{-1}$ gives $h_1(P,Q,X_4)=-h_5(X_1,X_2,X_4)$, so $h_1(P,Q,X_4)$ does not
involve $X_3$. By Lemma~\ref{lem:anick-tools}(ii), $h_1$ does not involve $X_1$ and
$X_2$, so $h_1=r(X_4)$ for some $r\in K[t]$. Then $h_1(P,Q,X_4)=r(X_4)$, and
$h_5=-h_1$. Thus $Q_5=Q_1^{-1}$, and
$A'=Q_1\circ A\circ Q_1^{-1}=Q_4\circ Q_3\circ Q_2$. Here
$A'=(X_1+\delta'X_4,\,X_2+\delta'X_3',\,X_3,\,X_4)$ with $X_3'=X_3-r(X_4)$ and
$\delta'=X_1X_3'-X_2X_4$, and $Q_2,Q_3,Q_4$ fix $X_3$ and $X_4$. Let $R=K[X_3,X_4]$. For
$(j,k)=(1,2)$, let $Q_2$ change $X_1$ by $p(X_2)$ and $Q_3$ change $X_2$ by $q(X_1)$, with
$p,q\in R[t]$. Then $A'_2=X_2+q(X_1+p(X_2))$, and substituting $X_1-p(X_2)$ for $X_1$
gives $q(X_1)=(X_1-p(X_2))X_3'^2-X_2X_3'X_4$. The coefficient of $X_2$ on the right is
$-X_3'(p_1X_3'+X_4)$, where $p_1\in R$ is the coefficient of $t$ in $p$. It vanishes only
if $X_3'$ divides $X_4$ in $R$, which is false, since $X_3'$ has degree $1$ in $X_3$. For
$(j,k)=(2,1)$ one obtains in the same way $q(X_2)=X_1X_3'X_4-(X_2-p(X_1))X_4^2$, whose
coefficient of $X_1$ is $X_4(X_3'+p_1X_4)\ne0$.

\emph{The sequences $(j,c,k,c,j)$.} Only $Q_2$ and $Q_4$ change $X_c$, and $h_3$
involves $X_c$ by the observation above. Comparing the $c$-th coordinates gives
$h_2\circ Q_1+h_4\circ Q_3\circ Q_2\circ Q_1=0$, so $h_4\circ Q_3\circ Q_2=-h_2$.
Replacing $X_c$ by $X_c-h_2$ gives
\[
h_4\bigl(X_j,\,X_k+h_3(X_j,X_c,X_{c'}),\,X_{c'}\bigr)=-h_2(X_j,X_k,X_{c'}),
\]
with the variables of $h_4$ written in the order $X_j,X_k,X_{c'}$. If $h_4$ has degree
$f\ge1$ in $X_k$, with leading coefficient $a\in K[X_j,X_{c'}]$, and $h_3$ has degree
$e\ge1$ in $X_c$, with leading coefficient $b\in K[X_j,X_{c'}]$, then the left side has
degree $fe$ in $X_c$, with leading coefficient $ab^f\ne0$, while the right side does not
involve $X_c$. So $h_4$ does not involve $X_k$, and
$h_4=-h_2$. Then $Q_4=Q_2^{-1}$, and $Q_4\circ Q_3\circ Q_2$ is elementary in the
direction $k$, since $h_2$ does not involve $X_k$. This gives a word of length $3$.

\emph{The sequences $(1,2,1,2,1)$ and $(2,1,2,1,2)$.} By
Lemma~\ref{lem:anick-tools}(i) it suffices to treat $(1,2,1,2,1)$. All factors fix
$X_3$ and $X_4$. With $R$ and $M$ as in the proof of Proposition~\ref{prop:anick-five},
let $Q_1,Q_3,Q_5$ change $X_1$ by $p_1(X_2),p_3(X_2),p_5(X_2)$ and $Q_2,Q_4$ change $X_2$
by $p_2(X_1),p_4(X_1)$, with $p_i\in R[t]$ of degree $d_i$. By the observation above,
$d_2,d_3,d_4\ge1$. Let $G=Q_4\circ\dots\circ Q_1$. Then $G_2=A_2$, $G_1=A_1-p_5(A_2)$, and
\[
G_1=X_1+p_1(X_2)+p_3\bigl(X_2+p_2(X_1+p_1(X_2))\bigr),\qquad
A_2=X_2+p_2(X_1+p_1(X_2))+p_4(G_1).
\]
Substituting $X_1-p_1(X_2)$ for $X_1$, putting $X_2=0$ and writing $u=X_1-p_1(0)$ gives
\[
\text{(a)}\ \ p_3(p_2(X_1))=X_3X_4X_1-\lambda p_1(0)-p_5(X_3^2u),\qquad
\text{(b)}\ \ p_4\bigl(\lambda u-p_5(X_3^2u)\bigr)=X_3^2u-p_2(X_1).
\]
Let $a\in R$ be the coefficient of $t$ in $p_5$. The right side of (a) has degree
$\max(1,d_5)$ in $X_1$, since for $d_5\le1$ its coefficient of $X_1$ is
$X_3(X_4-aX_3)\ne0$. The argument of $p_4$ in (b) has degree $\max(1,d_5)$, since for
$d_5\le1$ its coefficient of $X_1$ is $\lambda-aX_3^2\ne0$. Comparing the degrees in
$X_1$ in (a) and (b) gives
\[
d_2d_3=\max(1,d_5),\qquad d_4\max(1,d_5)\le\max(1,d_2).
\]
If $d_5\ge2$, then $d_2\ge1$ and $d_4d_5\le d_2=d_5/d_3$, so $d_3=d_4=1$. If $d_5\le1$,
then $d_2d_3=1$, so $d_2=d_3=1$, and $d_4\le1$. In both cases $d_3=d_4=1$. By Lemma~\ref{lem:anick-tools}(i),
the reversed word gives the same relations with $d_i$ replaced by $d_{6-i}$, so
$d_2=1$, and then $d_5\le1$ and $d_1\le1$. Thus all factors are affine maps of $R^2$ with
elementary linear parts, and $M\in\mathrm{E}_2(R)$, which contradicts the theorem of
Cohn as in the proof of Proposition~\ref{prop:anick-five}.
\end{proof}

The proof uses no assumption on $\operatorname{char}K$.

\begin{remark}\label{rem:anick}
The maps studied in this paper fall into three cases with respect to Theorem~\ref{thm:min-decomp}.
Vasyunin's map satisfies $\ell(F)=\#\supp(F)$, the maps of He and Yan satisfy
$\ell(F)=\#\supp(F)+1$, and $\ell(A)\ge\#\supp(A)+5$ by
Theorem~\ref{thm:anick-seven} below. The first two are not Pascal finite. The map $A$
is a quasi-translation, like its elementary factors, and needs at least seven of them. In characteristic zero
$A=\exp(\tfrac\delta2D_+)\circ\exp(\tfrac\delta2D_-)$ for two commuting triangular locally
nilpotent derivations $D_\pm$ annihilating $\delta$ \cite[Section~4]{MS}, and each factor
is obtained from a commutator of an elementary map with an exponential of the form
$\exp(sD_\pm)$, where $s$ is one of the variables $X_1,\dots,X_4$
\cite[Section~3]{MS}. In characteristic zero the maps $T_h$ and $R_k$ of the proof of (iv)
are the exponentials of the locally nilpotent derivations $h(z\partial_v+2v\partial_x)$
and $k(z\partial_u-2u\partial_x)$. In (iii) the same
commutator is formed with the additional coordinate $X_5$. The two constructions are
polynomial analogues of the Steinberg relations of Remark~\ref{rem:steinberg}. Over every
field of characteristic $\ne2$ the map $\tilde A$ of the proof of (iv) is moreover the
commutator of the product $P$ of the first eight factors with the signed permutation
$S=(-x,u,z,v)$, which satisfies $S\circ S=\Id$ \cite[Section~3]{MS}. So $\tilde A$ is the
product of $P$ and of the conjugate of $P^{-1}$ by $S$. Conjugation by $S$ preserves the set
of elementary maps, since $S$ is linear with a monomial matrix, and the shortened word of
the proof of (iv) consists of a word of six factors for $P$ and a word of six factors for
the conjugate.
\end{remark}

\subsection{Words of length six}\label{ssec:anick-six}

In this subsection we show that $A$ has no decomposition of length $6$
(Theorem~\ref{thm:anick-seven}). The arguments hold over every field.

\begin{lemma}\label{lem:anick-twice}
Let $r\ge0$. In every decomposition of $A\times\Id_r$ into elementary maps each
of the directions $1$ and $2$ occurs at least twice.
\end{lemma}

\begin{proof}
Let $A\times\Id_r=Q_l\circ\dots\circ Q_1$, and suppose that only $Q_i$ changes $X_1$,
by a tail $h$ not involving $X_1$. The map $B=Q_{i-1}\circ\dots\circ Q_1$ fixes $X_1$, so
the first four components of $B^{-1}$ are $(X_1,Y,T,S)$ for some $Y,T,S\in K[X]$, where
$X=(X_1,\dots,X_{4+r})$. The factors after $Q_i$ fix $X_1$, hence
$A_1=X_1+h\circ B$ and
\[
h=(A_1-X_1)\circ B^{-1}=S\,(X_1T-YS).
\]
Since $h\ne0$ does not involve $X_1$ and degrees in $X_1$ add, neither $S$ nor
$X_1T-YS$ involves $X_1$. Putting $X_1=0$ in $X_1T-YS$ gives $X_1T=S(Y-Y_0)$, where
$Y_0$ is obtained from $Y$ by the substitution $X_1=0$. The components of an
automorphism $F$ are irreducible, since a factorization $F_i=PR$ with $P,R$ nonconstant
would give the factorization $X_i=(P\circ F^{-1})(R\circ F^{-1})$ with nonconstant
factors. Two of them are not associate, since they are algebraically independent. So the prime element $S$ of $K[X]$ divides $X_1T$ and is
associate neither to $X_1$ nor to $T$, which is impossible. If only $Q_i$ changes $X_2$,
then the first four components of $B^{-1}$ are $(Y,X_2,T,S)$, and $h=(A_2-X_2)\circ B^{-1}=T\,(YT-X_2S)$. As before
$X_2S=T(Y-Y_0)$ with $Y_0=Y|_{X_2=0}$, and $T$ divides $X_2S$, which is impossible.
\end{proof}

\begin{lemma}\label{lem:anick-six-shape}
Let $A=Q_6\circ\dots\circ Q_1$ with $Q_i\in\cE(K,4)$ changing $X_{m_i}$. Then
$m_i\ne m_{i+1}$ for $1\le i\le5$, and there is $c\in\{3,4\}$ such that each of the
indices $1,2,c$ occurs exactly twice in $(m_1,\dots,m_6)$.
\end{lemma}

\begin{proof}
Adjacent directions differ by Lemma~\ref{lem:word}(ii) and
Proposition~\ref{prop:anick-six}. By Lemma~\ref{lem:anick-twice} the indices $1$ and $2$
occur at least twice, and by Proposition~\ref{prop:app-base}(a) some $c\in\{3,4\}$
occurs at least twice. Since the length is $6$, these bounds are equalities, and the
other index of $\{3,4\}$ does not occur.
\end{proof}

Throughout the rest of this subsection we assume that $A=Q_6\circ\dots\circ Q_1$ is as
in Lemma~\ref{lem:anick-six-shape}. Let $c'$ be the other element of $\{3,4\}$, and
write $x=X_1$, $y=X_2$, $s=X_c$ and $t=X_{c'}$. Let $Q$ and $Q'$ be the factors
changing $s$, with $Q$ applied first, and write
\[
A=U\circ Q'\circ V\circ Q\circ W,
\]
where $U$, $V$ and $W$ are the compositions of the factors applied after $Q'$, between
$Q$ and $Q'$, and before $Q$, respectively. Their factors change $x$ or $y$, and their
tails may involve $s$ and $t$. Let $m_V$ be the number of factors of $V$, and let $Q$ and
$Q'$ change $s$ by $h,g\in K[x,y,t]$. Since $V$ fixes $s$ and $t$, we write
$V(x,y,s,t)=(V_s(x,y),s,t)$, where $V_s$ is an automorphism of $K[s,t][x,y]$ over
$K[s,t]$. For a polynomial $r$ we denote by $V_r$ the result of the substitution of $r$
for $s$ in $V_s$.

\begin{lemma}\label{lem:anick-plane-family}
\textup{(a)} $g(V_r(x,y),t)=-h(x,y,t)$ in $K[r,t,x,y]$, where $r$ is a new
indeterminate.
\textup{(b)} $h\notin K[t]$ and $g\notin K[t]$.
\end{lemma}

\begin{proof}
(a) The maps $U$ and $W$ fix $s$, hence so does $Q'\circ V\circ Q=U^{-1}\circ A\circ
W^{-1}$. The map $Q$ replaces $s$ by $s+h$, then $V$ acts on $(x,y)$ by $V_{s+h}$, and
$Q'$ adds $g$ evaluated at the new coordinates. So the $s$-coordinate of $Q'\circ V\circ Q$ is
$s+h+g(V_{s+h}(x,y),t)$, and $g(V_{s+h}(x,y),t)=-h$. Since $h$ does not involve $s$, the
substitution of $s-h$ for $s$ is an automorphism of $K[x,y,t][s]$, and it gives
$g(V_s(x,y),t)=-h$.

(b) By (a), $h\in K[t]$ if and only if $g\in K[t]$, and then $g=-h$, since composition
with $V_s$ maps $K[s,t][x,y]\setminus K[s,t]$ into itself. In this case $Q'\circ V\circ Q$
is the map $V'$ obtained from $V$ by the substitution of $s+h(t)$ for $s$ in the tails,
a composition of $m_V$ elementary maps changing $x$ or $y$. Then $A=U\circ V'\circ W$ is
a word which changes neither $X_3$ nor $X_4$, contrary to
Proposition~\ref{prop:app-base}(a).
\end{proof}

\begin{proposition}\label{prop:anick-gaps}
$m_V\notin\{0,1,4\}$.
\end{proposition}

\begin{proof}
If $m_V=0$, then $Q$ and $Q'$ are adjacent, contrary to
Lemma~\ref{lem:anick-six-shape}.

Let $m_V=1$, and suppose that $V$ changes $x$ by $a(y,s,t)$. The case of $y$ is the same
with $x$ and $y$ exchanged. Put $\phi=a(y,r,t)-a(y,0,t)\in K[y,t][r]$. By
Lemma~\ref{lem:anick-plane-family}(a), $g(x+a(y,r,t),y,t)=-h=g(x+a(y,0,t),y,t)$, and the
substitution of $x-a(y,0,t)$ for $x$ gives $g(x+\phi,y,t)=g(x,y,t)$. Suppose that
$\phi\ne0$. Then $\phi$ is a polynomial of positive degree in $r$ over $K[y,t]$, hence
transcendental over $K(x,y,t)$, and the homomorphism $K(x,y,t)[w]\to K(x,y,t)[r]$,
$w\mapsto\phi$, is injective. It maps $g(x+w,y,t)-g(x,y,t)$ to $0$, so
$g(x+w,y,t)=g(x,y,t)$ in $K[x,y,t,w]$ and $g\in K[y,t]$. Then $h=-g$ by
Lemma~\ref{lem:anick-plane-family}(a), and $Q'\circ V\circ Q$ changes only $x$, by
$a(y,s+h(y,t),t)$, since $Q'$ restores $s$. If $\phi=0$, then $V$ does not depend on
$s$, and $Q'\circ V\circ Q=V$, since $Q'\circ V\circ Q$ fixes $s$. In both cases $A$ is
represented by a word of length $4$, contrary to Proposition~\ref{prop:anick-five}.

Let $m_V=4$, so that $U=W=\Id$ and the direction sequence is $(c,j,k,j,k,c)$ with
$\{j,k\}=\{1,2\}$. Comparing the $s$-coordinates in $A=Q'\circ V\circ Q$ gives
$h+g(A_1,A_2,t)=0$. Composing with $A^{-1}$ gives
$h\bigl((A^{-1})_1,(A^{-1})_2,t\bigr)=-g(x,y,t)$, which does not involve $s$. By Lemma~\ref{lem:anick-tools}(ii), $h\in K[t]$, contrary to
Lemma~\ref{lem:anick-plane-family}(b).
\end{proof}

So $m_V\in\{2,3\}$. By Lemma~\ref{lem:anick-tools}(i) we may assume $c=3$, and replace
the word by its reversal, which does not change $m_V$. The positions of the two
occurrences of $3$ are then $(2,6)$ if $m_V=3$, and the remaining positions carry
$1,1,2,2$ with no two equal neighbors. Hence for $m_V=3$ the direction sequence is
\begin{equation}\label{eq:anick-gap3}
(1,3,2,1,2,3)\quad\text{or}\quad(2,3,1,2,1,3).
\end{equation}
For $m_V=2$ the positions are $(1,4)$, $(2,5)$ or $(3,6)$. This gives $12$ sequences,
forming $6$ pairs under reversal, and every pair contains exactly one of
\begin{equation}\label{eq:anick-gap2}
\begin{gathered}
\text{(a) }(3,2,1,3,2,1),\quad \text{(b) }(3,1,2,3,2,1),\quad
\text{(c) }(3,2,1,3,1,2),\\
\text{(d) }(3,1,2,3,1,2),\quad \text{(e) }(1,3,1,2,3,2),\quad
\text{(f) }(1,3,2,1,3,2).
\end{gathered}
\end{equation}

\subsubsection*{The case $m_V=2$}
Let $c=3$. Then
\[
A_1=\ell_1:=(1+st)x-t^2y,\qquad A_2=\ell_2:=s^2x+(1-st)y .
\]
Let $S=K[s,t][x,y]$, and for $0\ne P\in S$ let $\deg P$ be its total degree in $x,y$ and
$P^{\mathrm{top}}$ its homogeneous component of degree $\deg P$, a form in $x,y$ with
coefficients in $K[s,t]$. We write $\deg_s$ for the degree in $s$ of an element of
$K[s,t]$.

\begin{lemma}\label{lem:anick-top}
Let $b\in K[t][x,w]$ with $\deg_wb\ge1$, $h\in K[t][x,y]$ with $\deg_yh\ge1$, and
$P=y+b(x,s+h)\in S$. If $\deg P\ge1$, then $P^{\mathrm{top}}\in K[t][x,y]$. The same
holds with $x$ and $y$ exchanged.
\end{lemma}

\begin{proof}
Write $h=p+q$, where $p\in K[t][x]$ and $q$ is the sum of the terms of $h$ divisible by
$y$. Then $q\ne0$. Put $\kappa=\deg q\ge1$ and $\eta=q^{\mathrm{top}}$, which is
divisible by $y$. Let $\hat b(x,T)=b(x,p(x)+T)\in K[t][x,T]$, so that
$\deg_T\hat b=\deg_wb\ge1$, and let $\hat B$ be the component of $\hat b$ of the largest
weight $N$ for the weights $1$ of $x$ and $\kappa$ of $T$. Then $N\ge\kappa\ge1$. Put
$R=b(x,s+h)=\hat b(x,s+q)$. Since $(s+q)^j$ has degree $\kappa j$ and top component
$\eta^j$, a term $\gamma x^iT^j$ of $\hat b$, $\gamma\in K[t]$, contributes to $R$ a
polynomial of degree at most $i+\kappa j$ whose component of degree $i+\kappa j$ is
$\gamma x^i\eta^j$. So $R=\hat B(x,\eta)+R'$ with $\deg R'<N$ or $R'=0$. Since
$\eta$ involves $y$, the polynomials $x$ and $\eta$ are algebraically independent
over $K(t)$, and $\hat B(x,\eta)\ne0$. Hence $R^{\mathrm{top}}=\hat B(x,\eta)\in
K[t][x,y]$ and $\deg R=N$. If $N\ge2$, then $P^{\mathrm{top}}=R^{\mathrm{top}}$. If $N=1$,
then $P-y-R^{\mathrm{top}}\in K[s,t]$, and $\deg P\ge1$ gives
$P^{\mathrm{top}}=y+R^{\mathrm{top}}$.
\end{proof}

For $\varepsilon\in\{1,-1\}$ let $\pi_\varepsilon\colon K[s,t]\to K[t,t^{-1}]$ be the
homomorphism over $K[t]$ with $\pi_\varepsilon(s)=\varepsilon t^{-1}$. It maps
$1-\varepsilon st$ to $0$ and is the identity on $K[t]$.

\begin{lemma}\label{lem:anick-tz}
\textup{(i)} If $D\in K[s,t]$ divides $st$, then $\pi_\varepsilon(D)=\gamma t^j$ with
$\gamma\in K^*$ and $j\in\{-1,0,1\}$.
\textup{(ii)} For every $\gamma\in K[s,t]$ the polynomials $1+st-\gamma s^2$,
$1-st+\gamma t^2$ and $s^2-\gamma(1+st)$ do not lie in $K[t]$, and the first one has
positive degree in $s$.
\end{lemma}

\begin{proof}
(i) Since $s$ and $t$ are prime in $K[s,t]$, $D=\gamma s^at^b$ with $\gamma\in K^*$ and
$a,b\in\{0,1\}$. (ii) If $\gamma\ne0$, then $\deg_s(\gamma s^2)\ge2$, so
$\deg_s(1+st-\gamma s^2)\ge1$, and this also holds for $\gamma=0$. The coefficient of $s$
in $1-st+\gamma t^2$ is $-t+t^2\gamma_1$, where $\gamma_1\in K[t]$ is the coefficient of
$s$ in $\gamma$, and it is nonzero. Finally $\pi_{-1}(s^2-\gamma(1+st))=t^{-2}\notin
K[t]$.
\end{proof}

We write the tails of the factors of $U$, $V$ and $W$ as follows. A factor of $V$ changes
$y$ by $b(x,s)$ or $x$ by $a(y,s)$, where $a\in K[t][y,w]$ and $b\in K[t][x,w]$. A factor
of $U$ changes $y$ by $e(x)$ or $x$ by $f(y)$, and a factor of $W$ changes $x$ by
$k(y)$, where $e\in K[s,t][x]$ and $f,k\in K[s,t][y]$. We write
$e=\sum_ie_ix^i$ with $e_i\in K[s,t]$, $\varepsilon=\deg_xe$, and similarly
$f=\sum_if_iy^i$, $\varphi=\deg_yf$, $k=\sum_ik_iy^i$, $\kappa=\deg_yk$, with the
convention that the degree of an element of $K[s,t]$ is $0$. Derivatives $e'$, $f'$,
$k'$ are taken with respect to $x$, $y$, $y$. Since $S$ is a domain, the total degree in $x,y$ is additive on
products. Recall that for $E\in K[s,t][T]$ and
$X\in S$ with $\deg X\ge1$ one has $\deg E(X)=\deg_TE\cdot\deg X$, and if $\deg_TE\ge1$,
then $E(X)^{\mathrm{top}}$ is the leading coefficient of $E$ times
$(X^{\mathrm{top}})^{\deg_TE}$.

\begin{proposition}\label{prop:anick-gap-two}
$m_V\ne2$.
\end{proposition}

\begin{proof}
Let $c=3$ and let the direction sequence be one of \eqref{eq:anick-gap2}.

\emph{The sequences \textup{(a)} to \textup{(d)}.} Here $W=\Id$, and
$\Theta=U^{-1}\circ A=Q'\circ V\circ Q$ fixes $s$ and $t$. As in the proof of
Lemma~\ref{lem:anick-plane-family}, the coordinates $(\Theta_1,\Theta_2)$ are
$V_{s+h}(x,y)$. The first factor of $V$ changes $y$ by $b$ in (a) and (c), and $x$ by $a$
in (b) and (d). In (a) and (b), $U$ changes $y$ by $e$ and then $x$ by $f$, and in (c)
and (d) it changes $x$ by $f$ and then $y$ by $e$. Computing $U^{-1}\circ A$ we obtain
the identities
\[
\begin{aligned}
\text{(a)}&\quad L:=\ell_2-e(X)=y+b(x,s+h),&& X=\ell_1-f(\ell_2),\\
\text{(b)}&\quad L:=\ell_1-f(\ell_2)=x+a(y,s+h),\\
\text{(c)}&\quad L:=\ell_2-e(\ell_1)=y+b(x,s+h),\\
\text{(d)}&\quad L:=\ell_1-f(Y)=x+a(y,s+h),&& Y=\ell_2-e(\ell_1).
\end{aligned}
\]
We show below that in each case
\begin{enumerate}[label=\textup{(\arabic*)},leftmargin=2em]
\item $\deg L\ge1$ and $L^{\mathrm{top}}\notin K[t][x,y]$, and
\item $\partial_yL\ne1$ in \textup{(a)} and \textup{(c)}, and $\partial_xL\ne1$ in
\textup{(b)} and \textup{(d)}.
\end{enumerate}
This gives a contradiction. In (a) and (c), if $\deg_wb\ge1$ and $\deg_yh\ge1$, then
Lemma~\ref{lem:anick-top} contradicts (1). Otherwise $y+b(x,s+h)\in y+K[s,t][x]$, so
$\partial_yL=1$, contrary to (2). The cases (b) and (d) are the same with $x$ and $y$
exchanged.

(a) If $\varphi\ge2$, then $X^{\mathrm{top}}=-f_\varphi\ell_2^\varphi$. If $\varphi\le1$,
then $X=\lambda-f_0$ with
\[
\lambda=(1+st-f_1s^2)\,x-\bigl(t^2+f_1(1-st)\bigr)\,y\ne0
\]
by Lemma~\ref{lem:anick-tz}(ii). In particular $\deg X\ge1$. If $\varepsilon=0$, then
$L^{\mathrm{top}}=\ell_2$, whose coefficient of $x$ is $s^2$. Let $\varepsilon\ge1$. If
$\varphi\ge2$, then
$L^{\mathrm{top}}=-e_\varepsilon(-f_\varphi)^\varepsilon\ell_2^{\varepsilon\varphi}$,
whose coefficient of $x^{\varepsilon\varphi}$ is divisible by $s^{2\varepsilon\varphi}$.
Let $\varphi\le1$. If $\varepsilon\ge2$, then
$L^{\mathrm{top}}=-e_\varepsilon\lambda^\varepsilon$, and its coefficient of
$x^\varepsilon$ has positive degree in $s$ by Lemma~\ref{lem:anick-tz}(ii). Let
$\varepsilon=1$, and put $E=e_1$, $F=f_1$, $\Pi=1+EF$. The linear part of $L$ is
$\ell_2-E\lambda=c_1x+c_2y$ with
\[
c_1=\Pi s^2-E(1+st),\qquad c_2=\Pi(1-st)+Et^2 .
\]
This linear system for $(\Pi,E)$ has determinant $s^2t^2+(1+st)(1-st)=1$, and
\[
\Pi=t^2c_1+(1+st)c_2,\qquad E=s^2c_2-(1-st)c_1 .
\]
Suppose that $c_1,c_2\in K[t]$. If $c_2\ne0$, then $\deg_sE=2>\deg_s(\Pi-1)$, so
$EF=\Pi-1$ gives $F=0$ and $\Pi=1$, while the coefficient of $s$ in $\Pi$ is $tc_2\ne0$.
If $c_2=0$ and $c_1\ne0$, then $\deg_sE=1>0=\deg_s(\Pi-1)$, so $\Pi=1$, that is
$t^2c_1=1$, which is impossible. If $c_1=c_2=0$, then $E=0$ and $\Pi=0$, while $\Pi=1$
for $E=0$. So $\deg L=1$ and $L^{\mathrm{top}}\notin K[t][x,y]$. This proves (1). For (2),
\[
\partial_yL=(1-st)+e'(X)\bigl(t^2+f'(\ell_2)(1-st)\bigr).
\]
If $\partial_yL=1$, then $e'(X)\bigl(t^2+f'(\ell_2)(1-st)\bigr)=st$. The right side has
total degree $0$ in $x,y$, so both factors have total degree $0$ in $x,y$. Since
$\deg X\ge1$ and $\deg\ell_2=1$, it follows that $e',f'\in K[s,t]$. Then $D=t^2+f'(1-st)$ divides $st$,
and $\pi_1(D)=t^2$, contrary to Lemma~\ref{lem:anick-tz}(i).

(b) If $\varphi\ge2$, then $L^{\mathrm{top}}=-f_\varphi\ell_2^\varphi$, whose coefficient of
$x^\varphi$ is divisible by $s^{2\varphi}$. If $\varphi\le1$, then $\deg L=1$, since the
coefficient of $x$ in $L$ is $1+st-f_1s^2$, of positive degree in $s$. If $\partial_xL=1+st-f'(\ell_2)s^2=1$,
then $f'(\ell_2)s^2=st$, so $f'(\ell_2)$ has total degree $0$ in $x,y$, hence
$f'\in K[s,t]$ and $f's=t$,
which is impossible.

(c) If $\varepsilon\ge2$, then $L^{\mathrm{top}}=-e_\varepsilon\ell_1^\varepsilon$, whose
coefficient of $x^\varepsilon$ is $-e_\varepsilon(1+st)^\varepsilon$. If
$\varepsilon\le1$, then $\deg L=1$, since the coefficient of $y$ in $L$ is
$1-st+t^2e_1\notin K[t]$. If
$\partial_yL=1-st+t^2e'(\ell_1)=1$, then $te'(\ell_1)=s$ has total degree $0$ in $x,y$,
so $e'\in K[s,t]$ and $te'=s$, which is
impossible.

(d) If $\varepsilon\ge2$, then $Y^{\mathrm{top}}=-e_\varepsilon\ell_1^\varepsilon$. If
$\varepsilon\le1$, then $Y=\mu-e_0$ with
\[
\mu=\bigl(s^2-e_1(1+st)\bigr)\,x+\bigl(1-st+t^2e_1\bigr)\,y\ne0 .
\]
So $\deg Y\ge1$. If $\varphi=0$, the coefficient of $x$ in $L^{\mathrm{top}}=\ell_1$ is
$1+st$. Let $\varphi\ge1$. If $\varepsilon\ge2$, then
$L^{\mathrm{top}}=-f_\varphi(-e_\varepsilon)^\varphi\ell_1^{\varepsilon\varphi}$, whose
coefficient of $x^{\varepsilon\varphi}$ is divisible by $(1+st)^{\varepsilon\varphi}$. If
$\varepsilon\le1$ and $\varphi\ge2$, then $L^{\mathrm{top}}=-f_\varphi\mu^\varphi$, and
its coefficient of $y^\varphi$, $-f_\varphi(1-st+t^2e_1)^\varphi$, has positive degree in
$s$ by Lemma~\ref{lem:anick-tz}(ii). Let $\varepsilon\le1$ and $\varphi=1$, and put
$E=e_1$, $F=f_1\ne0$ and $\Pi=1+FE$. The linear part of $L$ is $\ell_1-F\mu=c_1x+c_2y$
with
\[
c_1=(1+st)\Pi-s^2F,\qquad c_2=-t^2\Pi-(1-st)F .
\]
The determinant of this system for $(\Pi,F)$ is $-1$, and
\[
\Pi=(1-st)c_1-s^2c_2,\qquad F=-t^2c_1-(1+st)c_2 .
\]
In particular $(c_1,c_2)\ne(0,0)$, since $F\ne0$. Suppose that $c_1,c_2\in K[t]$. If
$c_2\ne0$, then $F$ has degree $1$ in $s$ with leading coefficient $-tc_2$, and $\Pi-1$
has degree $2$ in $s$ with leading coefficient $-c_2$. If $c_2=0$, then $F=-t^2c_1$ and
$\Pi-1=c_1-1-c_1st$ has degree $1$ in $s$ with leading coefficient $-tc_1$. In both
cases $EF=\Pi-1$ forces $E$ to have degree $1$ in $s$ with leading coefficient $t^{-1}$,
which is impossible. This proves (1). For (2),
\[
\partial_xL=1+st-f'(Y)\bigl(s^2-e'(\ell_1)(1+st)\bigr).
\]
If $\partial_xL=1$, then
$f'(Y)\bigl(s^2-e'(\ell_1)(1+st)\bigr)=st$, so both factors have total degree $0$ in
$x,y$. Since $\deg Y\ge1$ and $\deg\ell_1=1$, this gives $e',f'\in K[s,t]$, and
$D=s^2-e'(1+st)$ divides $st$. But $\pi_{-1}(D)=t^{-2}$, contrary to
Lemma~\ref{lem:anick-tz}(i).

\emph{The sequences \textup{(e)} and \textup{(f)}.} Here $W$ changes $x$ by $k(y)$, and
$U$ changes $y$ by $e(x)$. Put $x_1=x+k(y)$. The substitution $x\mapsto x_1$ is an
automorphism of $S$ over $K[s,t][y]$, and we compute in $S=K[s,t][x_1,y]$, with degrees
and top components with respect to $x_1,y$. In these variables
\[
\ell_1=(1+st)(x_1-k)-t^2y,\qquad \ell_2=s^2(x_1-k)+(1-st)y,
\]
and after $W$ and $Q$ the coordinates are $(x_1,y,s+h(x_1,y),t)$.

(e) The first factor of $V$ changes $x$ by $a$, the second changes $y$, and $U$ and $Q'$
do not change $x$. Comparing the first coordinates gives
\[
x_1+a\bigl(y,s+h(x_1,y)\bigr)=\ell_1=(1+st)x_1-(1+st)k-t^2y=:L .
\]
If $\kappa\ge2$, then $L^{\mathrm{top}}=-(1+st)k_\kappa y^\kappa$. If $\kappa\le1$, the
coefficient of $x_1$ in $L$ is $1+st$. In both cases $\deg L\ge1$ and
$L^{\mathrm{top}}\notin K[t][x_1,y]$, and $\partial_{x_1}L=1+st\ne1$. As above,
Lemma~\ref{lem:anick-top}, with $x_1$ and $y$ in the roles of $y$ and $x$, gives a
contradiction if $\deg_wa\ge1$ and $\deg_{x_1}h\ge1$, and otherwise
$\partial_{x_1}L=1$, a contradiction.

(f) The first factor of $V$ changes $y$ by $b$, the second changes $x$, and then $Q'$ and
$U$ follow. Since $U$ does not change $x$, the first coordinate after the second factor
of $V$ is $A_1=\ell_1$, and comparing the second coordinates gives
\[
y+b\bigl(x_1,s+h(x_1,y)\bigr)=\ell_2-e(\ell_1)=:L .
\]
If $\kappa\ge2$, then $\ell_1^{\mathrm{top}}=-(1+st)k_\kappa y^\kappa$ and
$\ell_2^{\mathrm{top}}=-s^2k_\kappa y^\kappa$. In this case, if $\varepsilon\ge2$, then
$L^{\mathrm{top}}=-e_\varepsilon\bigl(-(1+st)k_\kappa\bigr)^\varepsilon
y^{\varepsilon\kappa}$, if $\varepsilon=1$, then $L^{\mathrm{top}}=k_\kappa\bigl(e_1(1+st)-
s^2\bigr)y^\kappa$, and if $\varepsilon=0$, then $L^{\mathrm{top}}=-s^2k_\kappa y^\kappa$.
If $\kappa\le1$, then $\ell_1$ and $\ell_2$ have degree $1$, with coefficients $1+st$ and
$s^2$ of $x_1$. In this case, if $\varepsilon\ge2$, the coefficient of $x_1^\varepsilon$
in $L^{\mathrm{top}}$ is $-e_\varepsilon(1+st)^\varepsilon$, and if $\varepsilon\le1$, the
coefficient of $x_1$ in $L$ is $s^2-e_1(1+st)$. By Lemma~\ref{lem:anick-tz}(ii), in all
cases $\deg L\ge1$ and $L^{\mathrm{top}}\notin K[t][x_1,y]$. For (2),
\[
\partial_yL=1-st-s^2k'+e'(\ell_1)\bigl((1+st)k'+t^2\bigr).
\]
Suppose that $\partial_yL=1$, that is
\[
e'(\ell_1)\bigl((1+st)k'+t^2\bigr)=st+s^2k'
\]
in $K[s,t][x_1,y]$. The second factor on the left lies in $K[s,t][y]$ and is nonzero, since
$1+st$ does not divide $t^2$, and the right side lies in $K[s,t][y]$. Comparing degrees in
$x_1$, the factor $e'(\ell_1)$ does not involve $x_1$. Since $\ell_1$ has degree $1$ in
$x_1$ with coefficient $1+st$, its degree in $x_1$ is the degree of $e'$, so
$e'=\gamma\in K[s,t]$. Then $k'D=st-\gamma t^2$ with $D=\gamma(1+st)-s^2\ne0$, and
comparing degrees in $y$ gives $k'=\omega\in K[s,t]$. Multiplying by $1+st$ and using $\gamma(1+st)=D+s^2$ gives
\[
D\bigl(\omega(1+st)+t^2\bigr)=st .
\]
So $D$ divides $st$, while $\pi_{-1}(D)=-t^{-2}$, contrary to
Lemma~\ref{lem:anick-tz}(i). As in (a), if $\deg_wb\ge1$ and $\deg_yh\ge1$, then
Lemma~\ref{lem:anick-top}, applied in the variables $x_1,y$, contradicts (1), and
otherwise $y+b(x_1,s+h)\in y+K[s,t][x_1]$, so $\partial_yL=1$, a contradiction.
\end{proof}

\subsubsection*{The case $m_V=3$}

\begin{proposition}\label{prop:anick-gap-three}
$m_V\ne3$.
\end{proposition}

\begin{proof}
Let $c=3$ and let the direction sequence be one of \eqref{eq:anick-gap3}. Then $U=\Id$,
and $W$ consists of one factor, which changes $x$ by $a(y,s,t)$ for $(1,3,2,1,2,3)$ and
$y$ by $b(x,s,t)$ for $(2,3,1,2,1,3)$. Write $A(x,y,s,t)=(A_s(x,y),s,t)$ and
$W(x,y,s,t)=(W_s(x,y),s,t)$, where
\[
A_s(x,y)=\bigl((1+st)x-t^2y,\ s^2x+(1-st)y\bigr).
\]
In $A=Q'\circ V\circ Q\circ W$ the map $Q$ adds $h(W_s(x,y),t)$ to $s$, and $Q'$ adds $g$
evaluated at the final coordinates, whose $x,y$-part is $A_s(x,y)$. Since $A$ fixes $s$,
\begin{equation}\label{eq:anick-gap3-id}
g\circ A_s=H\circ W_s,\qquad H=-h,
\end{equation}
where $g,H\in K[t][x,y]$ are regarded as polynomials in $x,y$ over $K[t]$, and
\eqref{eq:anick-gap3-id} holds in $K[s,t][x,y]$. By
Lemma~\ref{lem:anick-plane-family}(b), $d=\deg g\ge1$. Write
$g_d=\sum_{k=0}^dc_kx^{d-k}y^k$ with $c_k\in K[t]$ for the homogeneous component of $g$
of degree $d$.

\emph{The sequence $(1,3,2,1,2,3)$.} Here $W_s(x,y)=(x+a(y,s,t),y)$, so the right side of
\eqref{eq:anick-gap3-id} has degree $\deg_xH$ in $x$, and its leading coefficient in $x$
lies in $K[t][y]$. On the left side the degree in $x$ is at most $d$, and the
coefficient of $x^d$ is
\[
g_d(1+st,s^2)=\sum_{k=0}^dc_k(1+st)^{d-k}s^{2k}.
\]
The summand of index $k$ has degree $d+k$ in $s$, with leading coefficient $c_kt^{d-k}$.
For the largest $k$ with $c_k\ne0$ it follows that $g_d(1+st,s^2)$ has degree $d+k\ge1$ in
$s$. So the leading coefficients in $x$ of the two sides differ, a contradiction.

\emph{The sequence $(2,3,1,2,1,3)$.} Here $W_s(x,y)=(x,y+b(x,s,t))$. Comparing the leading
coefficients in $y$ in the same way, the coefficient of $y^d$ on the left side is
\[
g_d(-t^2,1-st)=\sum_{k=0}^dc_k(-t^2)^{d-k}(1-st)^k,
\]
whose summand of index $k$ has degree $k$ in $s$ with leading coefficient
$c_k(-t^2)^{d-k}(-t)^k$. This coefficient is therefore nonzero, and it lies in $K[t]$
only if $c_k=0$ for $k\ge1$. Hence $g_d=c_0x^d$ with $c_0\ne0$. By
\eqref{eq:anick-gap3-id}, $g\circ A_s\circ W_s^{-1}=H$ does not depend on $s$, so
$g\circ\mathcal M=g$ for
\[
\mathcal M=A_s\circ W_s^{-1}\circ W_0\circ A_0^{-1},\qquad A_0^{-1}(x,y)=(x+t^2y,\,y).
\]
Put $w=x+t^2y$ and $\psi(w,s)=b(w,s,t)-b(w,0,t)$. A direct computation gives
\[
\mathcal M=\bigl((1+st)x+st^3y+t^2\psi(w,s),\ s^2x+(1-st+s^2t^2)y-(1-st)\psi(w,s)\bigr).
\]
Let $\beta$ be the degree of $\psi$ in $w$ over $K[s,t]$. If $\beta\ge2$, then both
components of $\mathcal M$ have degree at most $\beta$, the first one exactly $\beta$,
so $g_d\circ\mathcal M=c_0\mathcal M_1^d$ has degree $\beta d$, and the other homogeneous
components $g_k$ of $g$ give $\deg(g_k\circ\mathcal M)\le k\beta<\beta d$. Then
$\deg(g\circ\mathcal M)=\beta d>d=\deg g$, a contradiction. If $\beta\le1$, write
$\psi=\psi_0+\psi_1w$ with $\psi_0,\psi_1\in K[s,t]$. Then $\mathcal M$ is affine, and
the linear part of $\mathcal M_1$ is
$(1+st+t^2\psi_1)\,x+(st^3+t^4\psi_1)\,y$. The homogeneous component of degree $d$ of
$g\circ\mathcal M=g$ is $c_0$ times the $d$-th power of this linear form, and it equals
$c_0x^d$. So $st^3+t^4\psi_1=0$, that is $\psi_1=-s/t\notin K[s,t]$, a contradiction.
\end{proof}

\begin{theorem}\label{thm:anick-seven}
No word of length at most $6$ in $\cE(K,4)$ represents $A$. Hence $\ell(A)\ge7$ for every
field $K$ over which $A\in\EE(K,4)$. Moreover $\ell(A\times\Id_r)=6$ for every
field $K$ and every $r\ge1$, and every decomposition of $A\times\Id_r$ of length $6$
changes $X_1$, $X_2$ and exactly one coordinate $X_{4+i}$, $1\le i\le r$, each of them
exactly twice. If $\operatorname{char}K\ne2$, then
\[
7\le\ell(A)\le13,\qquad \ell(A\times\Id_1)=6<\ell(A).
\]
\end{theorem}

\begin{proof}
Words of length at most $5$ are excluded by Proposition~\ref{prop:anick-six}, and words
of length $6$ by Lemma~\ref{lem:anick-six-shape} and
Propositions~\ref{prop:anick-gaps}, \ref{prop:anick-gap-two}
and~\ref{prop:anick-gap-three}. The upper bounds are Propositions~\ref{prop:anick}(iii)
and~\ref{prop:app-13}. Let $A\times\Id_r=Q_l\circ\dots\circ Q_1$ with
$Q_i\in\cE(K,4+r)$ and $l\le6$. If no factor changes one of $X_5,\dots,X_{4+r}$, the
specialization of the proof of Proposition~\ref{prop:stab-bound} gives a word of length
at most $6$ in $\cE(K,4)$ representing $A$, which was just excluded. So some
$X_{4+i}$ is changed, and then at least twice by Lemma~\ref{lem:word}(i). By
Lemma~\ref{lem:anick-twice}, at least two factors change $X_1$ and at least two change
$X_2$. Hence $l=6$, and these are all the factors. The word of
Proposition~\ref{prop:anick}(iii), with the identity on the last $r-1$ coordinates, has
length $6$.
\end{proof}

So the addition of one variable shortens the minimal decompositions of $A$, and
$\ell(A\times\Id_1)<\ell(A)$ if $\operatorname{char}K\ne2$, whereas $\ell(F\times\Id_r)=\ell(F)$ whenever
$\ell(F)\le\#\supp(F)+2$ by Proposition~\ref{prop:stab-bound}. Here
$\ell(A)\ge\#\supp(A)+5$. Theorem~\ref{thmE} follows from
Theorem~\ref{thm:anick-seven}, Lemma~\ref{lem:anick-twice} and
Proposition~\ref{prop:app-base}(a).

\begin{remark}\label{rem:anick-seven}
The proof compares polynomials in $x,y$ over $K[s,t]$ and uses only that the tails of
$Q$ and $Q'$ do not involve $s$, while the coefficients $1\pm st$, $s^2$ and $t^2$ of
$A_1,A_2$ do. For $m_V=3$ the identity \eqref{eq:anick-gap3-id} involves only the outer
factor $W$, and the middle factors $V$ do not enter. For $m_V=2$ the identity involves the
first factor of $V$, and Lemma~\ref{lem:anick-top} is needed. Geometrically, the plane
automorphisms $V_s$ form a family parametrized by $s$, and the identity of
Lemma~\ref{lem:anick-plane-family}(a) says that $V_s$ maps the curve $\{h=-\varrho\}$ onto
the curve $\{g=\varrho\}$ for every $s$ and every $\varrho$. For $m_V=3$ the identity
\eqref{eq:anick-gap3-id} says that $A_s\circ W_s^{-1}$ maps $\{H=\varrho\}$ onto
$\{g=\varrho\}$, and the proof shows that the dependence of $A_s$ on $s$ is incompatible
with this, since $W_s$ is triangular.
\end{remark}

\section{The plane and higher dimensions}\label{sec:plane}

Let $\Aut_KK[X_1,X_2]^2=\Aff(K,2)\ast_{B(K,2)}J(K,2)$ be the Jung--van der Kulk
amalgam, with $J(K,2)$ the de~Jonqui\`eres group and $B(K,2)=\Aff(K,2)\cap J(K,2)$
\cite[Section~5.1]{vdE}. For $F$ in the plane let $\ell_J(F)$ be the number of letters from
$J(K,2)\setminus B(K,2)$ in a reduced alternating word for $F$. It does not depend on
the chosen reduced word \cite[Ch.~I, \S1.2]{Ser}, and it is the number of de~Jonqui\`eres letters of
the canonical normal form of \cite[Theorem~18]{AH_canon}. In terms of the Bass--Serre tree of the amalgam, if $v_0$ is the vertex fixed by
$\Aff(K,2)$, then the distance between $v_0$ and $F(v_0)$ is $2\ell_J(F)$, since a reduced
word with $r$ letters from $J(K,2)\setminus B(K,2)$ gives a path of length $2r$ without
backtracking \cite[Ch.~I, \S4]{Ser}. For automorphisms of $R[X,Y]$, $R$ a domain,
Kuttykrishnan defines a length by the number of elementary factors with coefficients in
the field of fractions of $R$ \cite[Definition~2.1]{Ku09}, and studies the automorphisms
of length three and four in connection with stable tameness \cite{Ku09,Ku}. His main
result assumes that $\mathrm{SL}_m(R[X_1,\dots,X_n])=\mathrm{E}_m(R[X_1,\dots,X_n])$ for
all $n\ge0$ and $m\ge3$ \cite[Theorem~4.1]{Ku09}. For $m=2$ the corresponding equality
fails for $K[X_3,X_4]$ \cite{Coh66}, and this is the obstruction used in
Section~\ref{sec:anick}. If $K$ is algebraically closed, then a plane automorphism $F$,
viewed in $\mathrm{Bir}(\mathbb{P}^2)$ by means of a linear embedding of the plane, has
length $\ell_J(F)$ in the Cremona group in the sense of Blanc and Furter
\cite[Proposition~4.2.2]{BF}.

\begin{lemma}\label{lem:lJ}
\textup{(i)} $\ell_J(F\circ G)\le\ell_J(F)+\ell_J(G)$. \textup{(ii)} For $Q\in\cE(K,2)$,
$\ell_J(Q)=0$ if $Q$ is affine and $\ell_J(Q)=1$ otherwise. \textup{(iii)}
$\ell_J(F)\le\ell(F)$ for $F\in\EE(K,2)$.
\end{lemma}

\begin{proof}
(i) Reducing the concatenation of reduced words for $F$ and $G$ multiplies adjacent
letters of the same factor and creates no letter of $J(K,2)\setminus B(K,2)$. (ii) A
non-affine elementary map changing $X_1$ lies in $J(K,2)\setminus B(K,2)$. A non-affine
elementary map changing $X_2$ is $\rho\circ Q'\circ\rho$ with $\rho=(X_2,X_1)$ and
$Q'\in J(K,2)\setminus B(K,2)$. (iii) follows from (i) and (ii).
\end{proof}

The linear parts are controlled by the following elementary bound, independent of $K$.
Let $e_{ij}(x)$ denote the elementary matrix $I+xE_{ij}$, $i\ne j$.

\begin{lemma}\label{lem:sln}
Let $n\ge2$. Every element of $\mathrm{SL}_n(K)$ is a product of at most $N_n$
elementary matrices, where $N_1=0$ and $N_n=2n+N_{n-1}$. In particular $N_2=4$ and
$N_3=10$.
\end{lemma}

\begin{proof}
Let $A\in\mathrm{SL}_n(K)$. Its first column is nonzero. If $a_{21}=0$, adding to the
second row a row with nonzero first entry makes $a_{21}\ne0$. Adding
$(1-a_{11})/a_{21}$ times the second row to the first makes $a_{11}=1$. Then $n-1$ row
operations clear the rest of the first column and $n-1$ column operations clear the rest
of the first row. These are at most $2n$ multiplications by elementary matrices, and they
give $\mathrm{diag}(1,A')$ with $A'\in\mathrm{SL}_{n-1}(K)$. The inverse of an elementary
matrix is elementary, and induction gives the bound.
\end{proof}

For $n=2$ and $a_{21}\ne0$ three factors suffice, namely
$A=e_{12}\bigl(\tfrac{a_{11}-1}{a_{21}}\bigr)\,e_{21}(a_{21})\,
e_{12}\bigl(\tfrac{a_{22}-1}{a_{21}}\bigr)$.

\begin{proposition}\label{prop:plane-length}
\textup{(i)} $\EE(K,2)=\{F\in\Aut_KK[X_1,X_2]^2:\det J_F=1\}$.
\textup{(ii)} For $F\in\EE(K,2)$, $\ell_J(F)\le\ell(F)\le4\,\ell_J(F)+9$ and
$\ell_J(F)\le\log_2\deg F$.
\textup{(iii)} If $Q,Q'\in\cE(K,2)$ are non-affine and change different coordinates,
then $\ell_J(Q\circ Q'\circ Q^{-1}\circ Q'^{-1})=4$.
\end{proposition}

\begin{proof}
Write $F$ in the normal form of
\cite[Theorem~18]{AH_canon}, with $l=\ell_J(F)$ letters $\Gamma_f=(X_1+f(X_2),X_2)$,
$f(0)=f'(0)=0$, $\deg f\ge2$, from the transversal $\mathcal C$, at most $l+1$ letters
$P(a)=(aX_1+X_2,-X_1)$ from the transversal $\mathcal A$
\cite[Section~5, Proposition~23]{AH_canon}, and a last letter $b_0\in B(K,2)$. Each letter of the
first kind is elementary. A linear map $X\mapsto AX$ with $A$ an elementary matrix is an
elementary map, and the matrix identity
\[
\begin{pmatrix}a&1\\-1&0\end{pmatrix}=
\begin{pmatrix}1&1-a\\0&1\end{pmatrix}
\begin{pmatrix}1&0\\-1&1\end{pmatrix}
\begin{pmatrix}1&1\\0&1\end{pmatrix}
\]
writes $P(a)=(X_1+(1-a)X_2,X_2)\circ(X_1,X_2-X_1)\circ(X_1+X_2,X_2)$ as a product of
at most three elementary maps, the first factor being omitted when $a=1$. All letters other than $b_0$ have Jacobian
determinant $1$, so $\det J_{b_0}=\det J_F$. (i)
Elementary maps have Jacobian determinant $1$. If $\det J_F=1$, the linear part of $b_0$
lies in $\mathrm{SL}_2(K)$, which is generated by elementary matrices, and the
translation part is a product of two elementary translations. (ii) The lower bound is
Lemma~\ref{lem:lJ}. By Lemma~\ref{lem:sln} the linear part of $b_0$ is a product of at
most four elementary matrices, so $\ell(F)\le l+3(l+1)+4+2$. If $l\ge1$, then $\deg F$ is
the product of the degrees $\deg f\ge2$ of the de~Jonqui\`eres letters
\cite[Lemma~9]{AH_canon}, so
$\deg F\ge2^l$. (iii) Up to exchanging $Q$ and $Q'$, which replaces the commutator by its
inverse, $Q=(X_1+g(X_2),X_2)$ and $Q'=\rho\circ(X_1+g'(X_2),X_2)\circ\rho$ with
$\deg g,\deg g'\ge2$. The commutator is the alternating word
$\gamma_1\rho\gamma_2\rho\gamma_3\rho\gamma_4\rho$ with
$\gamma_i\in J(K,2)\setminus B(K,2)$ and $\rho\in\Aff(K,2)\setminus B(K,2)$, which is
reduced.
\end{proof}

\begin{lemma}\label{lem:commutator}
Let $n\ge3$, $a\ne b$, and $f,g\in K[X_i:i\ne a,b]$. Let $C$ change $X_a$ by $X_bf$ and
$C'$ change $X_b$ by $g$. Then $C\circ C'\circ C^{-1}\circ C'^{-1}$ changes $X_a$ by $fg$
and fixes the other coordinates.
\end{lemma}

\begin{proof}
$C'^{-1}$ replaces $X_b$ by $X_b-g$, then $C^{-1}$ subtracts $(X_b-g)f$ from $X_a$, $C'$
restores $X_b$, and $C$ adds $X_bf$ to $X_a$. The polynomials $f$ and $g$ are not changed.
\end{proof}

For $n=3$, $f=X_3$ and $g=X_3^2$ the maps $C=(X_1+X_2X_3,X_2,X_3)$ and
$C'=(X_1,X_2+X_3^2,X_3)$ are non-affine, and their commutator $(X_1+X_3^3,X_2,X_3)$ is
elementary. By Proposition~\ref{prop:plane-length}(iii) this does not happen for $n=2$.

\begin{remark}\label{rem:steinberg}
For constant $f$ and linear $g=cX_k$, $k\ne a,b$, Lemma~\ref{lem:commutator} is the
Steinberg relation $[x_{ab}(f),x_{bk}(c)]=x_{ak}(fc)$ between elementary matrices. So
for $n\ge3$ the group $\EE(K,n)$ satisfies polynomial analogues of the Steinberg
relations, which shorten words. For $n=2$ there is no third index, and
Proposition~\ref{prop:plane-length}(iii) shows that no relation of this kind holds
between non-affine elementary maps. The same mechanism appears when variables are
added. Words for $F\times\Id_r$ may use the new coordinates, which is the setting of
stable tameness \cite{BEW,Ku}. The Nagata automorphism, which is Pascal finite, becomes
tame after adding one variable \cite{Sm}, and the argument of \cite{Sm} uses that it is
the exponential of a locally nilpotent derivation \cite[Introduction]{BEW}. Proposition~\ref{prop:stab-bound} shows that words of
length at most $\#\supp(F)+1$ never use them, and Question~\ref{q:stab} asks whether they
can shorten longer words.
\end{remark}

\begin{theorem}\label{thm:gap}
Let $U=(X_1+X_2^2,X_2)$, $L=(X_1,X_2+X_1^2)$ and $F_m=(U\circ L)^{(m)}$. Then
$\#\supp(F_m)=2$, $\deg F_m=4^m$ and $\ell(F_m)=2m$.
\end{theorem}

\begin{proof}
Since $L=\rho\circ U\circ\rho$, the map $F_m$ has the reduced alternating word
$(U\rho U\rho)^m$, with $2m$ de~Jonqui\`eres letters of degree $2$. Lemma~\ref{lem:lJ}
gives $\ell(F_m)\ge2m$, and $F_m$ is a product of $2m$ elementary maps. The degree is
multiplicative along reduced words \cite[Lemma~9]{AH_canon}. Finally
$\#\supp(F_m)=2$, since a map with $\#\supp\le1$ is $\Id$ or elementary and has
$\ell_J\le1$ by Lemma~\ref{lem:lJ}, while $\ell_J(F_m)=2m\ge2$.
\end{proof}

The same argument applies to the alternating products $W_r=U\circ L\circ U\circ\cdots$ of
$r\ge2$ factors, which satisfy $\#\supp(W_r)=2$ and $\ell(W_r)=r$. So every value $d\ge0$
of the defect $\ell(F)-\#\supp(F)$ occurs already in $\EE(K,2)$.

Theorem~\ref{thmC} follows from Proposition~\ref{prop:plane-length},
Lemma~\ref{lem:commutator} and Theorem~\ref{thm:gap}. By
Proposition~\ref{prop:plane-length}(ii) we have $\ell(F)\le4\log_2\deg F+9$ on
$\EE(K,2)$. The maps $(X_1+X_2^d,X_2)$ show that no lower bound in terms of the degree
holds.

For $F\in\EE(K,n)$, $F\ne\Id$, let $\omega(F)$ be the least value of
$\max_i\deg Q_i$ over all decompositions $F=Q_l\circ\dots\circ Q_1$ with
$Q_i\in\cE(K,n)$, and $\omega_0(F)$ the least value over the decompositions of length
$\ell(F)$. Then $\omega(F)\le\omega_0(F)$, and $\deg F\le\omega_0(F)^{\ell(F)}$, since
$\deg(G\circ G')\le\deg G\cdot\deg G'$.

\begin{proposition}\label{prop:plane-degrees}
Let $F\in\EE(K,2)$ be non-affine, and let $D_J(F)$ be the largest degree of a
de~Jonqui\`eres letter of a reduced alternating word for $F$. Then
$\omega(F)=D_J(F)$.
\end{proposition}

\begin{proof}
For $g=(aX_1+f(X_2),\,bX_2+c)\in J(K,2)$ put $\deg g=\max(1,\deg f)$. Composition on
either side with an element of $B(K,2)$ does not change $\deg g$, and
$\deg(g\circ g')\le\max(\deg g,\deg g')$ for $g,g'\in J(K,2)$. For $G$ in the plane let
$D(G)$ be the largest degree of a letter from $J(K,2)\setminus B(K,2)$ in a reduced
alternating word for $G$, and $D(G)=1$ if there is no such letter. Since the letters of
a reduced word are determined by $G$ up to multiplication by elements of $B(K,2)$
\cite[Ch.~I, \S1.2]{Ser}, $D(G)$ is well defined, and $D(F)=D_J(F)$. The reduction of
the concatenation of reduced words for $G$ and $G'$ multiplies adjacent letters lying in
the same factor of the amalgam. So every letter from $J(K,2)$ of the resulting reduced
word is a product of letters from $J(K,2)$ of the two words and of elements of
$B(K,2)$, and $D(G\circ G')\le\max(D(G),D(G'))$. By the proof of Lemma~\ref{lem:lJ}(ii),
$D(Q)=\deg Q$ for $Q\in\cE(K,2)$. Hence $D_J(F)\le\max_i\deg Q_i$ for every
decomposition of $F$. Conversely, the decomposition of $F$ in the proof of
Proposition~\ref{prop:plane-length} consists of the de~Jonqui\`eres letters of the
normal form of \cite[Theorem~18]{AH_canon}, of degree at most $D_J(F)$, and of affine
elementary maps.
\end{proof}

For Vasyunin's map $\omega(F)=2<\omega_0(F)=3$ by Theorem~\ref{thmB}, since $F$ is not
affine and every decomposition of length $4$ contains the cubic factor $\Qout$.

For $n=3$ and $\operatorname{char}K=0$ an upper bound in terms of the degree follows from
the reduction theory of Shestakov and Umirbaev in the formulation of Kuroda, and for
$\operatorname{char}K\ge7$ from its extension by Kuroda. Put
$|F|=\deg F_1+\deg F_2+\deg F_3$.

\begin{proposition}\label{prop:dim3-bound}
Let $\operatorname{char}K=0$ or $\operatorname{char}K\ge7$. Then
$\ell(F)\le3|F|+4\le9\deg F+4$ for every $F\in\EE(K,3)$.
\end{proposition}

\begin{proof}
In \cite{Kur,Kur26} the degree of an automorphism is its weighted degree for the weight
$(1,1,1)$, that is $|F|$ \cite[(1.5)]{Kur26}, and the elementary automorphisms are the
elements of $\cE(K,3)\cup\{\Id\}$. Their composition is written in the opposite order,
so a reduction of $F$ in the sense of \cite{Kur,Kur26} replaces $F$ by $E\circ F$ in our
notation, that is, it adds to one coordinate of $F$ a polynomial in the other two
\cite[Remark~6.1]{Kur26}. It suffices to show that every $F\in\EE(K,3)$ with $|F|>3$
admits elementary maps $E_1,\dots,E_r$, $r\le3$, with $|E_r\circ\dots\circ E_1\circ
F|<|F|$. Since $|F|\ge3$ with equality exactly for affine $F$, iterating gives at most
$|F|-3$ steps and an affine map $A$, and $F$ is the composition of $A$ with at most
$3(|F|-3)$ elementary maps.

By \cite[Theorem~2.1]{Kur}, see also \cite[Theorem~1.5]{Kur26}, if
$\operatorname{char}K=0$, and by \cite[Theorem~1.6]{Kur26} if $\operatorname{char}K\ge7$,
the automorphism $F$ admits an elementary reduction, so that $r=1$, or a
Shestakov--Umirbaev reduction $G$. In the latter case, after a permutation of the
coordinates, condition (SU1) of \cite[Definition~1.4]{Kur26} gives
$g_1\in f_1+Kf_3^2+Kf_3$, $g_2\in f_2+Kf_3$ and $g_3\in f_3+K[g_1,g_2]$. So $G$ is
obtained from $F$ by three elementary maps, which act as
\[
(f_1,f_2,f_3)\longmapsto(g_1,f_2,f_3)\longmapsto(g_1,g_2,f_3)\longmapsto(g_1,g_2,g_3),
\]
so that $r=3$. Moreover $|G|<|F|$, by \cite[Theorem~4.1]{Kur} if $\operatorname{char}K=0$.
If $\operatorname{char}K\ge7$, the pair $(F,G)$ is a weak Shestakov--Umirbaev pair in the
sense of \cite[Definition~6.3]{Kur26}, since every Shestakov--Umirbaev pair is one, and
$|G|<|F|$ by \cite[Section~6.2, Lemma~1(ii)]{Kur26}.

The iteration ends with an affine map $A\in\EE(K,3)$, so
$\det J_A=1$. Its linear part lies in $\mathrm{SL}_3(K)$ and is a product of at most
$10$ elementary matrices by Lemma~\ref{lem:sln}, and its translation part is a product
of three elementary translations. Hence $\ell(F)\le3(|F|-3)+13=3|F|+4$, and
$|F|\le3\deg F$.
\end{proof}

So for $n=3$ and $\operatorname{char}K\notin\{2,3,5\}$ the length is bounded linearly
in the degree, while for $n=2$ it is bounded logarithmically. We do not know whether a
logarithmic bound holds for $n=3$.

\begin{remark}\label{rem:quartic}
Wang and Yan prove that every Keller map $F\in\mathbb{C}[X_1,X_2,X_3]^3$ with $\deg F\le4$ is a
tame automorphism \cite[Theorem~1.1]{WY_quart}. If moreover $\det J_F=1$, then
$F\in\EE(\mathbb{C},3)$ by Proposition~\ref{prop:Gamma}, and Proposition~\ref{prop:dim3-bound}
gives $\ell(F)\le3|F|+4\le40$, since $|F|\le12$. So the length is bounded on the
set of Keller maps in $\mathbb{C}[X_1,X_2,X_3]^3$ of degree at most $4$ with Jacobian determinant
$1$. The maps of Section~\ref{sec:heyan} with $d=\deg g=1$ are quartic,
and for them $\ell(F)=4$.
\end{remark}

\section{Summary and questions}\label{sec:questions}

For $\ell(F)=\#\supp(F)$ all minimal decompositions of $F$ are obtained from $F^{-1}$
by extraction of left elementary factors (Theorem~\ref{thm:min-decomp}), and for
Vasyunin's map they still differ in the degrees of their factors (Theorem~\ref{thmB}).
For $\ell(F)=\#\supp(F)+1$ one tail remains unknown (Proposition~\ref{prop:length-k1}),
and for the maps of He and Yan it is determined by the extraction up to parameters, with
two move classes (Proposition~\ref{prop:heyan-min}). For the Anick automorphism
$\ell(A)-\#\supp(A)\ge5$ (Theorem~\ref{thm:anick-seven}), and no description of its
minimal decompositions is known. Every decomposition of $A$ changes $X_3$ or $X_4$, one of
them at least twice (Proposition~\ref{prop:app-base}). Relative to $K[X_3,X_4]$ the map
$A$ is wild \cite[p.~343, Example~3]{Coh}, see Lemma~\ref{lem:app-cohn}, relative to $K[X_3]$ it is tame if
$\operatorname{char}K\ne2$ \cite{MS}, and $A\times\Id_1$ has a decomposition of length $6$
whose factors fix $X_3$ and $X_4$ (Proposition~\ref{prop:anick}(iii)).

Table~\ref{tab:compare} compares the plane, after \cite{AH_canon}, with dimension
$n\ge3$. The first four rows list features of the plane normal form, fixed by $F$ up to
the choice of coset representatives of $B(K,2)$, which need not persist for $n\ge3$. The
last four rows compare the length with the degree, with the degrees of the factors, with
commutators and with iteration.

\begin{table}[ht]
\centering
\small
\renewcommand{\arraystretch}{1.2}
\begin{tabular}{>{\raggedright\arraybackslash}p{2.6cm}>{\raggedright\arraybackslash}p{5.6cm}>{\raggedright\arraybackslash}p{5.6cm}}
& $n=2$ & $n\ge3$\\ \hline
decompositions &
reduced word determined by $F$ up to $B(K,2)$, unique after a choice of coset
representatives \cite[Theorem~18]{AH_canon} &
no such normalization known, minimal decompositions not unique, six for Vasyunin's map
(Theorem~\ref{thm:V-minimal}), at most $k!$ of them for $\ell(F)=\#\supp(F)$,
determined by $F^{-1}$ (Theorem~\ref{thm:min-decomp}), infinitely many for the maps of
He and Yan, with $\ell(F)=\#\supp(F)+1$ (Proposition~\ref{prop:heyan-min})\\
invariants &
number of de~Jonqui\`eres letters and multidegree \cite[Proposition~23]{AH_canon} &
$\ell(F)$ is an invariant, different minimal decompositions of $F$ may have
different multisets of factor degrees
(Theorem~\ref{thm:V-minimal})\\
degree &
$\deg F$ is the product of the degrees $D_i$ of the de~Jonqui\`eres letters
\cite[Lemma~9]{AH_canon} &
$\deg F\le\prod_i\deg Q_i$, strict in general, $\deg F=2$ for Vasyunin's map and products
$24$ and $36$ for its minimal decompositions\\
Pascal finite factors &
canonical factorization \cite[Corollary~22]{AH_canon} &
six factorizations of Vasyunin's map, which is not Pascal finite, into factors of
depth $2$ (Theorem~\ref{thm:V-minimal})\\
length and degree &
$\ell(F)\le4\log_2\deg F+9$ (Proposition~\ref{prop:plane-length}) &
$\ell(F)\le9\deg F+4$ for $n=3$, $\operatorname{char}K\notin\{2,3,5\}$
(Proposition~\ref{prop:dim3-bound}), otherwise not known\\
degrees of factors &
$\omega(F)=D_J(F)$ (Proposition~\ref{prop:plane-degrees}), $\omega_0(F)$ not known in
general &
$\omega(F)=2<\omega_0(F)=3$ for Vasyunin's map (Theorem~\ref{thmB}), $\omega_0(F)=2d$
for the maps of He and Yan (Proposition~\ref{prop:heyan-min})\\
commutators &
of non-affine elementary maps in different directions, never elementary
(Proposition~\ref{prop:plane-length}) &
can be elementary (Lemma~\ref{lem:commutator})\\
iteration &
$d\ell(F)>0$ if and only if $\lambda(F)>1$ (Remark~\ref{rem:dyn-length}) &
open (Question~\ref{q:dyn-length}), $5\le\ell(F^{(m)})\le4m$ and $\lambda=\sqrt2$ for Vasyunin's map
(Theorem~\ref{thmD}), $\lambda=2d$ for the maps of He and Yan
(Proposition~\ref{prop:heyan})
\end{tabular}
\caption{Dimension two and dimensions at least three.}\label{tab:compare}
\end{table}

Table~\ref{tab:examples} lists the examples of this paper according to the defect.

\begin{table}[ht]
\centering
\small
\renewcommand{\arraystretch}{1.2}
\begin{tabular}{>{\raggedright\arraybackslash}p{2.6cm}cccc>{\raggedright\arraybackslash}p{5.2cm}}
map & field & $\#\supp$ & $\ell$ & defect & minimal decompositions\\ \hline
Vasyunin's map & $\operatorname{char}K\ne2$ & $4$ & $4$ & $0$ & six, one move class,
$\lambda=\sqrt2$ (Theorems~\ref{thmB}, \ref{thmD})\\
maps of He and Yan & $\operatorname{char}K=0$ & $3$ & $4$ & $1$ & infinitely many if $K$ is
infinite, two move classes, $\lambda=2d$ (Propositions~\ref{prop:heyan},
\ref{prop:heyan-min})\\
Anick automorphism & $\operatorname{char}K\ne2$ & $2$ & $7$ to $13$ & $5$ to $11$ &
unknown, $\ell(A\times\Id_r)=6$, $\lambda=1$ (Theorem~\ref{thm:anick-seven})\\
$W_r=U\circ L\circ\cdots$ & every $K$ & $2$ & $r$ & $r-2$ & $\lambda=4^{r/2}$ for even
$r$ (Theorem~\ref{thm:gap})
\end{tabular}
\caption{The examples ordered by the defect $\ell(F)-\#\supp(F)$.}\label{tab:examples}
\end{table}

\begin{question}\label{q:stab}
Let $n\ge3$. Is $\ell(F\times\Id_{n-2})=\ell(F)$ for $F\in\EE(K,2)$? By
Proposition~\ref{prop:stab-bound} this holds when $\ell(F)\le\#\supp(F)+2$. A word in $\cE(K,n)$
may use the additional coordinates, and Lemma~\ref{lem:commutator} shows that this can
shorten words. For the maps $F_m$ of Theorem~\ref{thm:gap},
$\ell(F_2\times\Id_{n-2})=4$, and for $m\ge3$ we only know
$4\le\ell(F_m\times\Id_{n-2})\le2m$. More generally, with the stable length $\ell_{\mathrm{st}}$ defined below, is
$\ell(F)-\ell_{\mathrm{st}}(F)$ bounded on $\EE(K,n)$ by a constant depending only on $n$? For the Anick automorphism
this difference is at least $1$ by Theorem~\ref{thm:anick-seven}.
\end{question}

For $F\in\EE(K,n)$ the sequence $\ell(F\times\Id_r)$, $r\ge0$, is nonincreasing, and
we call its limit $\ell_{\mathrm{st}}(F)$ the \emph{stable length} of $F$. It is the length
of the image of $F$ in the direct limit of the groups $\EE(K,n)\to\EE(K,n+1)$,
$F\mapsto F\times\Id_1$, with respect to the union of the sets $\cE(K,n)$. By
Proposition~\ref{prop:stab-bound},
$\ell_{\mathrm{st}}(F)\ge\min(\ell(F),\#\supp(F)+2)$, and by
Theorem~\ref{thm:anick-seven}, $\ell_{\mathrm{st}}(A)=6<7\le\ell(A)$. In particular
$\ell(F)-\ell_{\mathrm{st}}(F)\le\max\bigl(0,\ell(F)-\#\supp(F)-2\bigr)$, so the difference
vanishes for Vasyunin's map and the maps of He and Yan. Question~\ref{q:stab}
asks whether $\ell_{\mathrm{st}}(F)=\ell(F)$ on $\EE(K,2)$.

The length depends a priori on the field. For a field extension $K\subseteq L$ and
$F\in\EE(K,n)$ we write $\ell_K(F)$ and $\ell_L(F)$ for the lengths with respect to
$\cE(K,n)$ and $\cE(L,n)$. Clearly $\ell_L(F)\le\ell_K(F)$.

\begin{proposition}\label{prop:field-dependence}
Let $F\in\EE(k,n)$.
\begin{enumerate}[label=\textup{(\roman*)},leftmargin=2em]
\item $\ell_\Omega(F)$ is the same for all algebraically closed fields $\Omega\supseteq k$.
\item Suppose that $F$ has coefficients in $\mathbb{Z}$, and let $F_p$ be its reduction
modulo $p$. Then $\ell_{\overline{\FF}_p}(F_p)\le\ell_{\overline{\QQ}}(F)$ for all but
finitely many primes $p$.
\item For every extension $L\supseteq k$, $\ell_L(F)=\#\supp(F)$ if and only if
$\ell_k(F)=\#\supp(F)$. Hence $\ell_L(F)=\ell_k(F)$ whenever
$\ell_k(F)\le\#\supp(F)+1$.
\end{enumerate}
\end{proposition}

\begin{proof}
(i) Let $\Omega_1,\Omega_2\supseteq k$ be algebraically closed, and let $F$ have a
decomposition of length $l$ over $\Omega_1$. Let $(m_1,\dots,m_l)$ be its direction
sequence and $D$ a bound for the degrees of its tails. Replace the coefficients of the
tails by indeterminates $c$. The coefficients of the composition $Q_l\circ\dots\circ Q_1$
are polynomials in $c$ with coefficients in $k$, and comparing them with those of $F$
gives a finite system of polynomial equations over $k$. It has a solution in $\Omega_1$,
so the ideal $I$ it generates in $k[c]$ is proper. Since $\Omega_2[c]$ is a free
$k[c]$-module, $I\,\Omega_2[c]$ is proper, and by the Nullstellensatz the system has a
solution in $\Omega_2$. Omitting the factors equal to $\Id$ gives a decomposition of $F$
over $\Omega_2$ of length at most $l$. By symmetry $\ell_{\Omega_1}(F)=\ell_{\Omega_2}(F)$.

(ii) The coefficients of a decomposition of $F$ over $\overline{\QQ}$ of length
$l=\ell_{\overline\QQ}(F)$ lie in a ring $O_L[1/N]$, where $O_L$ is the ring of integers
of a number field $L$ and $N\ge1$. For a prime $p\nmid N$ and a prime ideal $\mathfrak p$
of $O_L$ over $p$, the reduction modulo $\mathfrak p$ maps the factors to elementary maps
or to $\Id$ over the finite field $O_L/\mathfrak p$, and their composition to $F_p$.

(iii) The tests of Theorem~\ref{thm:min-decomp} concern polynomials computed from
$F^{-1}$ by composition with the extracted factors, whose tails are such polynomials. If
$F$ is defined over $k$, all of them have coefficients in $k$, and the outcome of each
test does not depend on the field. So the admissible orderings over $L$ and over $k$
coincide. If $\ell_k(F)=\#\supp(F)+1$, then $\#\supp(F)\le\ell_L(F)\le\ell_k(F)$, and
$\ell_L(F)\ne\#\supp(F)$ by the first statement.
\end{proof}

By (iii) the lengths of Vasyunin's map and of the maps of He and Yan do not change under
field extensions. By (i) and Theorem~\ref{thm:anick-seven}, the length of $A$ over an algebraically closed
field depends only on the characteristic, and $\ell_K(A)\ge\ell_{\overline K}(A)\ge7$
for every field $K$.

\begin{question}\label{q:anick}
Let $\operatorname{char}K\ne2$. Then $\ell(A\times\Id_1)=6<7\le\ell(A)\le13$ by
Theorem~\ref{thm:anick-seven}, so the minimal decompositions of $A$ are shortened by the
addition of one variable, compare Question~\ref{q:stab}. What is $\ell(A)$? Can
$\ell_K(A)>\ell_{\overline K}(A)$ hold, and does $\ell_{\overline K}(A)$ depend on the
characteristic of $K$? If $\operatorname{char}K=2$, is $A\in\EE(K,4)$? By
Lemma~\ref{lem:anick-twice}, Lemma~\ref{lem:defect} and Proposition~\ref{prop:app-base}(a),
a decomposition of length $7$ changes $X_1$, $X_2$ and one $X_c$, $c\in\{3,4\}$, with
multiplicities $2,2,3$ in some order, and does not change the other base coordinate.
\end{question}

\begin{remark}\label{rem:degrees}
Let $\omega$ and $\omega_0$ be as in Section~\ref{sec:plane}.
\begin{enumerate}[label=(\arabic*),leftmargin=2em]
\item Conjugation by an affine map $X\mapsto MX+b$, with $M$ a monomial matrix, maps
$\cE(K,n)$ onto itself and preserves the degrees, so it preserves $\ell$, $\omega$ and
$\omega_0$. These are the only affine maps $T$ with $T\circ\cE(K,n)\circ T^{-1}=\cE(K,n)$.
Indeed $T$ conjugates the translation $X_j\mapsto X_j+1$ to the translation by $Me_j$,
which is elementary only if $Me_j$ is a multiple of a unit vector. A general linear conjugation does not preserve $\ell$. For
$F=(X_1+X_2^2,X_2)$ and $L=(X_1,X_2+X_1)$ one has
$L\circ F\circ L^{-1}=\bigl(X_1+(X_2-X_1)^2,\;X_2+(X_2-X_1)^2\bigr)$, with inverse
$\bigl(X_1-(X_2-X_1)^2,\;X_2-(X_2-X_1)^2\bigr)$. Both polynomials
$X_j-(L\circ F^{-1}\circ L^{-1})_j=(X_2-X_1)^2$ involve $X_j$, so
$\ell(L\circ F\circ L^{-1})=3$ by Theorem~\ref{thm:min-decomp}, while $\ell(F)=1$. In
this example $\omega_0=2$ for both maps. For Vasyunin's map $F$ and
$T(X)=(X_1+X_2+X_3+X_4+X_5,\,X_1,\,X_2,\,X_3,\,X_4)$ all five components of
$T^{-1}\circ F\circ T-X$ are nonzero, so $\ell(T^{-1}\circ F\circ T)\ge5>4=\ell(F)$ by
Lemma~\ref{lem:support-bound}, although the two maps have the same degree and the same
dynamical degree.
\item If $\ell(F)\le\#\supp(F)+1$, then $\omega_0(F\times\Id_r)=\omega_0(F)$ for
$r\ge1$. Indeed $\ell(F\times\Id_r)=\ell(F)$ by Proposition~\ref{prop:stab-bound},
and by its proof no factor of a decomposition of $F\times\Id_r$ of length
$\ell(F)$ changes a coordinate $X_k$ with $k>n$. The specialization of that proof gives
a decomposition of $F$ of length at most $\ell(F)$, hence equal to $\ell(F)$, with
factors of degree not larger. Conversely $Q_i\mapsto Q_i\times\Id_r$ preserves the
degrees.
\item For the maps of He and Yan, $\omega_0(F)=2d$ by
Proposition~\ref{prop:heyan-min}, hence $\omega(F)\le2d$. For the
Anick automorphism and $\operatorname{char}K\ne2$, $\omega(A)\le5$ by the word in the
proof of Proposition~\ref{prop:anick}(iv), and nothing is known about $\omega_0(A)$
beyond $\omega_0(A)\ge2$.
\end{enumerate}
\end{remark}

Since $\ell(F\circ G)\le\ell(F)+\ell(G)$, the limit
\[
d\ell(F)=\lim_{m\to\infty}\tfrac1m\,\ell\bigl(F^{(m)}\bigr)
\]
exists for $F\in\EE(K,n)$, as for the dynamical length in the Cremona group
\cite[Definition~1.2.1]{BF}. Recall the dynamical degree $\lambda(F)$ from
Section~\ref{sec:prel}.

\begin{remark}\label{rem:dyn-length}
Let $n=2$ and $d\ell_J(F)=\lim_m\frac1m\ell_J(F^{(m)})$. By
Proposition~\ref{prop:plane-length}(ii), $d\ell_J(F)\le d\ell(F)\le4\,d\ell_J(F)$.
Moreover $d\ell_J(F)>0$ if and only if $\lambda(F)>1$, over every field $K$. Indeed,
let $g_1\cdots g_r$ be a reduced alternating word for $F$. If $r\ge2$ and $g_1,g_r$ lie
in the same factor of the amalgam, conjugation by $g_1$ gives a shorter reduced word.
Repeating this, $F=g\circ G\circ g^{-1}$, where either $G$ lies in $\Aff(K,2)$ or in
$J(K,2)$, or $G$ has a reduced word of length at least $2$ whose first and last letters
lie in different factors. In the first case the iterates of $G$ have bounded degree and
$\ell_J(G^{(m)})\le1$. In the second case the word of $G^{(m)}$ is the $m$-fold
concatenation, which is reduced, so $\ell_J(G^{(m)})=m\,\ell_J(G)$ with $\ell_J(G)\ge1$,
and $\deg G^{(m)}=(\deg G)^m$ with $\deg G\ge2$ \cite[Lemma~9]{AH_canon}. Since
$\ell_J(g\circ G'\circ g^{-1})\le\ell_J(G')+2\ell_J(g)$ by Lemma~\ref{lem:lJ}, and
$\lambda(F)=\lambda(G)$ by Lemma~\ref{lem:bounded}(ii) applied in both directions, the
claim follows. Over $\mathbb{C}$ this dichotomy is the one of Friedland and Milnor
\cite{FM}. Hence $d\ell(F)>0$ on $\EE(K,2)$ if and only if $\lambda(F)>1$. More
precisely, in the second case $d\ell_J(F)=\ell_J(G)=s$ and
$\lambda(F)=\deg G=D_1\cdots D_s$, where $D_1,\dots,D_s\ge2$ are the degrees of the
de~Jonqui\`eres letters of $G$. Conversely, for integers $D_1,\dots,D_s\ge2$ the word
$U_{D_1}\rho\cdots U_{D_s}\rho$, with $U_D=(X_1+X_2^D,X_2)$ and $\rho=(X_2,X_1)$, gives
such a $G$. So the pairs $(d\ell_J(F),\lambda(F))$ with $\lambda(F)>1$ are exactly the
pairs $(s,\lambda)$ of integers with $\lambda\ge2$ and $1\le s\le\Omega(\lambda)$, where $\Omega(\lambda)$ is the number of prime factors of
$\lambda$ counted with multiplicity.
\end{remark}

\begin{remark}\label{rem:complex}
For $n\ge2$ and every field $K$, Lamy defines a simplicial complex $\mathcal C_n$ of
dimension $n-1$ on which $\mathrm{TA}_n(K)$ acts \cite[Section~1.1]{Lam19}. Its
vertices of type $r$ are the classes
$[f_1,\dots,f_r]=\{a\circ(f_1,\dots,f_r):a\in\Aff(K,r)\}$ of $r$-tuples of components of
tame automorphisms, vertices $v_1,\dots,v_n$ of types $1,\dots,n$ span a simplex if
$v_i=[f_1,\dots,f_i]$ for some tame $f$ and all $i$, and $g\in\mathrm{TA}_n(K)$ acts by
$g\cdot[f_1,\dots,f_r]=[f_1\circ g^{-1},\dots,f_r\circ g^{-1}]$. For $n=2$ the graph
$\mathcal C_2$ is the Bass--Serre tree of the amalgam \cite[Remark~1.2]{Lam19}. Let
$v_0=[X_1,\dots,X_n]$, and let $d$ be the combinatorial distance on the $1$-skeleton. Let
$Q$ change $X_j$ by $t$, and let $\pi$ be a permutation of the coordinates which moves
$X_j$ to the last place. The automorphisms $\pi$ and $\pi\circ Q^{-1}$ show that the
vertex $w=[X_i:i\ne j]$ of type $n-1$ is adjacent to $v_0=[\pi]$ and to
$Q\cdot v_0=[Q^{-1}]=[\pi\circ Q^{-1}]$. So $d(v_0,Q\cdot v_0)\le2$, in accordance with
the description of the type $3$ vertices at distance $2$ in \cite{Lam19}. Since
$d(v_0,GG'\cdot v_0)\le d(v_0,G\cdot v_0)+d(v_0,G'\cdot v_0)$, every $F\in\EE(K,n)$
satisfies
\[
d\bigl(v_0,F\cdot v_0\bigr)\le2\,\ell(F),\qquad
\lim_{m\to\infty}\tfrac1m\,d\bigl(v_0,F^{(m)}\cdot v_0\bigr)\le2\,d\ell(F).
\]
For $n=2$ the corresponding statement is the equality $d(v_0,F\cdot v_0)=2\ell_J(F)$ of
Section~\ref{sec:plane}. For $n=3$ and $\operatorname{char}K=0$ the complex $\mathcal C_3$
is Gromov-hyperbolic \cite[Theorem~B]{LP19}, and every $F\in\EE(K,3)$ which acts on it
as a loxodromic isometry satisfies $d\ell(F)>0$.
\end{remark}

\begin{question}\label{q:dyn-length}
Let $n\ge3$ and $F\in\EE(K,n)$ with $\lambda(F)>1$. Is $d\ell(F)>0$? In particular, is
$d\ell(F)>0$ for Vasyunin's map, for which $\lambda(F)=\sqrt2$ by
Theorem~\ref{thmD} and $5\le\ell(F^{(m)})\le4m$ for $m\ge2$ by Proposition~\ref{prop:iterate-left}? Is
$\ell(F^{(m)})$ unbounded? The same
questions arise for the maps of Proposition~\ref{prop:heyan}, for which $\lambda(F)=2d$
and $\ell(F^{(m)})\le4m$.
\end{question}

\appendix

\section{Data for the proofs}\label{app:record}

\subsection{The inverse}\label{ssec:inverse-data}

Expanding \eqref{eq:inverse-compact} gives the inverse $G=F^{-1}$ of Vasyunin's map in
the form of \cite[Section~5.3]{AH_sn}. Here $G_1=X_1$,
\begin{align*}
G_2&=X_2-X_1X_3+\tfrac12X_1X_2^2+X_1^2X_4-X_1^3X_5,\\
G_3&=X_3-X_1X_4-\tfrac12X_2^2+X_1^2X_5,\\
G_4&=X_4-X_1X_5+X_2X_3-X_1X_2X_4-\tfrac12X_1X_3^2-\tfrac12X_2^3
+X_1^2X_2X_5+X_1^2X_3X_4+\tfrac12X_1X_2^2X_3\\
&\quad-X_1^3X_3X_5-\tfrac12X_1^3X_4^2-\tfrac12X_1^2X_2^2X_4-\tfrac18X_1X_2^4
+X_1^4X_4X_5+\tfrac12X_1^3X_2^2X_5-\tfrac12X_1^5X_5^2,\\
G_5&=X_5-\tfrac12G_3^2 .
\end{align*}
The coefficients listed in Section~\ref{ssec:monomials} are read off these formulas.

\subsection{Coefficients used in the proof of Proposition~\ref{prop:ell2}}\label{ssec:monomials}

Let $Q_5$ change $X_c$ by $h$, with $\deg h\le2$, and $P=X_j-(G\circ Q_5)_j$. For the
pairs $(c,j)$ below, the listed monomial has the listed coefficient in $P$ for every
such $h$. It has degree at least $3$, so $P$ is not the tail of an elementary map of
degree at most $2$. The last column gives the term of $G_j$, after expansion of $G_5$,
from which the monomial comes.

\begin{table}[ht]
\centering
\renewcommand{\arraystretch}{1.15}
\begin{tabular}{ccccc}
$c$ & $j$ & monomial & coefficient in $P$ & term of $G_j$\\ \hline
$2$ & $2$ & $X_1X_2^2$ & $-\tfrac12$ & $\tfrac12X_1X_2^2$\\
$2$ & $4$ & $X_2^3$ & $\tfrac12$ & $-\tfrac12X_2^3$\\
$2$ & $5$ & $X_2^4$ & $\tfrac18$ & $-\tfrac18X_2^4$\\
$3$ & $2$ & $X_1^3X_5$ & $1$ & $-X_1^3X_5$\\
$3$ & $3$ & $X_1^2X_5$ & $-1$ & $X_1^2X_5$\\
$3$ & $4$ & $X_1X_3^2$ & $\tfrac12$ & $-\tfrac12X_1X_3^2$\\
$3$ & $5$ & $X_1^2X_3X_5$ & $1$ & $-X_1^2X_3X_5$\\
$4$ & $2$ & $X_1^2X_4$ & $-1$ & $X_1^2X_4$\\
$4$ & $4$ & $X_1X_2X_4$ & $1$ & $-X_1X_2X_4$\\
$4$ & $5$ & $X_1X_3X_4$ & $-1$ & $X_1X_3X_4$\\
$5$ & $2$ & $X_1^2X_4$ & $-1$ & $X_1^2X_4$\\
$5$ & $3$ & $X_1^2X_5$ & $-1$ & $X_1^2X_5$\\
$5$ & $4$ & $X_2^3$ & $\tfrac12$ & $-\tfrac12X_2^3$\\
$5$ & $5$ & $X_1X_3X_4$ & $-1$ & $X_1X_3X_4$
\end{tabular}
\caption{Monomials excluding the pairs $(c,j)$ in Step~2 of the proof of
Proposition~\ref{prop:ell2}.}\label{tab:monomials}
\end{table}

The entries of Table~\ref{tab:monomials} follow from the next lemma.

\begin{lemma}\label{lem:monomial-coeff}
Let $m$ be a monomial, $d$ the exponent of $X_c$ in $m$ and $m'=m/X_c^d$. Let
$h\in K[X_i:i\ne c]$ with $\deg h\le2$. If the coefficients of $m$ in $G_j$ and in
$G_j(X_1,\dots,X_c+h,\dots,X_5)$ differ, then $G_j$ has a term $\mu X_c^e$, with $\mu$
a monomial not involving $X_c$, such that $e>d$, $\mu$ divides $m'$, and
$\deg m'-\deg\mu\le2(e-d)$.
\end{lemma}

\begin{proof}
A term of $G_j$ not involving $X_c$ is unchanged by the substitution. A term $\mu X_c^e$
with $e\ge1$ becomes $\mu(X_c+h)^e=\mu X_c^e+\sum_{i=1}^{e}\binom ei\mu X_c^{e-i}h^i$, and
every monomial of the sum is $\mu X_c^{e-i}\nu$ with $\nu$ a monomial of $h^i$. Since
$\nu$ does not involve $X_c$ and $\deg\nu\le2i$, such a monomial equals $m$ only if
$e-i=d$, $\mu\nu=m'$ and $\deg m'-\deg\mu=\deg\nu\le2i=2(e-d)$.
\end{proof}

For every row of Table~\ref{tab:monomials}, no term of $G_j$ satisfies the three conditions of
Lemma~\ref{lem:monomial-coeff}, as one checks on the formulas of
Section~\ref{ssec:inverse-data} by inspecting the terms of $G_j$ in which $X_c$ occurs
with exponent $e>d$. Hence the coefficient of the listed monomial in $P$ is minus its
coefficient in $G_j$, which is the term in the last column. For instance, for $(c,j)=(3,2)$ the polynomial $G_2$
involves $X_3$ only through $-X_1X_3$, so
$P=X_1X_3+X_1h-\tfrac12X_1X_2^2-X_1^2X_4+X_1^3X_5$, and $X_1h$ contains $X_1^3X_5$ only
if $h$ contains $X_1^2X_5$, which has degree $3$.

\section{Comments on the Anick automorphism}\label{app:anick}

This appendix complements Section~\ref{sec:anick} in two directions. In
Section~\ref{ssec:app-13} the word of the proof of Proposition~\ref{prop:anick}(iv) is
shortened from $16$ to $13$ factors, which gives the bound $\ell(A)\le13$ stated there. In
Section~\ref{ssec:app-base} the theorem of Cohn is used in a form which applies to words
of every length and is insensitive to affine changes of coordinates.

As in the proof of Proposition~\ref{prop:anick}(iv) we use the coordinates
\[
(x,v,z,u)=\Phi(X):=\bigl(X_1,\;\tfrac12(X_2+X_4),\;X_3,\;\tfrac12(X_2-X_4)\bigr)
\]
of \cite{MS}, $\operatorname{char}K\ne2$. Put
\[
q_+=xz-v^2,\qquad q_-=xz+u^2,\qquad \tilde\delta=\delta\circ\Phi^{-1}=xz-v^2+u^2=q_++u^2=q_--v^2 .
\]
Then $\tilde A=\Phi\circ A\circ\Phi^{-1}$ is the map
\begin{equation}\label{eq:app-At}
\tilde A=\bigl(x+\tilde\delta\,(v-u),\;v+\tfrac12\tilde\delta z,\;z,\;u+\tfrac12\tilde\delta z\bigr).
\end{equation}
We write $[w\mathrel{+}=t]$ for the elementary map which changes the coordinate $w$ by the
tail $t$.

\subsection{A word of length thirteen}\label{ssec:app-13}

For $h\in K[u,z,q_+]$ and $k\in K[v,z,q_-]$ let
\[
T_h=\bigl(x+2hv+h^2z,\;v+hz,\;z,\;u\bigr),\qquad R_k=\bigl(x-2ku-k^2z,\;v,\;z,\;u+kz\bigr).
\]
For $h\in K[u,z]$ and $k\in K[v,z]$ these are the maps of the proof of
Proposition~\ref{prop:anick}(iv). Let $\beta=[u\mathrel{+}=\tfrac12q_+]$ and
$\gamma=[v\mathrel{+}=-\tfrac12q_-]$. These are elementary maps, since $q_+$ does not involve
$u$ and $q_-$ does not involve $v$.

\begin{lemma}\label{lem:app-TR}
\begin{enumerate}[label=\textup{(\roman*)},leftmargin=2em]
\item $T_h$ fixes $u$, $z$ and $q_+$, and $R_k$ fixes $v$, $z$ and $q_-$. Consequently
$T_h\circ T_{h'}=T_{h+h'}$ and $R_k\circ R_{k'}=R_{k+k'}$.
\item $\beta^{-1}\circ T_u\circ\beta=T_{u+q_+/2}$ and
$\gamma^{-1}\circ R_{-v}\circ\gamma=R_{q_-/2-v}$.
\item $\tilde A=T_{\tilde\delta/2}\circ R_{\tilde\delta/2}$.
\item For $h\in K[u,z]$ and $k\in K[v,z]$,
\begin{gather*}
T_h=[x\mathrel{+}=2hv-h^2z]\circ[v\mathrel{+}=hz]=[v\mathrel{+}=hz]\circ[x\mathrel{+}=2hv+h^2z],\\
R_k=[x\mathrel{+}=-2ku+k^2z]\circ[u\mathrel{+}=kz].
\end{gather*}
\end{enumerate}
\end{lemma}

\begin{proof}
(i) $(x+2hv+h^2z)z-(v+hz)^2=xz-v^2$ and $(x-2ku-k^2z)z+(u+kz)^2=xz+u^2$. Since $T_{h'}$
fixes $u$, $z$ and $q_+$, it fixes $h$, and
$x+2h'v+h'^2z+2h(v+h'z)+h^2z=x+2(h+h')v+(h+h')^2z$. The same computation applies to $R$.
(ii) The map $\beta$ replaces $u$ by $u+\tfrac12q_+$, then $T_u$ acts as $T_{u+q_+/2}$
and fixes $q_+$, so $\beta^{-1}$ restores $u$. Similarly $\gamma$ replaces $v$ by
$v-\tfrac12q_-$, $R_{-v}$ then acts as $R_{q_-/2-v}$ and fixes $q_-$, and $\gamma^{-1}$
restores $v$. (iii) $R_{\tilde\delta/2}$ fixes $v$ and $q_-$, hence $\tilde\delta$, and
$T_{\tilde\delta/2}$ fixes $u$ and $q_+$, hence $\tilde\delta$. So
$T_{\tilde\delta/2}\circ R_{\tilde\delta/2}$ sends $x$ to
$x-\tilde\delta u-\tfrac14\tilde\delta^2z+\tilde\delta v+\tfrac14\tilde\delta^2z$, $v$ to
$v+\tfrac12\tilde\delta z$ and $u$ to $u+\tfrac12\tilde\delta z$, which is
\eqref{eq:app-At}. (iv) is checked by composition, as in the proof of
Proposition~\ref{prop:anick}(iv).
\end{proof}

Part (ii) is the commutator relation of Smith \cite{Sm}, formed with the coordinate $u$
(respectively $v$), one of the four variables, compare Remark~\ref{rem:anick}.

\begin{proposition}\label{prop:app-13}
Let $\operatorname{char}K\ne2$. Then $A=Q_{13}\circ\dots\circ Q_1$ with the elementary maps
of Table~\ref{tab:app-13}, whose coefficients lie in $\mathbb{Z}[\tfrac12]$. Hence
$7\le\ell(A)\le13$.
\end{proposition}

\begin{table}[ht]
\centering
\small
\renewcommand{\arraystretch}{1.2}
\begin{tabular}{cclc}
factor & changes & tail & degree\\ \hline
$Q_1$ & $X_2$ & $X_4$ & $1$\\
$Q_2$ & $X_4$ & $\tfrac18X_2^2X_3-\tfrac12X_2X_3-\tfrac12X_2$ & $3$\\
$Q_3$ & $X_1$ & $\tfrac1{64}X_2^4X_3-\tfrac18X_2^3X_3+\tfrac14X_2^2X_3-\tfrac14X_2^2X_4+X_2X_4$ & $5$\\
$Q_4$ & $X_2$ & $-X_1X_3-X_4^2$ & $2$\\
$Q_5$ & $X_4$ & $\tfrac12X_2X_3$ & $2$\\
$Q_6$ & $X_1$ & $\tfrac14X_2^2X_3-X_2X_4$ & $3$\\
$Q_7$ & $X_2$ & $X_1X_3+X_3X_4^2+2X_3X_4+X_4^2$ & $3$\\
$Q_8$ & $X_1$ & $-\tfrac14X_3X_4^4-X_3X_4^3-X_3X_4^2+\tfrac12X_2X_4^2+X_2X_4$ & $5$\\
$Q_9$ & $X_4$ & $-\tfrac12X_1X_3+\tfrac18X_2^2$ & $2$\\
$Q_{10}$ & $X_1$ & $X_3X_4^2-X_2X_4$ & $3$\\
$Q_{11}$ & $X_2$ & $-2X_3X_4$ & $2$\\
$Q_{12}$ & $X_4$ & $\tfrac12X_1X_3-\tfrac18X_2^2+\tfrac12X_2$ & $2$\\
$Q_{13}$ & $X_2$ & $-X_4$ & $1$
\end{tabular}
\caption{A decomposition $A=Q_{13}\circ\dots\circ Q_1$ of the Anick automorphism,
$\operatorname{char}K\ne2$. The factor $Q_1$ is applied first.}\label{tab:app-13}
\end{table}

\begin{proof}
Put $h=\tfrac12u^2-u\in K[u,z]$ and $k=v-\tfrac12v^2\in K[v,z]$. Since
$\tfrac12\tilde\delta=(u+\tfrac12q_+)+h=(\tfrac12q_--v)+k$, Lemma~\ref{lem:app-TR}
gives
\begin{equation}\label{eq:app-word}
\tilde A=\beta^{-1}\circ T_u\circ\beta\circ T_h\circ\gamma^{-1}\circ R_{-v}\circ\gamma\circ R_k .
\end{equation}
By Lemma~\ref{lem:app-TR}(iv) this is a word of $12$ elementary maps. We write $T_h$ in
the first form of (iv), so that its first factor $[v\mathrel{+}=hz]$ is applied directly
after $\gamma^{-1}=[v\mathrel{+}=\tfrac12q_-]$. Neither tail involves $v$, so the two
factors merge into $[v\mathrel{+}=hz+\tfrac12q_-]$, and $\tilde A$ is a word $W$ of $11$
elementary maps. Its first factor $[u\mathrel{+}=kz]$ and its last factor
$\beta^{-1}=[u\mathrel{+}=-\tfrac12q_+]$ both change $u$, the fourth coordinate. The
eleven factors are listed, in the order of application, in Table~\ref{tab:app-W}.

Now factor $\Phi$ differently from the proof of Proposition~\ref{prop:anick}(iv). Let
$\Lambda'=[X_4\mathrel{+}=-\tfrac12X_2]\circ[X_2\mathrel{+}=X_4]$ and
$D'=\mathrm{diag}(1,\tfrac12,1,-1)$. Then $\Lambda'(X)=(X_1,X_2+X_4,X_3,\tfrac12(X_4-X_2))$
and $\Phi=D'\circ\Lambda'$, as the matrix identity
\[
\begin{pmatrix}\tfrac12&\tfrac12\\ \tfrac12&-\tfrac12\end{pmatrix}
=\begin{pmatrix}\tfrac12&0\\0&-1\end{pmatrix}
\begin{pmatrix}1&0\\-\tfrac12&1\end{pmatrix}
\begin{pmatrix}1&1\\0&1\end{pmatrix}
\]
on the coordinates $(X_2,X_4)$ shows. Hence
\[
A=\Lambda'^{-1}\circ\bigl(D'^{-1}\circ W\circ D'\bigr)\circ\Lambda',\qquad
\Lambda'^{-1}=[X_2\mathrel{+}=-X_4]\circ[X_4\mathrel{+}=\tfrac12X_2].
\]
Conjugation by $D'$ maps an elementary map changing the $j$-th coordinate to an
elementary map changing the $j$-th coordinate, so $D'^{-1}\circ W\circ D'$ is a word of
$11$ elementary maps whose first and last factors change $X_4$. The first one is applied
directly after $[X_4\mathrel{+}=-\tfrac12X_2]$, and the last one directly before
$[X_4\mathrel{+}=\tfrac12X_2]$. Merging these two pairs gives a word of length
$1+11+1=13$, namely the word of Table~\ref{tab:app-13}. The lower bound is
Theorem~\ref{thm:anick-seven}.
\end{proof}

\begin{table}[ht]
\centering
\small
\renewcommand{\arraystretch}{1.2}
\begin{tabular}{cll}
factor & changes & tail\\ \hline
$1$ & $u$ & $kz$, \ $k=v-\tfrac12v^2$\\
$2$ & $x$ & $-2ku+k^2z$\\
$3$ & $v$ & $-\tfrac12q_-$\\
$4$ & $u$ & $-vz$\\
$5$ & $x$ & $2uv+v^2z$\\
$6$ & $v$ & $hz+\tfrac12q_-$, \ $h=\tfrac12u^2-u$\\
$7$ & $x$ & $2hv-h^2z$\\
$8$ & $u$ & $\tfrac12q_+$\\
$9$ & $x$ & $2uv+u^2z$\\
$10$ & $v$ & $uz$\\
$11$ & $u$ & $-\tfrac12q_+$
\end{tabular}
\caption{The word $W$ of eleven elementary maps for $\tilde A$ in the proof of
Proposition~\ref{prop:app-13}, in the order of application.}\label{tab:app-W}
\end{table}

\begin{remark}\label{rem:app-13}
\begin{enumerate}[label=(\arabic*),leftmargin=2em]
\item The twelve-factor word of the proof of Proposition~\ref{prop:anick}(iv) has no
adjacent factors changing the same coordinate. The gain comes from two freedoms which
do not change $\tilde A$. First, the maps $T_f$, $f\in K[u,z,q_+]$, commute, so $T_h$ may
be placed next to $\gamma^{-1}$ instead of next to $T_u$. Second, $\Phi$ can be factored
so that the elementary maps adjacent to $W$ change the coordinate changed at both ends of
$W$.
\item The direction sequence of Table~\ref{tab:app-13} is
$(2,4,1,2,4,1,2,1,4,1,2,4,2)$. The coordinate $X_3$ is never changed, $X_4$ is changed
four times, in accordance with Proposition~\ref{prop:app-base}(a) below, and the factor
degrees are at most $5$, so the bound $\omega(A)\le5$ of Remark~\ref{rem:degrees}(3) is
unchanged.
\item One has $Q_{13}=Q_1^{-1}$. So $E\circ A\circ E^{-1}$, with
$E=[X_2\mathrel{+}=X_4]$, has a decomposition of length $11$. This is another instance of
Remark~\ref{rem:degrees}(1). We do not know whether $\ell(E\circ A\circ E^{-1})<\ell(A)$.
\item We could not shorten the word further. In \eqref{eq:app-word} the merging of adjacent
factors is only possible between maps from different blocks
$T_{u+q_+/2},T_h,R_{q_-/2-v},R_k$. The blocks of $T$ commute with each other, and so do
the blocks of $R$, but a block of $T$ does not commute with a block of $R$. So the
arrangement allows only one junction between the two groups. The summand $u^2$ of
$\tilde\delta=q_++u^2$ requires the block $T_h$, since $T_{u+q_+/2}$ cannot be replaced by
a single conjugate $\beta'^{-1}\circ T_{g(u)}\circ\beta'$, where $\beta'$ changes $u$
by $f\in K[z,q_+]$ and $g\in K[z][t]$. Indeed $g(u+f)=\tfrac12\tilde\delta$ forces $g$ to
be quadratic with leading coefficient $\tfrac12$, then $f\in K[z]$ by comparison of the
coefficients of $u$, and then the constant term in $u$ lies in $K[z]$, while it should
be $\tfrac12q_+$. By Theorem~\ref{thm:anick-seven}, $7\le\ell(A)\le13$, and the exact value is not
known
(Question~\ref{q:anick}).
\end{enumerate}
\end{remark}

\subsection{Words which do not change the base coordinates}\label{ssec:app-base}

Let $R=K[X_3,X_4]$ and $L=\operatorname{Frac}R$. We regard automorphisms fixing $X_3$ and
$X_4$ as automorphisms of $R[X_1,X_2]$ over $R$. Let
\[
J(R,2)=\bigl\{(aX_1+f(X_2),\,bX_2+c):a,b\in K^*,\ f\in R[t],\ c\in R\bigr\}
\]
be the de~Jonqui\`eres group of $R[X_1,X_2]$ \cite[Chapter~5]{vdE}, since $R^*=K^*$, and let $\mathrm{Aff}_R$ be the group of maps $X'\mapsto MX'+b$ on $X'=(X_1,X_2)$ with
$M\in\mathrm{GE}_2(R)$ and $b\in R^2$. Here $\mathrm{GE}_2(R)$ is the group generated by
the elementary and the invertible diagonal matrices. Let $G_R$ be the group generated by
$J(R,2)$ and $\mathrm{Aff}_R$. It contains every elementary map in $\cE(K,4)$ which changes $X_1$
or $X_2$, since $(X_1,X_2+g(X_1))$ is the conjugate of $(X_1+g(X_2),X_2)\in J(R,2)$ by the
exchange of $X_1$ and $X_2$, which lies in $\mathrm{Aff}_R$.

\begin{lemma}\label{lem:app-cohn}
One has $G_R\cap\Aff(L,2)=\mathrm{Aff}_R$. In particular, if $\varphi\in G_R$ is linear,
$\varphi(X')=MX'$ with $M\in\mathrm{GL}_2(R)$, then $M\in\mathrm{GE}_2(R)$, and the matrix of the Anick automorphism,
\[
M_A=\begin{pmatrix}1+X_3X_4&-X_4^2\\X_3^2&1-X_3X_4\end{pmatrix},
\]
is not the matrix of an element of $G_R$.
\end{lemma}

\begin{proof}
The inclusion $\mathrm{Aff}_R\subseteq G_R\cap\Aff(L,2)$ is clear. Let
$\varphi\in G_R\cap\Aff(L,2)$. Every element of $G_R$ can be written as $\varphi=a_0\circ t_1\circ a_1\circ\dots\circ
t_r\circ a_r$ with $a_i\in\mathrm{Aff}_R$ and $t_i\in J(R,2)$ not affine, since affine
elements of $J(R,2)$ lie in $\mathrm{Aff}_R$. Choose such an expression with $r$ minimal and
suppose $r\ge1$. If some $a_i$, $0<i<r$, is triangular, then its diagonal coefficients
are units of $R$, that is, elements of $K^*$, so $a_i\in J(R,2)$ and
$t_i\circ a_i\circ t_{i+1}\in J(R,2)$. Whether this element is affine or not, the number of
non-affine letters decreases, contrary to the minimality of $r$. So no $a_i$, $0<i<r$, lies in the
triangular subgroup $B(L,2)$, and $a_0\circ t_1\circ\dots\circ t_r\circ a_r$ is a reduced
word in $\Aut_LL[X_1,X_2]^2=\Aff(L,2)\ast_{B(L,2)}J(L,2)$ with $r\ge1$ letters from
$J(L,2)\setminus B(L,2)$. Hence $\varphi\notin\Aff(L,2)$ \cite[Ch.~I, \S1.2]{Ser}, a
contradiction. So $r=0$ and $\varphi=a_0\in\mathrm{Aff}_R$.

For the second statement, $\mathrm{GE}_2(R)=D\cdot\mathrm{E}_2(R)$, where $D$ is the group
of invertible diagonal matrices over $R$, with entries in $K^*$, because $D$ normalizes
$\mathrm{E}_2(R)$. If $M_A=dN$ with $d\in D$ and $N\in\mathrm{E}_2(R)$, then
$\det d=\det M_A=1$, so $d=\mathrm{diag}(c,c^{-1})\in\mathrm{E}_2(K)$, and
$M_A\in\mathrm{E}_2(R)$. This contradicts the theorem of Cohn \cite{Coh66}, as in the proof
of Proposition~\ref{prop:anick-five}.
\end{proof}

For $F\in\Aut_KK[X]^n$ let $\ell_{\mathrm{aff}}(F)$ be the least $k$ such that
$F=\alpha_k\circ Q_k\circ\alpha_{k-1}\circ\dots\circ Q_1\circ\alpha_0$ with affine
automorphisms $\alpha_i$ and $Q_i\in\cE(K,n)$. Unlike $\ell$, the number
$\ell_{\mathrm{aff}}$ is invariant under conjugation by affine maps, and
$\ell_{\mathrm{aff}}\le\ell$.

\begin{proposition}\label{prop:app-base}
Let $A=\alpha_k\circ Q_k\circ\dots\circ Q_1\circ\alpha_0$ with affine $\alpha_i$ and
$Q_i\in\cE(K,4)$ changing $X_{m_i}$. Let $c_i\in K^4$ be the image of the unit vector
$e_{m_i}$ under the linear part of $\alpha_k\circ\dots\circ\alpha_i$. Then the vectors
$c_1,\dots,c_k$ span a subspace of dimension at least $3$. In particular the following holds.
\begin{enumerate}[label=\textup{(\alph*)},leftmargin=2em]
\item Every word in $\cE(K,4)$ representing $A$ changes $X_3$ or $X_4$, and one of them
at least twice.
\item One has $\ell_{\mathrm{aff}}(A)\ge3$ for every field $K$, and $\ell_{\mathrm{aff}}(A)\le11$
if $\operatorname{char}K\ne2$.
\end{enumerate}
\end{proposition}

\begin{proof}
Put $\beta_i=\alpha_k\circ\dots\circ\alpha_i$ and $Q_i'=\beta_i\circ Q_i\circ\beta_i^{-1}$.
Then $A=Q_k'\circ\dots\circ Q_1'\circ\alpha$ with $\alpha=\beta_0$, and
$Q_i'=X+c_i\,g_i$ with $g_i\in K[X]$ and $\partial_{c_i}g_i=0$. Let $W$ be the span of the
$c_i$. Since $(X+c\,g)\circ(X+w)=X+w+c\,g(X+w)$, the map $G=A\circ\alpha^{-1}=Q_k'\circ\dots\circ
Q_1'$ satisfies $G-X\in W\otimes K[X]$. On the other hand $A=X+H$ with
$H=(\delta X_4,\delta X_3,0,0)$ homogeneous of degree $3$, so the homogeneous component of
degree $3$ of $G-X=\alpha^{-1}-X+H\circ\alpha^{-1}$ is
$e_1\,(\delta X_4)\circ N+e_2\,(\delta X_3)\circ N$, where $N$ is the linear part of
$\alpha^{-1}$. The two cubic forms are linearly independent, hence $e_1,e_2\in W$.

Suppose $\dim W\le2$, so $W=Ke_1+Ke_2$. Then every $Q_i'$ fixes $X_3$ and $X_4$, hence so
do $G$ and $\alpha=G^{-1}\circ A$. For each $i$ choose $N_i\in\mathrm{GL}_2(K)$, acting
on $(X_1,X_2)$, with $N_ie_1=c_i$. Then $N_i^{-1}\circ Q_i'\circ N_i$ changes $X_1$ by
$g_i\circ N_i$, which does not involve $X_1$, so it lies in $J(R,2)$, and $Q_i'\in G_R$. The
affine map $\alpha$ fixes $X_3$ and $X_4$, so it acts on $(X_1,X_2)$ by a matrix in
$\mathrm{GL}_2(K)\subseteq\mathrm{GE}_2(R)$ and a translation by an element of $R^2$, so
$\alpha\in\mathrm{Aff}_R$. Hence $A\in G_R$, which contradicts
Lemma~\ref{lem:app-cohn}. So $\dim W\ge3$.

(a) If all $\alpha_i$ are the identity, then $c_i=e_{m_i}$, so some $m_i$ lies in
$\{3,4\}$. Since $A$ fixes $X_3$ and $X_4$, the coordinate $X_{m_i}$ is changed at least twice by
Lemma~\ref{lem:word}(i). (b) The lower bound is clear, since $\dim W\le k$.
For the upper bound, $A=\Phi^{-1}\circ\tilde A\circ\Phi$ with $\Phi$ linear, and
$\tilde A$ is the word $W$ of $11$ elementary maps of the proof of
Proposition~\ref{prop:app-13}.
\end{proof}

The bound $3$ on the dimension is attained, since the factors of the word $W$ of
Table~\ref{tab:app-W} change only the coordinates $x$, $v$, $u$ of
Section~\ref{ssec:app-13}.

Part (a) contains the alternating cases $(1,2,1,2)$, $(2,1,2,1)$, $(1,2,1,2,1)$ and
$(2,1,2,1,2)$ of Propositions~\ref{prop:anick-five} and~\ref{prop:anick-six}, for words of
every length and with no degree comparison. So in every decomposition of $A$ the base
coordinates must be used, as in the Smith commutator of Proposition~\ref{prop:anick}(iii),
where the additional coordinate $X_5$ is changed twice, and in
Table~\ref{tab:app-13}, where $X_4$ is changed four times. The bound $3\le\ell_{\mathrm{aff}}(A)\le11$ leaves a large gap. The argument above only
uses the directions of the factors, and it gives no information on words in which three
independent directions occur.

\section*{Acknowledgements}
This research was supported by the AGH University of Krakow within a subsidy of the Polish Ministry of Science and Higher Education (grant no. 16.16.420.054).

\end{document}